\documentclass{amsart}
\usepackage{amssymb, amsmath, latexsym, amscd, graphicx, color}
\usepackage[all]{xy}
\usepackage[dvipdfm]{hyperref}
\usepackage[all]{hypcap}
\newtheorem{theorem}{Theorem}[section]
\newtheorem{lemma}[theorem]{Lemma}
\newtheorem{cor}[theorem]{Corollary}
\newtheorem{definition}[theorem]{Definition}
\newtheorem{example}[theorem]{Example}
\newtheorem{remark}[theorem]{Remark}
\newtheorem{proposition}[theorem]{Proposition}
\def\pagenumber{1}
\begin{document}
\setcounter{page}{\pagenumber}
\newcommand{\T}{\mathbb{T}}
\newcommand{\R}{\mathbb{R}}
\newcommand{\Q}{\mathbb{Q}}
\newcommand{\N}{\mathbb{N}}
\newcommand{\Z}{\mathbb{Z}}
\newcommand{\tx}[1]{\quad\mbox{#1}\quad}
\parindent=0pt
\def\SRA{\hskip 2pt\hbox{$\joinrel\mathrel\circ\joinrel\to$}}
\def\tbox{\hskip 1pt\frame{\vbox{\vbox{\hbox{\boldmath$\scriptstyle\times$}}}}\hskip 2pt}
\def\circvert{\vbox{\hbox to 8.9pt{$\mid$\hskip -3.6pt $\circ$}}}
\def\IM{\hbox{\rm im}\hskip 2pt}
\def\ES{\vbox{\hbox to 8.9pt{$\big/$\hskip -7.6pt $\bigcirc$\hfil}}}
\def\TR{\hbox{\rm tr}\hskip 2pt}
\def\TRS{\hbox{\rsmall tr}\hskip 2pt}
\def\GRAD{\hbox{\rm grad}\hskip 2pt}
\def\GRADS{\hbox{\rsmall grad}\hskip 2pt}
\def\bull{\vrule height .9ex width .8ex depth -.1ex}
\def\VLLA{\hbox to 25pt{\leftarrowfill}}
\def\VLRA{\hbox to 25pt{\rightarrowfill}}
\def\DU{\mathop{\bigcup}\limits^{.}}
\setbox2=\hbox to 25pt{\rightarrowfill}
\def\DRA{\vcenter{\copy2\nointerlineskip\copy2}}
\def\RANK{\hbox{\rm rank}\hskip 2pt}
\font\rsmall=cmr7 at 7truept \font\bsmall=cmbx7 at 7truept
\newfam\Slfam
\font\tenSl=cmti10 \font\neinSl=cmti9 \font\eightSl=cmti8
\font\sevenSl=cmti7 \textfont\Slfam=\tenSl
\scriptfont\Slfam=\eightSl \scriptscriptfont\Slfam=\sevenSl
\def\Sl{\fam\Slfam\tenSl}

\vskip -1cm
\title[Extended Crystal PDE's]{\mbox{}\\[1cm] EXTENDED CRYSTAL PDE'S}
\author{Agostino Pr\'astaro}
\maketitle
\vspace{-.5cm}
{\footnotesize
\begin{center}
Department of Methods and Mathematical Models for Applied
Sciences, University of Rome ''La Sapienza'', Via A.Scarpa 16,
00161 Rome, Italy. \\
E-mail: {\tt Prastaro@dmmm.uniroma1.it}
\end{center}
\vspace{.5cm} {\bsmall ABSTRACT.} In this paper we show that between
PDE's and crystallographic groups there is an unforeseen relation.
In fact we prove that integral bordism groups of PDE's can be
considered extensions of crystallographic subgroups. In this respect
we can consider PDE's as {\em extended crystals}. Then an
algebraic-topological obstruction ({\em crystal obstruction}),
characterizing existence of global smooth solutions for smooth
boundary value problems, is obtained. Applications of this new theory to the Ricci-flow equation and Navier-Stokes equation are given that solve some well-known fundamental problems. These results, are also
extended to singular PDE's, introducing ({\em extended crystal
singular PDE's}). An application to singular MHD-PDE's, is given
following some our previous results on such equations, and showing
existence of (finite times stable smooth) global solutions crossing
critical nuclear energy production zone.\footnote{See also \cite{PRA30, PRA31, PRA32, PRA33, PRA34, PRA35}.\\Work partially supported by Italian grants MURST ''PDE's Geometry and Applications''.}

\vskip 0.5cm

{\bsmall AMS (MOS) MS CLASSIFICATION. 55N22, 58J32, 57R20; 20H15.}
{\rsmall KEY WORDS AND PHRASES. Integral bordisms in PDE's;
Existence of local and global solutions in PDE's; Conservation laws;
Crystallographic groups; singular PDE's; singular MHD-PDE's.}}

\section{\bf Introduction}
New points of view were recently introduced by us in the geometric
theory of PDE's, by adopting some algebraic topological approaches.
In particular, integral (co)bordism groups are seen very useful to
characterize global solutions. The methods developed by us, in the
category of (non)commutative PDE's, in order to find integral
bordism groups, allowed us to obtain, as a by-product, existence
theorems for global solutions, in a pure geometric way. Another
result that is directly related to the knowledge of PDE's integral
bordism groups, is the possibility to characterize PDE's by means of
some important algebras, related to the conservation laws of these
equations (PDE's Hopf algebras). Moreover, thanks to an algebraic
characterization of PDE's, one has also a natural way to recognize
quantized PDE's as quantum PDE's, i.e., PDE's in the category of
quantum (super)manifolds, in the sense introduced by us in some
previous works. These results have opened a new sector in Algebraic
Topology, that we can formally define the {\em PDE's Algebraic
Topology}. (See Refs.\cite{PRA1, PRA2, PRA3, PRA4, PRA5, PRA6, PRA7,
PRA8, PRA9, PRA10, PRA11, PRA12, PRA13, PRA14, PRA15, PRA16, PRA17,
PRA18, PRA19, PRA20, PRA21, PRA22, PRA23, PRA24, PRA25, PRA26,
PRA27, PRA28, PRA29}, and related works \cite{AG-PRA1, AG-PRA2, L-P, PRA-RAS1, PRA-RAS2, PRA-RAS3}.)

Aim of the present paper is, now, to show that PDE's can be
considered as extended crystals, in the sense that their integral
bordism groups, characterizing the geometrical structure of PDE's,
can be considered as extended groups of crystallographic subgroups.
This fundamental relation gives new insights in the PDE's geometrical structure understanding, and opens also new possible
mathematical and physical interpretations of the same PDE's
structure. In particular, we get a new general workable criterion
for smooth global solutions existence satisfying smooth boundary
value problems. In fact, we identify the obstruction for existence
of such global solutions, with an algebraic-topologic object ({\em
crystal obstruction}). Since it is easy to handle this method in all
the concrete PDE's of interest, it sheds a lot of light on all the
PDE's theory.

Finally we extend above results also to singular PDE's, and we
recognize {\em extended crystal singular PDE's}. For such equations
we identify algebraic-topological obstructions to the existence of
global (smooth) solutions solving boundary value problems and
crossing singular points too. Applications to MHD-PDE's, as
introduced in some our previous papers \cite{PRA24, PRA26}, and
encoding anisotropic incompressible nuclear plasmas dynamics are
given.

The paper, after this Introduction, has three more sections,  and
four appendices. In Section 2 we consider some fundamental
mathematical properties of crystallographic groups that will be used
in Section 3. There we recall some our results about PDE's
characterization by means of integral bordism groups. Furthermore we
relate these groups to crystallographic groups. The main results are
Theorem \ref{main1}, Theorem \ref{main2} and Theorem \ref{main3}.
The first two relate formal integrability and complete integrability
of PDE's to crystallographic groups. (It is just this theorem that
allows us to consider PDE's as {\em extended crystallographic
structures}.) The third main theorem identifies an obstruction,
({\em crystal obstruction}), characterizing existence of global
smooth solutions in PDE's. Applications to some under focus PDE's
of Riemannian geometry (e.g., Ricci-flow equation) and Mathematical Physics (e.g., Navier-Stokes equation) are given that solve some well-known fundamental mathematical problems. (Further applications are given in Refs.\cite{PRA23, PRA24, PRA25,
PRA26}.) Section 4 is devoted to extend above results also to
singular PDE's. The main result in this section is Theorem
\ref{main-singular1} that identifies conditions in order to
recognize global (smooth) solutions of singular PDE's crossing
singular points. There we characterize {\em $0$-crystal singular
PDE's}, i.e., singular PDE's having smooth global solutions crossing
singular points, stable at finite times. Applications of these
results to singular MHD-PDE's, encoding anisotropic incompressible
nuclear plasmas dynamics are given in Example \ref{MHD-PDE}. Here,
by using some our previous recent results on MHD-PDE's, we
characterize global (smooth) solutions crossing critical nuclear
zone, i.e., where solutions pass from states without nuclear energy
production, to states where there is nuclear energy production. The
stability of such solutions is also considered on the ground of our
recent geometric theory on PDE's stability and stability of global
solutions of PDE's \cite{PRA22, PRA23, PRA24, PRA25, PRA26}.

In appendices are collected and organized some standard informations
about crystallographic groups and their subgroups in order to give
to the reader a general map for a more easy understanding of their
using in the examples considered in the paper.

\section{\bf CRYSTALLOGRAPHIC GROUPS}

\begin{definition}
Let $(E,({\bf E},\bar g),\alpha)$ be a $d$-dimensional Euclidean
affine space, where $\alpha:{\bf E}\times E\to E$ is the action
mapping of the $n$-dimensional vector space ${\bf E}$ on the set $E$
of points, and $\bar g$ is an Euclidean metric on ${\bf E}$. Let us
denote by $A(E)={\bf E}\rtimes GL({\bf E})$ the {\em affine group}
of $E$, i.e., the symmetry group of the above affine structure and
by $M(E)= {\bf E}\rtimes O({\bf E})$ the {\em group of Euclidean
motions} of $E$, i.e., the symmetry group of the above euclidean
affine structure. (The symbol $\rtimes$ denotes {\em semidirect
product}, i.e., the set is the cartesian product and the
multiplication is defined as
$(a,u)(b,v)=(ab,u+av)$.)\footnote{Recall that given two groups $A$
and $B$ and an homomorphism $\alpha: B\to Aut(A)$, the {\em
semidirect product} is a group denoted by $A\rtimes_\alpha B$, that
is the cartesian product $A\times B$, with product given by
$(a,b).(\bar a,\bar b)=(a.\alpha(b)(\bar a),b.\bar b)$. The
semidirect product reduces to the direct product, i.e.,
$A\rtimes_\alpha B=A\times B\equiv A\oplus B$, when $\alpha(b)=1$,
for all $b\in B$. In the following we will omit the symbol
$\alpha$.} Let us denote by $R(E)={\bf E}\rtimes SO({\bf E})$ the
group of all rigid motions of $E$, i.e., the symmetry group of the
above oriented euclidean affine structure, where the orientation is
the canonical one induced by the metric. One has the following
monomorphisms of inclusions: $R(E)< M(E)<A(E)$. One has the natural
split exact commutative diagram
\begin{equation}\label{comdiag1}
\xymatrix{0\ar[r]&{\bf E}\ar@{=}[d]\ar[r]&A(E)\ar@<0.3ex>[r]&GL({\bf
E})\ar@<0.3ex>[l]\ar[r]&1\\
0\ar[r]&{\bf
E}\ar@{=}[d]\ar[r]&M(E)\ar@{^{(}->}[u]\ar@<0.3ex>[r]&O({\bf
E})\ar[r]\ar@{^{(}->}[u]\ar@<0.3ex>[l]&1\\
0\ar[r]&{\bf E}\ar[r]&R(E)\ar@{^{(}->}[u]\ar@<0.3ex>[r]&SO({\bf
E})\ar[r]\ar@{^{(}->}[u]\ar@<0.3ex>[l]&1\\}
\end{equation}
The quotient groups $A(E)/{\bf E}\cong GL({\bf E})$, $M(E)/{\bf
E}\cong O({\bf E})$ and $R(E)/{\bf E}\cong SO({\bf E})$, are called
{\em point groups} of the corresponding groups $A(E)$, $M(E)$ and
$R(E)$ respectively. ${\bf E}$ is called the {\em translations
group}.
\end{definition}

\begin{definition} A {\em crystallographic group} is a cocompact
\footnote{A discrete subgroup $H\subset G$ of a topological group
$G$, is {\em cocompact} if there is a compact subset $K\subset G$
such that $HK=G$.}, discrete subgroup of the isometries of some
Euclidean space.
\end{definition}

\begin{definition}
A $d$-dimensional {\em affine crystallographic group} $G(d)$ is a
subgroup of $M(E)$, such that its subgroup ${\bf T}\equiv
G(d)\cap{\bf E}$ of all pure translations is a discrete normal
subgroup of finite index.

If the rank of ${\bf T}$ is $d$, i.e., ${\bf T}\cong
\mathbb{Z}^d$, $G(d)$ is called a {\em space group}. The point
group $G\equiv G(d)/{\bf T}$ of a space group is finite, and
isomorphic to a subgroup of $O(\mathbf{E})$: $G<O(\mathbf{E})$.
\end{definition}

\begin{remark}
Note that the point group $G$ of a crystallographic group $G(d)$
does not necessarily can be identified with a subgroup of $G(d)$. In
other words, $G(d)$ is in general an {\em upward extension} of ${\bf
T}$, as well as an {\em downward extension} of $G$, i.e., for a
crystallographic group $G(d)$ one has the following short exact
sequence:\footnote{Usually one denotes such an extension simply with
$G(d)/\mathbf{T}$, whether we are not interested to emphasize the
notation for $G\cong G(d)/\mathbf{T}$.}
\begin{equation}\label{exactseq1}
\xymatrix{0\ar[r]&{\bf T}\ar[r]^{i}&G(d)\ar[r]^{\pi}&G\ar[r]&1\\}
\end{equation}
\end{remark}

Thus we can give the following definition.

\begin{definition}\label{symcrygr}
We call {\em symmorphic crystallographic groups} ones $G(d)$ such
that the exact sequence {\em(\ref{exactseq1})} splits.\end{definition}

\begin{theorem}[Characterization of symmorphic crystallographic
groups] The following propositions are equivalent.

{\em(i)} $G(d)$ is a $d$-dimensional symmorphic crystallographic
group.

{\em(ii)} $G(d)$ has a subgroup $\widetilde{G}< G(d)$ mapped by
$\pi$ isomorphically onto $G$, i.e., $G(d)=i(\mathbf{T}).G$ and
$i(\mathbf{T})\cap\widetilde{G}=\{1\}$.

{\em(iii)} $G(d)$ has a subgroup $\widetilde{G}< G(d)$ such that
every element $a\in G(d)$ is uniquely expressible in the form
$a=i(h)\tilde g$, $h\in\mathbf{T}$, $\tilde g\in\widetilde{G}$.

{\em(iv)} The short exact sequence {\em(2)} is equivalent to the
extension

\begin{equation}\label{exactseq2}
\xymatrix{0\ar[r]&{\bf T}\ar[r]^{i'}&\mathbf{T}\rtimes
G\ar[r]^{\pi{}'}&G\ar[r]&1\\}
\end{equation}
\end{theorem}

\begin{proof}
The equivalence of the propositions (i)--(iv) follows from the
Definition \ref{symcrygr} and standard results of algebra. (See, e.g.,
\cite{BOURB}.)
\end{proof}

\begin{theorem}[Cohomology symmorphic crystallographic
groups classes] The symmorphic crystallographic groups $G(d)$,
having point group $G$ and translations group $\mathbf{T}$, are
classified by $\mathbf{T}$-conjugacy classes that are in 1-1
correspondence with the elements of $H^1(G;\mathbf{T})$.
\end{theorem}

\begin{proof}
Even if this result refers to standard subjects in homological
algebra, let us enter in some details of the proof, in order to
better specify how the theorem works. In fact these details will be
useful in the following. In the case $G$ acts trivially on
$\mathbf{T}$, so that the group $ G(d)=\mathbf{T}\times G$, the
splitting of (\ref{exactseq1}) are in 1-1 correspondence with
homomorphisms $G\to\mathbf{T}$. In general case the splitting
correspond to derivations ({\em crossed homomorphisms})
$d:G\to\mathbf{T}$ satisfying $d(ab)=da+a.db$, for all $a,b\in G$.
In fact, let us consider the extension (\ref{exactseq2}). A section
$s:G\to \mathbf{T}\rtimes G$ has the form $s(g)=(dg,g)$, where $d$
is a function $G\to\mathbf{T}$. One has $s(g)s(g')=(dg+g.dg',gg')$.
So $s$ will be a homomorphism iff $d$ is a derivation. Two splitting
$s_1$, $s_2$, are called {\em$\mathbf{T}$-conjugate} if there is an
element $h\in\mathbf{T}$ such that $s_1(g)=i(h)s_2i(h)^{-1}$, for
all $g\in G$. Since $(h,1)(k,g)(h,1)^{-1}=(h+k-gh,g)$ in $\mathbf{T}
\rtimes G$, the conjugacy relation becomes $d_1g=h+d_2g-gh$, in
terms of the derivations $d_1$, $d_2$ corresponding to $s_1$ and
$s_2$ respectively. Thus $d_1$ and $d_2$ correspond to $\mathbf{T}
$-conjugate splittings iff their difference $d_2-d_1$ is a function
({\em principal derivation}) $G\to\mathbf{T}$ of the form $g\mapsto
gh-a$ for some fixed $h\in\mathbf{T}$. Therefore
$\mathbf{T}$-conjugacy classes of splittings of a split extension of
$G$ by $\mathbf{T}$ correspond to the elements of the quotient group
$Der(G,\mathbf{T})/P(G,\mathbf{T})$, where $Der(G,\mathbf{T})$ is
the abelian group of derivations $G\to \mathbf{T}$, and
$P(G,\mathbf{T})$ is the group of principal derivations. On the
other hand, considering the cochain complex
$C^\bullet(G,\mathbf{T})$, we see that $Der(G,\mathbf{T})$ is the
group of $1$-cocycles and $P(G,\mathbf{T})$ is the group of
$1$-coboundaries. Thus we get the theorem.
\end{proof}

\begin{theorem}[Cohomology crystallographic group classes]
The cohomological classification of $d$-dimensional crystallographic
groups $G(d)$, with point group $G$, is made by means of the first
cohomology group $H^1(G;\mathbb{R}^d/\mathbb{Z}^d)$. One has the
natural isomorphism: $H^1(G;\mathbb{R}^d/\mathbb{Z}^d)\cong
H^2(G;\mathbb{Z}^d)$. Two cohomology classes define equivalent
crystallographic groups iff they are transformed one another by the
normalizer of $G$ in $GL_d(\mathbb{Z})$. Two crystallographic groups
$G(d)$, $G(d)'$, belong to the same class {\em(arithmetical class)}
if their point groups, respectively $G$, $G'$, are conjugate in
$GL_d(\mathbb{R})$, (in $GL_d(\mathbb{Z})$).
\end{theorem}

\begin{proof}
The proof follows directly by the following standard lemmas of
(co)homological algebra.
\begin{lemma}[\cite{BROW}]\label{lem1}
Let $G$ be a group and $A$ a $G$-module. Let $\mathcal{E}(G,A)$ be
the set of equivalence classes of extensions of $G$ by $A$
corresponding to the fixed action of $G$ on $A$. Then, there is a
bijection $\mathcal{E}(G,A)\cong H^2(G,A)$.
\end{lemma}

\begin{lemma}[\cite{BROW}]\label{lem2}
For any exact sequence
\begin{equation}\xymatrix{0\ar[r]&M'\ar[r]&M\ar[r]&M''\ar[r]&0\\}
\end{equation}
of $G$-modules and any integer $n$ there is a natural map
$\delta:H^n(G;M'')\to H^{n+1}(G;M')$ such that the sequence
\begin{equation}
\xymatrix{0\ar[r]&H^0(G;M')\ar[r]&H^0(G;M)\ar[r]&H^0(G;M'')\ar[r]^{\delta}&H^1(G;M')\ar[r]
&\dots\\}
\end{equation}
is exact. Furthermore if $P$ is a projective (resp. if $Q$ is an
injective) $\mathbb{Z}G$-module then $H_n(G;P)=0$ (resp.
$H^n(G;Q)=0$) for $n>0$.\footnote{$\mathbb{Z}G$ is the free
$\mathbb{Z}$-module generated by the elements of $G$. The
multiplication in $G$ extends uniquely to a $\mathbb{Z}$-bilinear
product $\mathbb{Z}G\times\mathbb{Z}G\to \mathbb{Z}G$ so that
$\mathbb{Z}G$ becomes a ring called the {\em integral group ring} of
$G$. A {\em$G$-module} $A$, is just a (left)$\mathbb{Z}G$-module on
the abelian group $A$. This action can be identified by a group
homomorphism $G\to Hom_{group}(A)$.}
\end{lemma}
In fact, it is enough to take $A=\mathbf{T}=\mathbb{Z}^d$,
$M'=\mathbb{Z}^d$, $M=\mathbb{R}^d$ and
$M''=\mathbb{R}^d/\mathbb{Z}^d$.
\end{proof}

\begin{definition}
We call {\em lattice subgroups} of a crystallographic group $G(d)$
the set of its subgroups.\footnote{The lattice subgroups of a group
$G$ is a lattice under inclusion. The identity $1\in G$ identifies
the {\em minimum} $\{1\}<G$ and the {\em maximum} is just $G$. In
the following we will denote also by $e=1$ the identity of a group
$G$. The subgroup $\mathbf{T}\subset G(d)$ is called also the {\em
Bieberbach lattice} of $G(d)$.}
\end{definition}

\begin{definition}
We call {\em $d$-dimensional crystal} any topological space
${}^{\maltese}E\subset E$, i.e., contained in the $d$-dimensional
Euclidean affine space $E$, having the crystallographic group $G(d)$
as (extrinsic) symmetry group. We call {\em unit cell} of
${}^{\maltese}E$, the compact quotient space ${}^{\maltese}E/{\bf
T}$, having the point group as (extrinsic) symmetry
group.\footnote{For example if ${}^{\maltese}E$ is identified by
means of an infinite graph, (see, e.g., \cite{SUN}), then the unit
cell is the corresponding fundamental finite graph. Of course
$d$-dimensional chains can be associated to such graphs, so that
also the corresponding unit cells can be identified with compact
$d$-chains.}
\end{definition}

\begin{theorem}[Hilbert's 18th problem] For any dimension $d$, there are, up to
equivalence, only finitely many-dimensional crystallographic space
groups. These are called {\em space-group types}. We denote by
$[G(d)]$ the space-group type identified by the space group
$G(d)$.\footnote{Note that a space group is characterized other than
by translational and point symmetry, also by the metric-parameters
characterizing the unit cell. Thus the number of space groups is
necessarily infinite.}
\end{theorem}

\begin{proof}
The 18th Hilbert's problem put at the beginning of 1900, has been
first proved in 1911 by L. Bieberbach \cite{BIEB}. For modern proofs
see also L. S. Charlap. \cite{CHARL}.\footnote{See also
Refs.\cite{ABELS, ALP-BEL, ARFK, BERN-SCHW, K-N, P-P, RAGH, RUH,
SCHW}, and works by M. Gromov \cite{GRO} and E.Ruh \cite{RUH} on the
almost flat manifolds, that are related to such crystallographic
groups.}
\end{proof}

\begin{example}[$3$-dimensional crystal-group types]\label{ex1}
Up to isomorphisms, there are 17 crystallographic groups in
dimension $2$ and $219$ in dimension $3$. However, if the spatial
groups are considered up to conjugacy with respect to
orientation-preserving affine transformations, their number is 230.
These last can be called {\em affine space-group types}.
 These are just the F.S. Fedorov and A. Sch\"onflies groups {\em\cite{FED, SCH}}. The
 corresponding translation-subgroup types, or {\em Bravais lattices}, are 14.
 In Appendix A, Tab. \ref{the-32-crystallographic-point-groups-of-the-space-group-types}, are reported the 32 affine crystallographic point-group types, and in
 Tab. \ref{230-affine-crystallographic-space-group-types-3-dim} the 230 affine crystallographic space-group types. There are 66 symmorphic space-group
 types in $[G(3)]$. For the other 164 $G(3)$ cannot be identified with the semidirect
product ${\bf T}\rtimes G(3)/{\bf T}$. It is useful, also, to know
the subgroups corresponding to the 32 crystallographic point groups,
i.e., their subgroup-lattices, in relation to the results of the
following section where integral bordism groups of PDE's will be
related to crystallographic subgroups. For this we report in
Appendix B a full list for such subgroups, and in Appendix C a list
of amalgamated free products in $[G(3)]$.\footnote{Let $A$, $B$ and
$C$ groups and $B\lhd A$, $B\lhd C$, the {\em amalgamated free
product} $A\bigstar_BC$ is generated by the elements of $A$ and $C$
with the common elements from $B$ identified.} See also in Appendix
A, Tab. \ref{the-17-affine-crystallographic-space-group-types-2-dim}
 for the $2$-dimensional 17 crystal-group types, and Appendix D
 for the corresponding list of subgroups.\footnote{Let us recall that the index of a
subgroup $H\subset G$, denoted $[G:H]$, is the number of {\em left
cosets} $aH=\{ah : h\in H\}$, (resp. {\em right cosets} $Ha=\{ha :
h\in H\}$), of $H$. For a finite group $G$ one has the following
formula: ({\em Lagrange's formula}) $[G:H]=\frac{o(G)}{o(H)}$, where
$o(G)$, resp. $o(H)$, is the order of $G$, resp. $H$. If $aH=Ha$,
for any $a\in G$, then $H$ is said to be a {\em normal subgroup}.
Every subgroup of index $2$ is normal, and its cosets are the
subgroup and its complement.}
\end{example}

\begin{definition}[Generalized crystallographic group]
 Let $K$ be a principal ideal domain, $H$ a finite group, and $M$ a
 $KH$-module which as $K$-module is free of finite rank, and on
 which $H$ acts faithfully. A {\em generalized crystallographic group} is a group
 $G(d)$ which has a normal subgroup isomorphic to $M$ such that the following conditions are satisfied:

{\em(i)} $G(d)/M\cong  H$.

{\em(ii)} Conjugation in $G(d)$ gives the same action of $H$ on
$M$.

{\em(iii)} The extension
\begin{equation}\xymatrix{0\ar[r]&M\ar[r]&G(d)\ar[r]&H\ar[r]&0\\}
\end{equation}
does not split. We define {\em dimension} of $G(d)$ the $K$-free
rank $d$ of $M$. We define {\em holonomy group} of $G(d)$ the group
$H$.
\end{definition}

\begin{example}\label{ex2}
With $K=\mathbb{Z}$ one has the crystallographic groups.
\end{example}

\begin{example}\label{ex3}
A {\em complex crystallographic group} is a discrete group of affine
transformations of a complex affine space $(V,{\bf V})$, such that
the quotient $X\equiv V/G$ is compact. This is a generalized
crystallographic group with $K=\mathbb{C}$.\footnote{See also the
recent work by Bernstein and Schwarzman on the complex
crystallographic groups \cite{BERN-SCHW}. }
\end{example}

In the following we give some examples of crystallographic subgroups
in dimension $d=3$. In fact, it will be useful to know such
subgroups in relation to results of the next section.

\begin{example}[The group $\mathbb{Z}\bigoplus\mathbb{Z}$]\label{ex4}
This is the crystallographic subgroup, of the crystallographic group
$G(3)=\mathbb{Z}^3\times\mathbb{Z}_1$, with point group
$C_1=\mathbb{Z}_1$ (triclinic syngony). Let us emphasize that
$\mathbb{Z}\bigoplus\mathbb{Z}$ is crystallographic since it can be
identified with the $2$-dimensional crystallographic group
$G(2)=\mathbb{Z}^2\times\mathbb{Z}_1=p1$, generated by translations
parallel to the $x$ and $y$-axes, with point group $\mathbb{Z}_1$.
\end{example}

\begin{example}[The groups $\mathbb{Z}^2\bigoplus\mathbb{Z}_n$, $n=2,3,4,6$]\label{ex5}
These are subgroups of the crystallographic groups
$G(3)=\mathbb{Z}^3\rtimes\mathbb{Z}_n$, $n=2,3,4,6$, (point groups
$\mathbb{Z}_n$, monoclinic, hexagonal, tetragonal, trigonal syngony
respectively). Furthermore, $\mathbb{Z}^2\times\mathbb{Z}_n$,
$n=2,3,4,6$, can be considered also subgroups of the $2$-dimensional
crystallographic groups $G(2)=\mathbb{Z}^2\rtimes\mathbb{Z}_n=pn$,
$n=2,3,4,6$, (point groups $\mathbb{Z}_n$, oblique, trigonal,
square, hexagonal syngony respectively), generated by translations
parallel to the $x$ and $y$-axes, and a rotation by $\pi/n$,
$n=2,3,4,6$, about the origin.
\end{example}

\begin{example}[The groups $\mathbb{Z}_2\bigoplus\mathbb{Z}_2$]\label{ex6}
This group coincides with the amalgamated free product
$\mathbb{Z}_2\bigstar_e\mathbb{Z}_2$ that is generated by reflection
over $x$-axis and reflection over $y$-axis. It is a subgroup of the
crystallographic group $G(3)=\mathbb{Z}^3\rtimes D_4$, (point group
$D_4$, tetragonal syngony). (See Appendix A and Appendix C.) The
group $\mathbb{Z}_2\bigoplus\mathbb{Z}_2$ is  also a subgroup of the
$2$-dimensional crystallographic groups $G(2)=\mathbb{Z}^2\rtimes
D_4=p4m$ and $G(2)=\mathbb{Z}^2\rtimes D_4=p4g$ that have both point
group $D_4=\mathbb{Z}_2\times\mathbb{Z}_2$, and square syngony. (See
in Appendix A, Tab. \ref{diehdral-groups} and Tab. \ref{the-17-affine-crystallographic-space-group-types-2-dim}.)
\end{example}

\begin{example}[The groups $\mathbb{Z}_4\bigstar_{\mathbb{Z}_2}D_2$]\label{ex7}
(Note that $D_2\cong \mathbb{Z}_2$, see in Appendix A, Tab. \ref{diehdral-groups}). This
group is generated by rotation of $\pi$, rotation by $\pi/2$ and
reflection over the $x$-axis. It can be considered a subgroup of
some crystallographic group $G(3)$, (see Appendix C).
\end{example}

\begin{definition}
The subgroups of $GL_d(\mathbb{Z})$ that are lattice
symmetry groups are called {\em Bravais subgroups}. Every maximal
finite subgroup of $GL_d(\mathbb{Z})$ is a {\em Bravais subgroup}.
\end{definition}

\begin{definition}
The {\em geometrical (arithmetical) holohedry} of a crystallographic
group $G(d)$ is the smallest Bravais subgroup $\widehat{G}(d)$
containing the point group $G\equiv G(d)/\mathbf{T}$ of $G(d)$.
\end{definition}

\begin{definition}
A crystallographic group $G(d)$ is said in {\em general position} if
 there is no affine transformation $\phi\in A(E)$ such that $\phi
 G(d)\phi^{-1}\equiv{}^\phi G(d)$ is a crystallographic group whose
 lattice of parallel translations has lower symmetry.
\end{definition}

\begin{proposition}
If the crystallographic group $G(d)$ is in general position, then
its holohedry $\widehat{G}(d)$ is the lattice symmetry group of the
parallel translations $\mathbf{T}$ of $G(d)$.
\end{proposition}

\begin{definition}
Two crystallographic groups belong to the same {\em syngony},
{\em(Bravais type)}, if their geometrical (arithmetical) holohedries
coincide.
\end{definition}

\begin{example}\label{ex8}
For the $3$-dimensional case there are $73$ arithmetic crystal
classes. The space group types with the same point group symmetry
and the same type of centering belong to the same arithmetic crystal
class. An arithmetic crystal class is indicated by the crystal
symbol of the corresponding point group followed by the symbol of
the lattice. Furthermore $G(3)$ has $7$ syngonies and $14$ Bravais
types of crystallographic groups. (See in Appendix A, Tab. \ref{230-affine-crystallographic-space-group-types-3-dim}. The number between brackets {\em()} after the symbol of the point group
is the number of space-group types with that point group.) The
geometric crystal classes with the point symmetries of the lattices
are called {\em holohedries} and are 7. The other 25 geometric
crystal classes are called {\em merhoedries}. The {\em Bravais
classes} (or {\em Bravais arithmetic crystallographic classes}) are
the arithmetic crystal classes with the point symmetry of the
lattice. The Bravais types of lattices and the Bravais classes have
the same point symmetry.
\end{example}

\section{\bf INTEGRAL BORDISM GROUPS VS. CRYSTALLOGRAPHIC GROUPS}

In the following we shall relate above crystallographic groups to
the geometric structure of PDE's. More precisely in some previous
works we have characterized the structure of global solutions of
PDE's by means of {\em integral bordism groups}. Let us start with
the bordism groups.\footnote{For general informations on bordism
groups, and related problems in differential topology, see, e.g.,
Refs.\cite{HIR, M-M, M-S, PRA7, QUIL, RUD, S-W, STO, SWI, THO1,
THO2, THO3, WAL1, WAL2}.}

\begin{theorem}[Bordism groups vs. crystallographic groups]\label{bordismgrcry}
Bordism groups of closed compact smooth manifolds can be
considered as subgroups of crystallographic ones. More precisely
one has the following:

{\em(i)} To each nonoriented bordism group $\Omega_n$, can be
canonically associated a crystallographic group $G(q)$,
{\em(crystal-group of $\Omega_n$)}, for a suitable integer $q$,
{\em(crystal dimension of $\Omega_n$)}, such that one has the
following split short exact sequence:
\begin{equation}\label{exactseq3}
\xymatrix{0\ar[r]&\mathbb{Z}^q\ar@<0.5ex>[r]&G(q)\ar[l]\ar@<0.5ex>[r]&\Omega_n\ar[l]\ar[r]&0\\}
\end{equation}
So $\Omega_n$ is at the same time a subgroup of its crystal-group,
as well as an extension of this last. If there are many crystal
groups satisfying condition in {\em(\ref{exactseq3})}, then we call
respectively {\em crystal dimension} and {\em crystal group} of
$\Omega_n$ the littlest one.

{\em(ii)} To each oriented bordism group ${}^+\Omega_n$,
$n\not\equiv 0$ mod $4$, can be canonically associated a
crystallographic group $G(d)$, {\em(crystal-group)} of
${}^+\Omega_n$, for a suitable integer $d$, {\em crystal
dimension)} of ${}^+\Omega_n$, such that one has the following
short exact sequence:
\begin{equation}
\xymatrix{0\ar[r]&{}^+\Omega_n\ar@<0.3ex>[r]&G(q)\ar@<0.3ex>[r]&G(d)/{}^+\Omega_n\ar[r]&0\\}
\end{equation}
So ${}^+\Omega_n$ is a subgroup of its crystal-group.
\end{theorem}

\begin{proof}
Let us recall the structure of bordism groups.

\begin{lemma}[Pontrjagin-Thom-Wall \cite{STO, SWI, WAL1}]\label{lem3}
A closed $n$-dimensional smooth manifold $V$, belonging to the
category of smooth differentiable manifolds, is bordant in this
category, i.e., $V=\partial M$, for some smooth $(n+1)$-dimensional
manifold $M$, iff the Stiefel-Whitney numbers $<w_{i_1}\cdots
w_{i_p},\mu_V>$ are all zero, where $i_1+\cdots +i_p=n$ is any
partition of $n$ and $\mu_V$ is the fundamental class of $V$.
Furthermore, the bordism group $\Omega_n$ of $n$-dimensional smooth
manifolds is a finite abelian torsion group of the form:
\begin{equation}
\Omega_n\cong\mathbb{Z}_2\bigoplus\cdots_{q}\cdots\bigoplus\mathbb{Z}_2,
\end{equation}
 where $q$ is the number of nondyadic partitions of $n$.\footnote{A
partition $(i_1,\cdots,i_r)$ of $n$ is nondyadic if none of the
$i_\beta$ are of the form $2^s-1$.} Two smooth closed
$n$-dimensional manifolds belong to the same bordism class iff all
their corresponding Stiefel-Whitney numbers are equal.
Furthermore, the bordism group ${}^+\Omega_n$ of closed
$n$-dimensional oriented smooth manifolds is a finitely generated
abelian group of the form:
\begin{equation}
{}^+\Omega_n\cong\mathbb{Z}\bigoplus\cdots
\bigoplus\mathbb{Z}\bigoplus\mathbb{Z}_2\bigoplus\cdots\bigoplus\mathbb{Z}_2,
\end{equation}
where infinite cyclic summands can occur only if $n\equiv 0$ mod
$4$. Two smooth closed oriented $n$-dimensional manifolds belong to
the same bordism class iff all their corresponding Stiefel-Whitney
and Pontrjagin numbers are equal.\footnote{{\em Pontrjagin numbers}
are determined by means of homonymous characteristic classes
belonging to $H^\bullet(BG,\mathbb{Z})$, where $BG$ is the
classifying space for $G$-bundles, with $G=S_p(n)$. See, e.g.,
Refs.\cite{M-M, M-S, RUD, STO, SWI, THO1, THO2, WAL1, WAL2}.}
\end{lemma}

Let us recall that a group $G$ is {\em cyclic} iff it is generated
by a single element. All cyclic groups $G$ are isomorphic either to
$\mathbb{Z}_p$, $p\in\mathbb{Z}^+$ or $\mathbb{Z}$.
$G\cong\mathbb{Z}_p$ iff there is some finite integer $q$ such that
$g^q=e$, for each $g\in G$. Here $e$ is the unit of $G$. A group $G$
is {\em virtually cyclic} if it has a cyclic subgroup $H$ of finite
index. All finite groups are virtually cyclic, since the trivial
subgroup $H=\{e\}$ is cyclic. Therefore, to compute the finite
virtually cyclic subgroups of the $d$-dimensional crystallographic
groups is equivalent to compute the finite subgroups. Of course can
be there also infinite virtually cyclic subgroups. For these we can
use the following lemma.

\begin{lemma}[Scott \& Wall \cite{S-W}]\label{lem4}
Given a group $G$, the only infinite virtually cyclic subgroups of
$G$ will be semidirect products $H\rtimes_\alpha\mathbb{Z}$ and
amalgamated free products $A\bigstar_BC$ for $H,A,B,C< G$.
\end{lemma}

The finite subgroups of the crystallographic groups may be derived
exclusively from their point groups. The following lemmas are useful
to gain the proof.

\begin{lemma}\label{lem5}
{\em(a)} Let $f=(0,b)$ be a finite-ordered element of a crystallographic group $G(d)$. Furthermore, let $x=(v,1)$ be a translation. Then, $f$ commutes with $x$ iff $b(v)=v$, i.e., $b$ fixes $v$.

{\em(b)} $F\times\mathbb{Z}$ is a subgroup of some $d$-dimensional
crystallographic group iff $F$ is a subgroup of some
$(d-1)$-dimensional crystallographic group.
\end{lemma}

\begin{proof}
(a) It follows directly by computation.

(b) Let assume that $F$ is a subgroup of some $(d-1)$-dimensional
crystallographic group $G(d-1)$, then $F\times\mathbb{Z}$ is a
subgroup of the $d$-dimensional crystallographic group
$G(d-1)\times\mathbb{Z}$. Vice versa, let $F\times\mathbb{Z}$ be a
subgroup of some $d$-dimensional crystallographic group $G'(d)$.
Since in this direct product the generator $x\in\mathbb{Z}$ commutes
with all $f\in F$ and $x$ belongs to the Bieberbach lattice
$\mathbf{T}'$  of $G'(d)$, there are at most $(d-1)$ independent
elements $t\in\mathbf{T}'$ which do not commute with a given $f\in
F$. Therefore, $F$ must be a subgroup of some $(d-1)$-dimensional
crystallographic group.
\end{proof}

\begin{lemma}\label{lem6}
If a crystallographic group $G(d)$ admits a subgroup $F\rtimes_\alpha
\mathbb{Z}$ for some finite group $F$ and some homomorphism
$\alpha:\mathbb{Z}\to Aut(F)$, then $G(d)$ also admits a subgroup
$F\times\mathbb{Z}$.
\end{lemma}

\begin{proof}
Since $\mathbb{Z}$ is cyclic, $\alpha$ is completely determined by
the $\alpha(1)$. Since $F$ is finite, $Aut(F)$ is also finite, hence
one has $\alpha(1)^q=1$ for some finite $q$. Therefore the elements
$y=x^{pq}\in\mathbb{Z}$, $\forall p\in\mathbb{Z}$, commute with any
$f\in F$, hence $F\rtimes_\alpha \mathbb{Z}$ contains as a subgroup
$F\times\mathbb{Z}$, identified with the couples $(f,y)$, were $y$
are above defined elements.
\end{proof}

\begin{lemma}\label{lem7}
Let $H=F\rtimes_\alpha\mathbb{Z}$, $\mathbb{Z}=<x>$, be a subgroup
of a crystallographic group $G(d)$. The elements $f\in F$ fix the shift
vectors $x^q\in \mathbf{T}$.
\end{lemma}

\begin{proof}
According to Lemma \ref{lem6} we can consider the subgroup $K=F\times\mathbb{Z}<H$.
Then the elements $z\in\mathbb{Z}$ of $K$
correspond to translation vectors $x^q\in\mathbb{Z}<H$. On the other
hand, since $K$ is a direct product $F\times\mathbb{Z}$, all $f\in
F$ commute with all $z\in K$. By Lemma \ref{lem5}(a), all $f\in F$ fix all
$x^q\in H$. \end{proof}
\begin{lemma}[Alperin \& Bell \cite{ALP-BEL}]\label{lem8}
Let $F$ and $H$ be groups, let $\alpha:H\to Aut(F)$ b a
homomorphism, and $\phi\in Aut(F)$. If $\hat\phi$ is the inner
automorphism of $Aut(F)$ induced by $\phi$, then
$F\rtimes_{\hat\phi\circ\alpha}H\cong F\rtimes_\alpha H$.

(This means that if we identify the conjugacy classes of
authomorphisms of a given group, we need to consider only one
element of each class to evaluate the candidacy of all automorphisms
in that class.)
\end{lemma}

\begin{lemma}\label{lem9}
If the presentation of the amalgamated free product contains two or
more elements of order two that do not commute, then the amalgamated
free product is not a subgroup of any three dimensional
crystallographic group.
\end{lemma}

\begin{proof}
In three dimension there are only three possible elements of order
two, inversions, $\pi$ rotation, and reflection. All these symmetry
commute. Therefore, an amalgamated free product with two order two
elements that do not commute cannot exist in $3$-dimensions.
\end{proof}

Let us first note that for bordism groups identified with some
finite or infinite cyclic groups, theorem is surely true by
considering the following two standard lemmas.

\begin{lemma}\label{lem10}
If $H$ is a finite subgroup of a group $G$, every element $a\in H$
generates a finite cyclic subgroup $<a>\equiv\mathbb{Z}_n\subset
H$, where $n$ is the order of $a$, and $a^{-1}=a^{n-1}$, or
equivalently $a^n=e$, where $e$ is the unit of $H$ (and also that
of $G$.)
\end{lemma}

\begin{lemma}\label{lem11}
Every element $a$ of a group $G$ generates a cyclic subgroup
$<a>\subset G$. If $a$ has {\em infinite order}, then
$<a>\cong\mathbb{Z}$.
\end{lemma}

Thus if $\Omega_p\cong \mathbb{Z}_2$, or ${}^+\Omega_p\cong
\mathbb{Z}$, it follows that theorem is proved.

Now, let us consider the more general situation. We shall consider
the following lemma.
\begin{lemma}\label{lem12}
The group $\mathbb{Z}^s\rtimes\mathbb{Z}^s_2$ can be considered a
crystallographic group in the Euclidean space $\mathbb{R}^s$.
\end{lemma}

\begin{proof}
The $\mathbb{Z}^s$-conjugacy classes of splittings of the split
extension
\begin{equation}\label{exactseq4}
\xymatrix{0\ar[r]&\mathbb{Z}^s\ar[r]&\mathbb{Z}^s\rtimes\mathbb{Z}_2^s\ar[r]&\mathbb{Z}_2^s\ar[r]&0\\}
\end{equation}
are in 1-1 correspondence with the elements of
\begin{equation}
  \left\{\begin{array}{ll}
           H^1(\mathbb{Z}_2^s;\mathbb{Z}^s)& =Hom_{\mathbb{Z}}(H_1(\mathbb{Z}_2^s;\mathbb{Z});\mathbb{Z}^s)\bigoplus Ext_{\mathbb{Z}}(H_0(\mathbb{Z}_2^s;\mathbb{Z});\mathbb{Z}^s)\\
           & \cong
 \left(Hom_{\mathbb{Z}}(H_1(\mathbb{Z}_2^s;\mathbb{Z});\mathbb{Z})\right)^{s}\bigoplus Ext_{\mathbb{Z}}(\mathbb{Z};\mathbb{Z})^s\\
           & \cong  Hom_{\mathbb{Z}}(\mathbb{Z}_2;\mathbb{Z})^{s^2}.\\
\end{array}
         \right.
\end{equation}

We have used the fact that for any finite cyclic group $K$ of order
$i$ one has
\begin{equation}
H_r(K;\mathbb{Z})=\left\{\begin{array}{ll}
                                      \mathbb{Z}&\hbox{\rm if $r=0$} \\
                                      \mathbb{Z}_i&\hbox{\rm if $r=$odd} \\
                                      0&\hbox{\rm if $r>0$ even.}
                                    \end{array}\right\}
\Rightarrow H_r(\mathbb{Z}_2;\mathbb{Z})=\left\{\begin{array}{ll}
                                      \mathbb{Z}&\hbox{\rm if $r=0$}  \\
                                      \mathbb{Z}_2&\hbox{\rm if $r=$odd} \\
                                      0&\hbox{\rm if $r>0$ even.}
                                    \end{array}\right\}.
\end{equation}
Furthermore, the K\"unneth theorem for groups allows us to write the
unnatural isomorphism:
\begin{equation}
H_r(G_1\times
G_2;\mathbb{Z})\cong\bigoplus_{p+q=r}H_p(G_1;\mathbb{Z})\otimes_{\mathbb{Z}}H_q(G_2;\mathbb{Z})
\bigoplus_{p+q=r-1}Tor^{\mathbb{Z}}(H_p(G_1;\mathbb{Z}),H_q(G_2;\mathbb{Z}))
\end{equation}
for any two groups $G_1$ and $G_2$. Taking into account that
$Tor^{\mathbb{Z}}(A,B)=0$ for projective $\mathbb{Z}$-module $A$ (or
$B$), and $Ext_{\mathbb{Z}}(\mathbb{Z};\mathbb{Z})=0$ since $\mathbb{Z}$ is a projective $\mathbb{Z}$-module, we get that $H^1(\mathbb{Z}_2^s;\mathbb{Z}^s)=Hom_{\mathbb{Z}}(\mathbb{Z}_2;\mathbb{Z})^{s^2}$. Thus we
conclude that there is not an unique $\mathbb{Z}^s$-conjugacy class in the splitting in (\ref{exactseq4}). Furthermore, the set
of equivalence classes of $s$-dimensional crystallographic groups
with such a point group are in 1-1 correspondence with
\begin{equation}
  \left\{\begin{array}{ll}
H^2(\mathbb{Z}_2^s;\mathbb{Z}^s)&\cong
Hom_{\mathbb{Z}}(H_2(\mathbb{Z}_2^s;\mathbb{Z});\mathbb{Z}^s)\bigoplus Ext_{\mathbb{Z}}(H_1(\mathbb{Z}_2^s;\mathbb{Z});\mathbb{Z}^s)\\
&\cong
Hom_{\mathbb{Z}}(H_2(\mathbb{Z}^s_2;\mathbb{Z});\mathbb{Z})^s\bigoplus Ext_{\mathbb{Z}}(\mathbb{Z}_2^{s};\mathbb{Z}^s)\\
&\cong \bigoplus_{1\le r\le (s-1)}Hom_{\mathbb{Z}}(\mathbb{Z}_2\otimes_{\mathbb{Z}}\mathbb{Z}_2^{s-r};\mathbb{Z})\bigoplus \mathbb{Z}_2^{s^2}.\\
\end{array}\right.
\end{equation}
The above calculation holds for $s>1$. Instead, for  $s=1$, one has $H^2(\mathbb{Z}_2;\mathbb{Z})=\mathbb{Z}_2$.
\end{proof}

Therefore if
$\Omega_n=\mathbb{Z}^s_2=\mathbb{Z}_2\times\cdots_s\cdots\times\mathbb{Z}_2$
it can be identified with a point group $G$ of a crystallographic
group $G(q)$, belonging to one of these equivalence classes, such
that $G(q)<M(\mathbb{R}^q)$. So one has the exact sequence
\begin{equation}\label{exactseq5}
\xymatrix{0\ar[r]&\mathbb{Z}^q\ar[r]&G(q)\ar[r]&\Omega_n\ar[r]&1
\\}
\end{equation}
that proves that $\Omega_n$ admits a crystallographic group as its
extension. Now, since $G(q)$ contains also as subgroup
$\mathbb{Z}^q\times\mathbb{Z}^q_2$, that contains as subgroup
$\mathbb{Z}^q_2$, it follows that $G(q)$ contains also as subgroup
$\Omega_n$. So one has also the following exact sequence:
\begin{equation}\label{exactseq6}
\xymatrix{0\ar[r]&\Omega_n\ar[r]&G(q)\ar[r]&G(q)/\Omega_n\ar[r]&0\\}
\end{equation}
Since $G(q)/\Omega_n\cong\mathbb{Z}^q$, it follows that sequence
(\ref{exactseq6}) is the split sequence of (\ref{exactseq5}), and vice versa.

Let us consider, now, the oriented case. Let us exclude the case
$n\equiv 0$ mod $4$. Then we can write
${}^+\Omega_n=\mathbb{Z}^r\times\mathbb{Z}^s_2$. Let us assume $r\ge
s$. Then we can consider in $\mathbb{R}^r$ the crystallographic
group $G(r)=\mathbb{Z}^r\rtimes\mathbb{Z}^s_2$. This contains as
subgroup $\mathbb{Z}^r\times\mathbb{Z}^s_2={}^+\Omega_n$. Therefore
one has the following exact sequence:
\begin{equation}
\xymatrix{0\ar[r]&{}^+\Omega_n\ar[r]&G(r)\ar[r]&G(r)/{}^+\Omega_n\ar[r]&0\\}
\end{equation}
Let us assume, now, that $r<s$, then one has the following sequence
of subgroups:
\begin{equation}
{}^+\Omega_n=\mathbb{Z}^r\times\mathbb{Z}^s_2<\mathbb{Z}^s\times\mathbb{Z}^s_2<\mathbb{Z}^s\rtimes\mathbb{Z}^s_2
=G(s).
\end{equation}
So we have the following short exact sequence
\begin{equation}
\xymatrix{0\ar[r]&{}^+\Omega_n\ar[r]&G(s)\ar[r]&G(s)/{}^+\Omega_n\ar[r]&0\\}
\end{equation}
Therefore we can conclude that ${}^+\Omega_n$ is a subgroup of the
crystallographic group $G(d)=\mathbb{Z}^d\rtimes\mathbb{Z}^d_2$,
with $d=\max\{r,s\}$.
\end{proof}

\begin{remark}
It is important emphasize that crystallographic groups and bordism
groups, even if related by above theorem, are in general different
groups. In other words, it is impossible identify any
crystallographic group with some bordism group, since the first
one is in general nonabelian, instead the bordism groups are
abelian groups.
\end{remark}

We can extend above proof also by including bordism groups
relatively to some manifold.

\begin{theorem}
Bordism groups relative to smooth manifolds can be considered as
extensions of crystallographic subgroups.
\end{theorem}

\begin{proof}
Let us recall that a {\em$k$-cycle} of $M$ be a couple $(N,f)$,
where $N$ is a $k$-dimensional closed (oriented) manifold and
$f:N\to M$ is a differentiable mapping. A {\em group of cycles}
$(N,f)$ of an $n$-dimensional manifold $M$ is the set of formal
sums $\sum_i(N_i,f_i)$, where $(N_i,f_i)$ are cycles of $M$. The
quotient of this group by the cycles equivalent to zero, i.e., the
boundaries, gives the {\em bordism groups}
$\underline{\Omega}_s(M)$. We define {\em relative bordisms}
$\underline{\Omega}_s(X,Y)$, for any pair of manifolds $(X,Y)$,
$Y\subset X$, where the boundaries are constrained to belong to
$Y$. Similarly we define the {\em oriented bordism groups}
${}^+\underline{\Omega}_s(M)$ and ${}^+\underline{\Omega}_s(X,Y)$.
One has $\underline{\Omega}_s(*)\cong\Omega_s$ and
${}^+\underline{\Omega}_s(*)\cong {}^+\Omega_s$. For bordisms, the
theorem of invariance of homotopy is valid. Furthermore, for any
CW-pair $(X,Y)$, $Y\subset X$, one has the isomorphisms:
$\underline{\Omega}_s(X,Y)\cong{\Omega}_s(X/Y)$, $s\ge 0$. One has
a natural group-homomorphism $\underline{\Omega}_s(X)\to
H_s(X;\mathbb{Z}_2)$. This is an isomorphism for $s=1$. In
general, $\underline{\Omega}_s(X)\not= H_s(X;\mathbb{Z}_2)$. In
fact one has the following lemma.

\begin{lemma}[Quillen \cite{QUIL}]\label{lem13}
One has the canonical isomorphism:
\begin{equation}
\underline{\Omega}_p(X)\cong \bigoplus_{r+s=p}
H_r(X;\mathbb{Z}_2)\otimes_{\mathbb{Z}_2}\Omega_s. \end{equation} In
particular, as $\Omega_0=\mathbb{Z}_2$ and $\Omega_1=0$, we get
$\underline{\Omega}_1(X)\cong H_1(X;\mathbb{Z}_2)$. Note that for
contractible manifolds, $H_s(X)=0$, for $s>0$, but
$\underline{\Omega}_s(X)$ cannot be trivial for any $s>0$. So, in
general, $\underline{\Omega}_s(X)\not= H_s(X;\mathbb{Z}_2)$.
\end{lemma}
So we get the following short exact sequence:
\begin{equation}
\xymatrix{0\ar[r]&\underline{K}_p(X)\ar[r]&\underline{\Omega}_p(X)\ar[r]&\Omega_p\ar[r]&0\\}
\end{equation}
where $\underline{K}_p(X)\cong \bigoplus_{r+s=p, r>0}
H_r(X;\mathbb{Z}_2)\otimes_{\mathbb{Z}_2}\Omega_s$. Therefore by
using Theorem \ref{bordismgrcry} we get the proof soon.
\end{proof}

Let us, now, consider a relation between PDE's and crystallographic
groups. This will give us also a new classification of PDE's on the
ground of their integral bordism groups.

\begin{definition}
We say that a PDE $E_k\subset J^k_n(W)$ is an {\em extended
0-crystal PDE}, if its integral bordism group is zero.
\end{definition}

The first main theorem is the following one relating the
integrability properties of a PDE to crystallographic groups.

\begin{theorem}[Crystal structure of PDE's]\label{main1}
Let $E_k\subset J^k_n(W)$ be a formally integrable and completely
integrable PDE, such that $\dim E_k\ge 2n+1$. Then its integral
bordism group $\Omega_{n-1}^{E_k}$ is an extension of some
crystallographic subgroup. Furthermore if $W$ is contractible, then
$E_k$ is an extended $0$-crystal PDE, when $\Omega_{n-1}=0$.
\end{theorem}

\begin{proof}
Let us first recall some definitions and results about integral
bordism groups of PDE's. (For details see Refs.\cite{PRA7, PRA9,
PRA10, PRA11, PRA13, PRA15, PRA16}, and the following related works
\cite{AG-PRA1, AG-PRA2}.) Let $W$ be an $(n+m)$-dimensional smooth
manifold with fiber structure $\pi:W\to M$ over a $n$-dimensional
smooth manifold $M$. Let $E_k\subset J^k_n(W)$ be a PDE of order
$k$, for $n$-dimensional submanifolds of $W$. For an "admissible"
$p$-dimensional, $p\in\{0,\cdots,n-1\}$, integral manifold $N\subset
E_k$, we mean a $p$-dimensional smooth submanifold of $E_k$,
contained in an admissible integral manifold $V\subset E_k$, of
dimension $n$, i.e. a solution of $E_k$, that can be deformed into
$V$, in such a way that the deformed manifold $\widetilde{N}$ is
diffeomorphic to its projection
$\widetilde{X}\equiv\pi_{k,0}(\widetilde{N})\subset W$. In such a
case the $k$-prolongation $\widetilde{X}^{(k)}=\widetilde{N}$. The
existence of $p$-dimensional admissible integral manifolds $N\subset
E_k$ is obtained solving Cauchy problems of dimension
$p\in\{0,\cdots,n-1\}$, i.e., finding $n$-dimensional admissible
integral manifolds (solutions) of a PDE $E_k\subset J^k_n(W)$, that
contain $N$. Existence theorems for such solutions can be studied in
the framework of the geometric theory of PDE's. For a modern
approach, founded on webs structures see also \cite{AG-PRA1,
AG-PRA2}. A geometric way to study the structure of global solutions
of PDE's, is to consider their integral bordism groups. Let
$N_i\subset E_k$, $i=1,2$ be two $(n-1)$-dimensional compact closed
admissible integral manifolds. Then, we say that they are
{\em$E_k$-bordant} if there exists a solution $V\subset E_k$, such
that $\partial V=N_1\DU N_2$ (where $\DU$ denotes disjoint union).
We write $N\thicksim_{E_k}N_2$. This is an equivalence relation and
we will denote by $\Omega^{E_k}_{n-1}$ the set of all $E_k$-bordism
classes $[N]_{E_k}$ of $(n-1)$-dimensional compact closed admissible
integral submanifolds of $E_k$. The operation of taking disjoint
union defines a structure of abelian group on $\Omega^{E_k}_{n-1}$.
In order to distinguish between integral bordism groups where the
bording manifolds are smooth, (resp. singular, resp. weak), we shall
use also the following symbols $\Omega^{E_k}_{n-1}$, (resp.
$\Omega^{E_k}_{n-1,s}$, resp. $\Omega^{E_k}_{n-1,w}$).\footnote{Let
us recall that {\em weak solutions}, are solutions $V$, where the
set $\Sigma(V)$ of singular points of $V$, contains also
discontinuity points, $q,q'\in V$, with
$\pi_{k,0}(q)=\pi_{k,0}(q')=a\in W$, or $\pi_{k}(q)=\pi_{k}(q')=p\in
M$. We denote such a set by $\Sigma(V)_S\subset\Sigma(V)$, and, in
such cases we shall talk more precisely of {\em singular boundary}
of $V$, like $(\partial V)_S=\partial V\setminus\Sigma(V)_S$.
However for abuse of notation we shall denote $(\partial V)_S$,
(resp. $\Sigma(V)_S$), simply by $(\partial V)$, (resp.
$\Sigma(V)$), also if no confusion can arise. Solutions with such
singular points are of great importance and must be included in a
geometric theory of PDE's too \cite{PRA15}. Let us also emphasize
that singular solutions can be identified with integral
$n$-chains in $E_k$, and in this category can be considered also
{\em fractal solutions}, i.e., solutions with sectional fractal or
multifractal geometry. (For fractal geometry see, e.g., \cite{FALC,
LAP-FRA, MANDE}.)} Let us first consider the integral bordism group
$\Omega_{n-1,w}^{E_k}$ (or $\Omega_{n-1,s}^{E_k}$), for weak
solutions (or for singular solutions). We shall use Theorem 2.15 in
\cite{PRA15}, that we report below to be more direct. (See also
\cite{AG-PRA1, AG-PRA2}.)

{\em Theorem 2.15 in \cite{PRA15}. Let $E_k\subset J^k_n(W)$ be a
formally and also completely integrable PDE, such that $\dim E_k\ge
2n+1$. Then one has the following canonical isomorphism:
\begin{equation}
    \Omega^{E_k}_{p,w}\cong\bigoplus_{r,s,r+s=p}H_r(W;\mathbb{Z}_2)\oplus_{\mathbb{Z}_2}\Omega_s.
\end{equation}

Furthermore, if $W$ is an affine fiber bundle $\pi:W\to M$ over a
$n$-dimensional manifold $M$, one has the isomorphisms:\footnote{The
bording solutions considered for the bordism groups are singular
solutions if the symbols $g_k$ and $g_{k+1}$ are different from
zero, and for singular-weak solutions in the general case. Here we
have denoted $\underline{\Omega}_p(X)$ the $p$-bordism group of a
manifold $X$. For informations on such structure of the algebraic
topology see. e.g. \cite{HIR, PRA7, QUIL, RUD, STO, SWE, SWI, THO1,
THO2, WAL1, WAL2}.}}
\begin{equation}
\Omega_p^{E_k}\cong\underline{\Omega}_p(M)\cong\bigoplus_{r,s,
r+s=p}H_r(M;\mathbb{Z}_2)\oplus_{\mathbb{Z}_2}\Omega_s.
\end{equation}

So we can write for the weak integral bordism group
\begin{equation}
\Omega_{n-1,w}^{E_k}\cong \bigoplus_{r+s=n-1}
H_r(W;\mathbb{Z}_2)\otimes_{\mathbb{Z}_2}\Omega_s.
\end{equation}
Then one has that $\Omega_{n-1,w}^{E_k}$ is an extension of the
bordism group $\Omega_{n-1}$. Hence by using Theorem
\ref{bordismgrcry} we get the proof.

Let us now consider the integral bordism group $\Omega_{n-1}^{E_k}$
for smooth solutions. This bordism group is related to previous one,
and to singular bordism group $\Omega_{n-1,s}^{E_k}$, by means of
the exact commutative diagram (\ref{comdiag2}). Furthermore, the
relation between $\Omega_{n-1}^{E_k}$, $\Omega_{n-1,w}^{E_k}$ and
$\Omega_{n-1}$ is given by means of the exact commutative diagram
(\ref{comdiag3}), where $K^{E_k}_{n-1,w}=\ker(a)$,
$K^{E_k}_{n-1,w;n-1}=\ker(b)$, $\overline{K}^{E_k}_{n-1}=\ker(c)$,
with $c=b\circ a$. From this we get that also $\Omega^{E_k}_{n-1}$
can be considered an extension of $\Omega_{n-1}$ if
$\Omega^{E_k}_{n-1,w}$ is so. Therefore we can apply Theorem
\ref{bordismgrcry} also to the integral bordism group for smooth
solutions whether it can be applied to the integral bordism group
for weak (or singular) solutions. Therefore theorem is proved.

\begin{equation}\label{comdiag2}
\xymatrix{
&0\ar[d]&0\ar[d]&0\ar[d]&\\
0\ar[r]&K^{E_k}_{n-1,w/(s,w)}\ar[d]\ar[r]&
K^{E_k}_{n-1,w}\ar[d]\ar[r]&K^{E_k}_{n-1,s,w}\ar[d]\ar[r]&0\\
0\ar[r]&K^{E_k}_{n-1,s}\ar[d]\ar[r]&
\Omega^{E_k}_{n-1}\ar[d]\ar[r]&\Omega^{E_k}_{n-1,s}\ar[d]\ar[r]&0\\
&0\ar[r]&\Omega^{E_k}_{n-1,w}\ar[d]\ar[r]&\Omega^{E_k}_{n-1,w}\ar[d]\ar[r]&0\\
&&0&0&}
\end{equation}

\begin{equation}\label{comdiag3}
\xymatrix{0\ar[rd]&&&0\ar[d]&\\
&\overline{K}^{E_k}_{n-1}\ar[rd]\ar[rr]&&K^{E_k}_{n-1,w;n-1}\ar[d]\ar[r]&0\\
0\ar[r]&K^{E_k}_{n-1,w}\ar[u]\ar[r]&\Omega^{E_k}_{n-1}\ar[rd]^{c}\ar[r]^{a}&\Omega^{E_k}_{n-1,w}\ar[d]^{b}\ar[r]&0\\
&0\ar[u]\ar[rr]&&\Omega_{n-1}\ar[d]&\\
&&&0&\\}
\end{equation}
\end{proof}

\begin{definition}
We say that a PDE $E_k\subset J^k_n(W)$ is an {\em extended crystal
PDE}, if conditions of  above theorem are verified. We define {\em
crystal group of $E_k$} the littlest crystal group such Theorem \ref{main1}
is satisfied. The corresponding dimension will be called {\em
crystal dimension of $E_k$}.
\end{definition}

In the following we relate crystal structure of PDE's to the
existence of global smooth solutions, identifying an
algebraic-topological obstruction.

\begin{theorem}\label{main2}
Let $E_k\subset J^k_n(W)$ be a formally integrable and completely
integrable PDE. Then, in the algebra $\mathbf{H}_{n-1}(E_k)\equiv
Map(\Omega_{n-1}^{E_k};\mathbb{R})$, {\em(Hopf algebra of $E_k$)},
there is a subalgebra, {\em(crystal Hopf algebra)} of $E_k$. On such
an algebra we can represent the algebra $\mathbb{R}^{G(d)}$
associated to the crystal group $G(d)$ of $E_k$. (This justifies the
name.) We call {\em crystal conservation laws} of $E_k$ the elements
of its crystal Hopf algebra.
\end{theorem}

\begin{proof}
In fact the short exact sequence
\begin{equation}
\xymatrix{0\ar[r]&\overline{K}^{E_k}_{n-1}\ar[r]^{d}&\Omega_{n-1}^{E_k}\ar[r]^{e}&\Omega_{n-1}\ar[r]&0\\}
\end{equation}
obtained from the commutative diagram in (\ref{comdiag3}),
identifies for duality the following sequence

\begin{equation}\label{dual-sequence-a}
\xymatrix{0&\ar[l]\mathbf{K}_{n-1}^{E_k}&\ar[l]_{d_*}\mathbf{H}_{n-1}(E_k)&\ar[l]_{e_*}\mathbb{R}^{\Omega_{n-1}}&\ar[l]0\\}
\end{equation}
On the other hand from the short exact sequence
\begin{equation}
\xymatrix{0\ar[r]&\Omega_{n-1}\ar[r]^{f}&G(d)\ar[r]^{g}&G(d)/\Omega_{n-1}\ar[r]&1\\}
\end{equation}
identifying $\Omega_{n-1}$ with a crystallographic subgroup, we get
for duality the following sequence
\begin{equation}\label{dual-sequence-b}
\xymatrix{0&\ar[l]\mathbb{R}^{\Omega_{n-1}}&\ar[l]_{f_*}\mathbb{R}^{G(d)}&\ar[l]_{g_*}\mathbb{R}^{G(d)/\Omega_{n-1}}&\ar[l]0\\}.
\end{equation}
So we can identify the crystal Hopf algebra of $E_k$ with
$\mathbb{R}^{\Omega_{n-1}}$. On such an algebra we can represent all
the Hopf algebra $\mathbb{R}^{G(d)}$ associated to the crystal group
$G(d)$ of $E_k$.\footnote{Let us remark that sequences (\ref{dual-sequence-a}) and (\ref{dual-sequence-b}) do not necessitate to be exact, but are always partially exact. In fact, even if it does not necessitate that $d_\circ e_*=0$ and $f_*\circ g_*=0$, (hence neither $\IM(e_*)=\ker(d_*)$ and $\IM(g_*)=\ker(f_*)$), one has that $e_*$ and $g_*$ are monomorphisms and $d_*$ and $f_*$ are epimorphisms. This is enough for our proof.
Let us recall also that $\mathbf{H}_{n-1}(E_k)\equiv
Map(\Omega_{n-1}^{E_k};\mathbb{R})$ is an Hopf algebra in extended
sense, i.e. it contains the Hopf algebra $\mathbb{R}^{\Omega_{n-1}}$
as a subalgebra. (See also \cite{PRA10}.)}
\end{proof}

\begin{theorem}\label{main3}
Let $E_k\subset J^k_n(W)$ be a formally integrable and completely
integrable PDE. Then, the obstruction to find global smooth
solutions of $E_k$ can be identified with the quotient
$\mathbf{H}_{n-1}(E_\infty)/\mathbb{R}^{\Omega_{n-1}}$.
\end{theorem}

\begin{proof}
Let us first consider the following lemma that gives some criteria
to recognize $(n-1)$-dimensional admissible integral manifolds in
$E_k$.

\begin{lemma}[Cauchy problem solutions criteria]\label{cauchy-criteria}
{\em 1)} Let $E_k\subset J^k_n(W)$ be a formally integrable and completely
integrable PDE on the fiber bundle $\pi:W\to M$, $\dim W=m+n$, $\dim
M=n$. Let $N\subset E_k$ be a smooth $(n-1)$-dimensional integral
manifold, diffeomorphic to $X\equiv\pi_{k,0}(N)\subset W$, such that $X\subset Y$, where $Y$ is a smooth $n$-dimensional submanifold of $W$, satisfying the condition that its $k$-order prolongation $Y^{(k)}\subset J^k_{n-1}(W)$, contains $N$, {\em($0$-order admissibility condition)}. Then,
there exists a (weak, or singular, or smooth) solution $V\subset
E_k$, such that $N\subset V$.

{\em 2)} Furthermore, if the symbols $g_k$ and $g_{k+1}$, of $E_k$ and
$E_{k+1}$ respectively, are different from zero, then $V$ can be a
singular or smooth solution. Moreover, if there exists a nonzero
smooth vector field $\zeta:E_k\to TE_k$, transversal
to $N$, and characteristic, at least for a sub-equation of $E_k$, then a smooth solution $V$ passing through $N$ can be built by means of the flow $\phi$ associated to $\zeta$. For suitable initial Cauchy integral manifolds, solutions can be built by using infinitesimal symmetries of the closed ideal encoding $E_k$. (For details see Theorem 2.10 in \cite{PRA10}.)

{\em 3)} Finally, let $E_k\to E_{k-1}\equiv\pi_{k,k-1}(E_k)$, be an affine
subbundle of $\pi_{k,k-1}:J^k_n(W)\to J^{k-1}_n(W)$, with associated
vector bundle $\bar\pi_{k-1}:g_k\to E_{k-1}$, where
$\bar\pi_{k,k-1}=\pi_{k,k-1}\circ\bar\pi_k$, with $\bar\pi_k:g_k\to
E_k$ is the canonical projection. Let $N\subset E_k$ be a smooth
$(n-1)$-dimensional integral manifold, diffeomorphic to
$X\equiv\pi_{k,0}(N)\subset W$ and satisfying the $0$-order admissibility condition. Then, there exists a (singular, or smooth) solution $V\subset E_k$, such that $N\subset V$.

An integral manifold $N\subset E_k$, as above defined and contained
into a solution $V\subset E_k$, is called {\em admissible}.
\end{lemma}

\begin{proof}
Since $N$ is an $(n-1)$-dimensional smooth integral manifold,
diffeomorphic to $X\equiv\pi_{k,0}(N)\subset W$, satisfying the $0$-order admissibility condition, we can consider $N\subset Y^{(k)}\subset J^k_n(W)$, where $Y$ is just a $n$-dimensional smooth manifold of $W$, containing $X$. In general $Y^{(k)}\not\subset E_k$, but taking into account that $E_k$ is
formally integrable and completely integrable, we get that
$(E_k)_{+r}$ is a strong retract of $J^{k+r}_n(W)$, $\forall r>0$. Then, we can
deform $Y^{(k+r)}\subset J^{k+r}_n(W)$ into $E_{k+r}$, obtaining
a (weak, or singular, or smooth) solution $\tilde Y\subset E_{k+r}$,
passing for $N^{(r)}\cong X^{(k+r)}$. Then $\pi_{k+r,k}(\tilde
Y)\equiv V\subset E_k$ is a solution of $E_k$, passing through $N$.
In particular, if $g_k\not=0$, and $g_{k+r}\not=0$, $\tilde Y$ is a
singular (or smooth) solution, and so can be also $V$. Moreover, in
the case that $N$ is transversal to a characteristic smooth vector
field $\zeta:E_k\to TE_k$, then
$V=\bigcup_{\lambda\in]-\epsilon,\epsilon[}\phi_\lambda(N)$ is a
smooth solution of $E_k$ passing through $N$. Under suitable conditions on the Cauchy integral manifold, the vector field used to build the solution can be an infinitesimal symmetry. (For a detailed proof see the one of Theorem 2.10 in \cite{PRA10}.)

Finally, whether $E_k\to E_{k-1}\equiv\pi_{k,k-1}(E_k)$, is an
affine subbundle of $\pi_{k,k-1}:J^k_n(W)\to J^{k-1}_n(W)$, with
associated vector bundle $\bar\pi_{k-1}:g_k\to E_{k-1}$, then also
$E_k$ is a strong retract of $J^k_n(W)$, so we can reproduce above
strategy used to build a solution passing for $N$, without the necessity to prolong $E_k$, (if it is not differently required by the structure of this equation).
\end{proof}

\begin{example}[Fourier's heat equation]\label{fourier-heat-equation}
Let us consider the second order PDE
\begin{equation}\label{Fourier}
    (F)\subset J{\it D}^2(W)\subset J^2_2(W): F\equiv u_t-u_{xx}=0
\end{equation}
on the fiber bundle $\pi:W\equiv\mathbb{R}^3\to\mathbb{R}^2$,
$(t,x,u)\mapsto(t,x)$. In this case $(F)\to \pi_{2,1}(F)\equiv(F)_{-1}=J{\it D}(W)$
is an affine fiber subbundle of the affine fiber bundle
$\pi_{2,1}:J{\it D}^2(W)\to J{\it D}(W)$, with associated vector
bundle identified with the symbol $g_2$: $\zeta=\zeta_{tt}\partial
u_{tt}+\zeta_{tx}\partial u_{tx}\in g_2$. So we can apply Lemma
\ref{cauchy-criteria}(3). Therefore, if $N\subset (F)$ is a
(compact) $1$-dimensional integral manifold, diffeomorphic to its
projection $\pi_{2,0}(N)\equiv X\subset W$, we can find solutions
$V\subset (F)$, passing from $N$. In particular, whether $N$ is
diffeomorphic to a smooth space-like curve $u=h(x)$, (at $t=0$), we get that $N$
is the image of a mapping, say
$\gamma:I\equiv[0,1]\subset\mathbb{R}\to (F)$, iff
$\gamma^*\mathcal{C}_2=0$, where $\mathcal{C}_2=<dF,\omega\equiv
du-u_tdt-u_xdx,\omega_t\equiv du_t-u_{tt}dt-u_{tx}dx,\omega_x\equiv
du_x-u_{xt}dt-u_{xx}dx>$ is the Pfaffian contact ideal encoding
solutions of $(F)$. Then one can see that the integral curve
$\gamma$ is represented in coordinates
$(t,x,u,u_t,u_x,u_{tt},u_{tx},u_{xx})$ on $ J{\it D}^2(W)$, by the following
equations:
\begin{equation}\label{integral-curve}
    \left\{
\begin{array}{l}
  t\circ\gamma =0;\hskip 0.5cm
   x\circ\gamma =x;\hskip 0.5cm
    u\circ\gamma =h(x);\hskip 0.5cm
    u_t\circ\gamma =\frac{d^2h}{dx^2}(x);\hskip 0.5cm
      u_x\circ\gamma =\frac{dh}{dx}(x);\\
u_{xx}\circ\gamma =\frac{d^2h}{dx^2}(x);\hskip 0.5cm u_{tx}\circ\gamma =\frac{d^3h}{dx^3}(x);\hskip 0.5cm
u_{tt}\circ\gamma =\kappa(x).\\
\end{array}
    \right.
\end{equation}
where $\kappa(x)$ is an arbitrary smooth function. Then, by
considering that the Fourier's heat equation is a formally
integrable and completely integrable PDE, we can see, by taking the
first and second prolongations of $(F)$, that must necessarily be
$\kappa(x)=\frac{d^4h}{dx^4}(x)$. In fact, from the first prolongation we get that must be $u_{tt}=u_{xxt}$. From the second prolongation of $(F)$ one has $u_{txx}=u_{xxxx}$. Since must be $u_{xxt}=u_{txx}$ we get $u_{tt}=u_{xxxx}$. By restriction on $t=0$, one has $u_{tt}|_{t=0}=u_{xxxx}|_{t=0}=h_{xxxx}$.
So we can build solution also by
using Lemma \ref{cauchy-criteria}(1). For example, if we are interested to a solution $V_S\subset (F)$, obtained by means of a rigid propagation of the initial Cauchy space-like integral curve {\em(\ref{integral-curve})}, one can see that must necessarily be  $h(x)=\alpha x+\beta$, with $\alpha,\beta\in\mathbb{R}$, $x\in[0,1]\subset\mathbb{R}$. In other words such type of steady-state solution, $u(t,x)=\alpha x+\beta$, determines also the admissible structure of the integral initial Cauchy line $\Gamma\subset (F)_{t=0}$. In such a case this must be given by the following parametric equations in $ J{\it D}^2(W)$:
\begin{equation}\label{integral-curve-steady-state}
    (\Gamma):\hskip 2pt\left\{
\begin{array}{l}
t\circ\gamma=0;\hskip 0.5cm x\circ\gamma=x;\hskip 0.5cm u\circ\gamma=\alpha x+\beta;\hskip 0.5cm u_t\circ\gamma=0;\hskip 0.5cm u_x\circ\gamma=\alpha;\\
u_{xx}\circ\gamma=0;\hskip 0.5cm u_{tx}\circ\gamma=0;\hskip 0.5cm u_{tt}\circ\gamma=0.\\
\end{array}
    \right.
\end{equation}
Remark that $\zeta=\partial t$ is the smooth vector field propagating the initial Cauchy curve $\Gamma$ in $V_S$, i.e., $\zeta$ is the {\em characteristic vector field} of $V_S$. This is not a characteristic vector field for the Fourier's heat equation and neither it is a characteristic vector field for the steady-state sub-equation $\left\{u_{xx}=0\right\}$. Instead,  $\zeta=\partial t$ is an infinitesimal symmetry for $(F)$. (See Lemma \ref{cauchy-criteria}(2).) Let us also emphasize that the above steady-state solution is not the unique solution passing for $\Gamma$. In order to see this it is enough to generalize the concept of solution and to consider also weak-singular solutions. In fact, let us find regular perturbations of the steady-state solution $V_S$ that at $t=0$, in correspondence of the boundary points $\{A,B\}=\partial\Gamma$ of $\Gamma$, have values respectively $a$ and $a+b$. These can be obtained by considering the Fourier's heat equation, that is a linear equation, as its linearized at the steady-state solution $V_S$, $(F)[V_S]=(F)$, i.e., the equation for perturbations of $V_S$. Then we get:
\begin{equation}\label{solution-for-two-initial-points}
\nu(t,x|\mu,a,b)\equiv e^{-\mu^2t}\left[\frac{a+b-b\cos(\mu)}{\sin(\mu)}\sin(\mu x)+b\cos(\mu x)\right].
\end{equation}
 $\nu(t,x|\mu,a,b)$ depends on an arbitrary positive parameter $\mu>0$, and two other parameters $a,b\in\mathbb{R}$, and has the following limits: $$\lim_{t\to+\infty}\nu(t,x|\mu,a,b)=0;\hskip 2pt\lim_{\mu\to 0}\nu(t,x|\mu,a,b)=a x+b;\hskip 2pt\lim_{a,b\to 0}\nu(t,x|\mu,a,b)=0. $$
 This means that with respect to these perturbations the steady-state solution is asymptotically stable, and all such perturbations produce deformations of the steady-state solution, with deformation parameters just $(\mu,a,b)$. Let us denote by $V[\mu,a,b]$ the $2$-dimensional integral manifold, contained in $(F)$, representing the deformed steady-state solution by $\nu(t,x|\mu,a,b)$, $0\le x\le 1$. More precisely $V[\mu,a,b]$ is the integral manifold of $(F)$ corresponding to the solution $u(t,x|\mu,a,b)\equiv \alpha x+\beta+\nu(t,x|\mu,a,b)$. Set $V=V_S\bigcup V[\mu,a,b]$. $V$ is a weak-singular solution of $(F)$, that for $t=0$ passes for $\Xi^{(2)}=\Gamma^{(2)}\bigcup \Gamma[\mu,a,b]^{(2)}$, a $1$-dimensional admissible integral manifold of $(F)$, such that $\pi_{2,0}(\Xi^{(2)})\equiv\Xi=\Gamma\bigcup \Gamma[\mu,a,b]$, where $\Gamma[\mu,a,b]$ is the curve identified by $ V[\mu,a,b]$ for $t=0$.  $\Xi$ is a space-like curve contained in $W$, that is projected, by means of $\pi$, on the interval $[0,1]$ of the $x$-axis. The $2$-dimensional manifold $Y\equiv\pi_{2,0}(V)=\pi_{2,0}(V_S)\bigcup\pi_{2,0}(V[\mu,a,b])\equiv Y_S\bigcup Y[\mu,a,b]$, passes for $\Xi$. Furthermore $V$ converges, for $t\to+\infty$, to the steady-state solution $V_S$.\footnote{Singular solutions, like those described in this example, are very important in many physical applications too, since they represent complex phenomena related to perturbations of some fixed dynamic background. For example, for suitable values of the parameters $a$ and $b$ in {\em(\ref{solution-for-two-initial-points})}, manifolds $Y[\mu,a,b]$ intersect $Y_S$ along common characteristic lines, and the singular solutions $V$ are piecewise $\mathbb{Z}_2$-manifolds. (For complementary informations on such singular manifolds, and singular solutions of PDE's, see also \cite{BUON-ROU-SAND, MO-SU, PRA10, PRA17}.)}

 Let us, now, show also an explicit construction of solution for $(F)$ by means of the method of the retraction given in Lemma  \ref{cauchy-criteria}. For example, let $f(P)\equiv V\subset J^2_2(W)$, $P\subset \mathbb{R}^2$, be a $2$-dimensional integral manifold representing the smooth function $u(t,x)\equiv \sum_{n\ge 0}t^n a_n(x)$, passing for the integral curve given in {\em(\ref{integral-curve})}. The coefficients $a_n(x)$ are suitable functions on $x\in[0,1]$, such that the series converges in $t\in\mathbb{R}^+$, and such that $a_0(x)=h(x)$, $a_1(x)=h_{xx}$, $a_2(x)=\frac{1}{2}h_{xxxx}$. These last conditions assure that $V$ passes for the integral curve {\em(\ref{integral-curve})}. The integral manifold $V$ is not a solution of $(RF)$ since there are not restrictions on the other coefficients $a_n$, $n\ge 3$. The parametric equations of $V$ are given in {\em(\ref{parametric-equation-integral-manifold-solution-jet-space-example-fourier-equation}(a))}.
  \begin{equation}\label{parametric-equation-integral-manifold-solution-jet-space-example-fourier-equation}
\scalebox{0.8}{$\hbox{\rm(a)}\left\{
\begin{array}{l}
  t\circ f=t\\
  x\circ f=x\\
  u\circ f=\sum_{0\le n\le 2}t^n\frac{1}{n!} \frac{d^{2nh}}{dx^{2n}}(x)+\sum_{n\ge 3}t^n a_n(x) \\
  \\
  u_t\circ f=\sum_{1\le n\le 2}t^{n-1}\frac{d^{2nh}}{dx^{2n}}(x)+\sum_{n\ge 3}nt^{n-1}a_n(x)\\
  \\
  u_x\circ f=\sum_{0\le n\le 2 }\frac{t^n}{n!}\frac{d^{2n+1}h}{dx^{2n+1}}(x)+\sum_{n\ge 3}t^n\frac{d a_n}{dx}(x)\\
  \\
  u_{xx}\circ f=\sum_{0\le n\le 2}\frac{t^n}{n!}\frac{d^{2n+2}h}{dx^{2n+2}}(x)+\sum_{n\ge 3}t^n\frac{d^2 a_n}{dx^2}(x)\\
  \\
   u_{tx}\circ f=\sum_{1\le n\le 2}t^{n-1}\frac{d^{2n+1}h}{dx^{2n+1}}(x)+\sum_{n\ge 3n}nt^{n-1}\frac{d a_n}{dx}(x)\\
   \\
    u_{tt}\circ f=\frac{d^{4}h}{dx^{4}}(x)+\sum_{n\ge 3}n(n-1)t^{n-2}a_n(x).\\
\end{array}
\right\}
\Rightarrow
\hbox{\rm(b)}\left\{
\begin{array}{l}
  t\circ r\circ f=t\\
  x\circ r\circ f=x\\
  u\circ r\circ f=\widetilde{u}(t,x)\\
  u_t\circ r\circ f=\widetilde{u}_t(t,x)\\
  u_x\circ r\circ f=\widetilde{u}_x(t,x)\\
  u_{tt}\circ r\circ f=\widetilde{u}_{tt}(t,x)\\
  u_{tx}\circ r\circ f=\widetilde{u}_{tx}(t,x)\\
  \end{array}
\right\}.$}
\end{equation}

 Taking into account that on $(F)$ we can consider the following coordinate functions $\{t,x,u,u_t,u_x,u_{tt},u_{tx}\}$, it follows that the retraction mapping $r(f(P))=r(V)\equiv \widetilde{V}\subset (F)$ has the parametric equation
given in {\em(\ref{parametric-equation-integral-manifold-solution-jet-space-example-fourier-equation}(b))},
where the function $\widetilde{u}(t,x)$ is determined starting from the function $u(t,x)$, by imposing the condition to belong to $(F)$.  So we get $\widetilde{u}(t,x)=\sum_{n\ge 0}t^n \widetilde{a}_n(x)$, with $\widetilde{a}_n\equiv \frac{1}{n!}\frac{d^{2n}h}{dx^{2n}}(x)$.\footnote{This function $\widetilde{u}(t,x)$, is well defined and limited in all $(t,x)\in\mathbb{R}^+\times[0,1]$, since the function $h(x)$ is smooth in $[0,1]$ and with $|\frac{d^{2n}h}{dx^{2n}}(x)|\le C\in\mathbb{R}$, $\forall x\in[0,1]$. In fact one has $\sum_{n\ge 0}t^n\widetilde{a}_n\le C\sum_{n\ge 0}\frac{t^n}{n!}\equiv C\sum_{n\ge 0}b_n$. The last series converges, with convergence radius $r=\infty$ since $\lim_{n\to\infty}\left|\frac{b_{n+1}}{a_n}\right|=\lim_{n\to\infty}\frac{t}{n+1}=0=\frac{1}{r}$. Therefore the series $\sum_{n\ge 0}t^n \widetilde{a}_n(x)$ is convergent as it is absolutely convergent.}
$\widetilde{V}\subset(F)$, represents the solution $\widetilde{u}(t,x)$. This can be considered as a deformation of $u(t,x)$. In fact, the functions defined in {\em(\ref{deformations-retractions-solutions-example-fourier-equation})}, is just the explicit deformation connecting $\widetilde{u}(t,x)$ with $u(t,x)$.
\begin{equation}\label{deformations-retractions-solutions-example-fourier-equation}
\begin{array}{l}
\widetilde{u}(t,x|\lambda)=\sum_{n\ge 0}t^n \widetilde{a}_n(\lambda|x),\hskip 2pt \widetilde{a}_n(\lambda|x)\equiv a_n(x)-\lambda[a_n(x)+\widetilde{a}_n(x)]\\
\Downarrow\\
\left\{
\begin{array}{l}
\widetilde{u}(t,x|0)=\sum_{n\ge 0}t^n a_n(x)\\
\\
\widetilde{u}(t,x|1)=\sum_{n\ge 0}t^n \widetilde{a}_n(x).\\
 \end{array}
\right.\\
\end{array}
\end{equation}
The integral manifolds $\widetilde{V}_{\lambda}\subset J^2_2(W)$, $\lambda\in[0,1]$, corresponding to $\widetilde{u}(t,x|\lambda)$, are not contained in $(F)$ for all $\lambda\in[0,1]$, hence are not solutions of $(F)$.
Therefore, $\widetilde{V}\subset(F)$ is a solution for the corresponding Cauchy-problem, obtained by means of the retraction method given in Lemma \ref{cauchy-criteria}. The parametric equation of $\widetilde{V}$ is given in {\em(\ref{parametric-equation-solution-cauchy-problem-fourier-equation-b})}

\begin{equation}\label{parametric-equation-solution-cauchy-problem-fourier-equation-b}
\left\{
\begin{array}{l}
 t\circ r\circ f=t\\
  x\circ r\circ f=x\\
  u\circ r\circ f=\sum_{n\ge 0}t^n \frac{1}{n!}\frac{d^{2n}h}{dx^{2n}}(x)\\
  \\
  u_t\circ r\circ f=\sum_{n\ge 1}t^{n-1} \frac{1}{(n-1)!}\frac{d^{2n}h}{dx^{2n}}(x)\\
  \\
  u_x\circ r\circ f=\sum_{n\ge 0}t^n \frac{1}{n!}\frac{d^{2n+1}h}{dx^{2n+1}}(x)\\
  \\
  u_{tt}\circ r\circ f=\sum_{n\ge 2}t^{n-2} \frac{1}{(n-2)!}\frac{d^{2n}h}{dx^{2n}}(x)\\
  \\
  u_{tx}\circ r\circ f=\sum_{n\ge 1}t^{n-1} \frac{1}{(n-1)!}\frac{d^{2n+1}h}{dx^{2n+1}}(x)\\
 \end{array}
\right\}\hskip 2pt x\in[0,1].
\end{equation}

Let us emphasize, that also in this case the solution so obtained is not unique. In fact, similarly to the previous case, where we have considered weak-singular solutions by means of perturbations of Cauchy data for the steady state solution $V_S$, now we have the following weak-singular solution $\widehat{V}\equiv\widetilde{V}\bigcup\widetilde{V}[\mu,a,b]\subset(F)$, where $\widetilde{V}[\mu,a,b]$ is the $2$-dimensional integral manifold identified by the solution of $(F)$ given in {\em(\ref{perturbed-solution-b})}.
\begin{equation}\label{perturbed-solution-b}
\widetilde{u}(t,x|\mu,a,b)=\sum_{n\ge 0}t^n \frac{1}{n!}\frac{d^{2n}h}{dx^{2n}}(x)+e^{-\mu^2t}\left[\frac{a+b-b\cos(\mu)}{\sin(\mu)}\sin(\mu x)+b\cos(\mu x)\right].
\end{equation}
Let us denote by $\widetilde{\Gamma}[\mu,a,b]\subset W_{t=0}$ the space-like curve identified by $\widetilde{u}(0,x|\mu,a,b)$. Then $\widetilde{V}[\mu,a,b]$ passes for $$\widetilde{\Gamma}^{(2)}[\mu,a,b]\subset(F)_{t=0},$$ and the singular solution $\widehat{V}$ passes for the space-like curve $$\widetilde{\Xi}^{(2)}=\widetilde{\Gamma}^{(2)}\bigcup\widetilde{\Gamma}^{(2)}[\mu,a,b]\subset\widehat{V}\subset(F),$$ where $$\widetilde{\Xi}=\widetilde{\Gamma}\bigcup\widetilde{\Gamma}[\mu,a,b]\subset W_{t=0}$$ is the space-like curve containing the fixed Cauchy data, i.e., the curve $\widetilde{\Gamma}$. Furthermore the weak-singular solution $\widehat{V}$ asymptotically converges {\em($t\to+\infty$)} to the regular solution $\widetilde{V}$. Since the singular solution $\widetilde{V}$ depends on the parameters $(\mu,a,b)\in\mathbb{R}^+\times\mathbb{R}^2$, the Cauchy problem identified by the curve $\widetilde{\Gamma}$ has more than one solution.
\end{example}

Let us, now, $N_0, N_1\subset E_k$ be two closed compact $(n-1)$-dimensional
admissible integral manifolds of $E_k$. Then there exists a weak,
(resp. singular, resp. smooth) solution $V\subset E_k$, such that
$\partial V=N_0\DU N_1$, iff $X\equiv N_0\DU
N_1\in[0]\in\Omega_{n-1,w}^{E_k}$, (resp.
$X\in[0]\in\Omega_{n-1,s}^{E_k}$, resp.
$X\in[0]\in\Omega_{n-1}^{E_k}$). On the other hand there exists such
a smooth solution iff $X\in[0]\in\Omega_{n-1}$ and $X$ has zero all
its integral characteristic numbers, i.e., are zero on $X$ all the
conservation laws of $E_k$. Since these last can be identified with
the Hopf algebra
$\mathbf{H}_{n-1}(E_k)\cong\mathbf{H}_{n-1}(E_\infty)$, where
$E_\infty$ is the infinity prolongation of $E_k$, it follows that
the quotient $\mathbf{H}_{n-1}(E_\infty)/\mathbb{R}^{\Omega_{n-1}}$
measures the amount of how the conservation laws of $E_k$ differ
from the crystal conservation laws, identified with the elements of
the Hopf algebra $\mathbb{R}^{\Omega_{n-1}}$.
\end{proof}

\begin{definition}
We define {\em crystal obstruction} of $E_k$ the above quotient of
algebras, and put: $ cry(E_k)\equiv
\mathbf{H}_{n-1}(E_\infty)/\mathbb{R}^{\Omega_{n-1}}$. We call
{\em$0$-crystal PDE} one $E_k\subset J^k_n(W)$ such that
$cry(E_k)=0$.\footnote{An extended $0$-crystal PDE $E_k\subset
J^k_n(W)$ does not necessitate to be a $0$-crystal PDE. In fact
$E_k$ is an extended $0$-crystal PDE if $\Omega_{n-1,w}^{E_k}=0$.
This does not necessarily imply that $\Omega_{n-1}^{E_k}=0$.}
\end{definition}

\begin{cor}
Let $E_k\subset J^k_n(W)$ be a $0$-crystal PDE. Let $N_0, N_1\subset
E_k$ be two closed compact $(n-1)$-dimensional admissible integral
manifolds of $E_k$ such that $X\equiv N_0\DU
N_1\in[0]\in\Omega_{n-1}$. Then there exists a smooth solution
$V\subset E_k$ such that $\partial V=X$.
\end{cor}

\begin{example}[The Ricci-flow equation]\label{ex9}
The Ricci-flow equation
\begin{equation}
F_{ij}\equiv(\partial t.g_{ij})-\kappa R_{ij}=0
\end{equation}

on a Riemannian $n$-dimensional manifold $(M,g)$, can be encoded by
means of a second order differential equation $(RF)\subset J{\it
D}^2(E)\subset J^2_4(E)$ over the following fiber bundle $\pi:E\equiv
\mathbb{R}\times\widetilde{S^0_2M}\to\mathbb{R}\times M$, where
$\widetilde{S^0_2M}\subset S^0_2M$ is the open subbundle of
non-degenerate Riemannian metrics on $M$. In {\em\cite{PRA15}} we
have calculated the integral bordism group of the equation $(RF)\subset J^2_4(E)$. In
particular if $M$ is a $3$-dimensional closed compact smooth simply
connected manifold, we get
$\Omega_{3,w}^{(RF)}\cong\Omega_{3,s}^{(RF)}\cong\mathbb{Z}_2$. Thus
$(RF)$ is not an
 extended $0$-crystal PDE. Taking into account exact commutative
 diagram {\em(\ref{comdiag3})}, we get also
 the short exact sequence
 $0\to K_{3,w}^{(RF)}\to \Omega_3^{(RF)}\to \mathbb{Z}_2\to 0$. Taking into account Example \ref{ex5} we can consider
 $\Omega_3^{(RF)}$ as an extension of a subgroup of the crystallographic group $G(3)=\mathbb{Z}^3\rtimes\mathbb{Z}_2$ or
 $G(2)=\mathbb{Z}^2\rtimes\mathbb{Z}_2$. Therefore the integral
 bordism group of the Ricci-flow equation on $S^3$ is an extended crystal PDE,
 with crystal group $G(2)=\mathbb{Z}^2\rtimes\mathbb{Z}_2=p2$ and crystal dimension $2$.
Furthermore, from the exact commutative
 diagram (\ref{comdiag3}) we get also the following short exact sequence
 $0\to \overline{K}_{3}^{(RF)}\to \Omega_3^{(RF)}\to \Omega_3\to 0$.
 Taking into account that $\Omega_3=0$, we can have $cry(RF)\not=0$.
On the other hand let us
 consider admissible only space-like integral Cauchy manifolds
 satisfying the following conditions: {\em(i)} They are
 diffeomorphic to $S^3$ or to $M$, assumed any smooth $3$-dimensional Riemannian,
 compact, closed, orientable, simply connected manifold; {\em(ii)}
 $M$ is homotopy equivalent to $S^3$.
Such integral manifolds surely
 exist, since we can embed in $E$ both manifolds $M$ and $S^3$ and
 identify these, for example, with space-like smooth admissible integral manifolds of the
subequation $(RF)_t\subset(RF)$, where $(RF)_t=\hat\pi_2^{-1}(t)$,
$t\in\mathbb{R}$. Here $\hat\pi_2$ is the canonical projection $(RF)\to
\mathbb{R}$, induced by $\pi_2:J{\it D}^2(E)\to \mathbb{R}\times M$, i.e., $\hat\pi_2\equiv pr_1\circ\pi_2:J{\it D}^2(E)\to \mathbb{R}$, where $pr_1$ is the canonical projection $pr_1:\mathbb{R}\times M\to\mathbb{R}$. In
 fact, in $(RF)$, $(\partial t.g_{ij})$ is solved with respect to
 $R_{ij}$.  More precisely, starting from a $3$-dimensional compact closed, orientable, simply connected Riemannian manifold $(M,\gamma)$, we can identify a space-like integral Cauchy $3$-dimensional manifold $N_0\subset(RF)_{t=t_0=0}$, diffeomorphic to its projection  $\pi_{2,0}(N_0)\equiv Y_0\subset W_{t=t_0=0}\cong \widetilde{S^0_2M}$, by means of the mapping $f:M\to J{\it D}^2(E)$ defined by the parameter equation in {\em(\ref{3-dimensional-integral-cauchy-manifold-ricci-flow-equation})}.
 \begin{equation}\label{3-dimensional-integral-cauchy-manifold-ricci-flow-equation}
   \left\{
   \begin{array}{l}
     t\circ f=t_0=0\\
     x^k\circ f=x^k\\
     g_{ij}\circ f=\gamma_{ij}\\
      g_{ij,t}\circ f=\kappa R_{ij}(\gamma)\\
       g_{ij,h}\circ f=\gamma_{ij,h}\\
        g_{ij,th}\circ f=\kappa R_{ij}(\gamma)_{,h}\\
         g_{ij,hk}\circ f=\gamma_{ij,hk}\\
          g_{ij,tt}\circ f=\Phi_{ij}(\gamma),\\
   \end{array}
   \right\},\hskip 2pt 1\le i,j,h,k\le 3.
 \end{equation}
where $\{t,x^k,g_{ij},g_{ij,t},g_{ij,h},g_{ij,th},g_{ij,hk},g_{ij,tt}\}_{1\le i,j,h,k\le 3}$ are local coordinates on $J{\it D}^2(E)$, $\Phi_{ij}(\gamma)$, $i,j\in\{,2,3\}$, are known analytic functions of $\gamma_{ij}$ and its derivatives up to fourth order, symmetric in the indexes. In fact, since $(RF)$ is a formally integrable and completely integrable PDE \cite{PRA15}, from its first prolongation we get
$$\{g_{ij,tt}=\kappa R_{ij}(g)_{,t},g_{ij,th}=\kappa R_{ij}(g)_{,h}\}_{1\le i,j,h\le 3},$$ where $R_{ij}(g)_{,t}$ and $R_{ij}(g)_{,h}$ denote the first prolongation of $R_{ij}(g)$, with respect the $t$-variable and $x^h$-variable respectively. Since $R_{ij}(g)=R_{ij}(g_{rs},g_{rs,h},g_{rs,hk})$, $i,j\in\{,2,3\}$, are analytic functions, we get also the following analytic functions
$$R_{ij}(g)_{,t}=K_{ij}(g_{rs},g_{rs,h},g_{rs,hk},g_{rs,t},g_{rs,ht},g_{rs,hkt}),\hskip 2pt 1\le i,j\le 3$$
where $K_{ij}$ are analytic functions of their arguments. Taking into account the expression of $(RF)$ and second order prolongation, we get the initial values ($3$-dimensional Cauchy integral manifold), given in {\em(\ref{3-dimensional-integral-cauchy-manifold-ricci-flow-equation-a})}.
\begin{equation}\label{3-dimensional-integral-cauchy-manifold-ricci-flow-equation-a}
   \left\{
   \begin{array}{l}
   g_{rs,h}|_{t=0}=\gamma_{rs,h}\\
     g_{rs,t}|_{t=0}=\kappa R_{rs}(\gamma)\\
    g_{rs,th}|_{t=0}=\kappa R_{rs,h}|_{t=0}\\
    g_{rs,hk}|_{t=0}=\gamma_{rs,hk}\\
    g_{rs,tt}|_{t=0}=\Phi_{ij}(\gamma)\\
   \end{array}
   \right\}_{1\le r,s,h,k\le 3}\\
    \end{equation}
    where $R_{rs,h}$ is the first prolongation of $R_{rs}$ with respect to the coordinates $x^h$. Taking into account that one has the following functional dependence $R_{rs}=R_{rs}(g_{ij},g_{ij,p},g_{ij,pq})$, we get
    $$R_{rs,h}=(\partial g^{ij}.R_{rs})g_{ij,h}+(\partial g^{ij,p}.R_{rs})g_{ij,ph}+(\partial g^{ij,pq}.R_{rs})g_{ij,pqh}.$$
    As by-product we get that $R_{rs,h}|_{t=0}$ are functions of $\gamma_{ij}$ and their derivative up to third order, that we shortly denote by $R_{rs,h}|_{t=0}$. Furthermore, from the first prolongation of $(RF)$ we get also
    \begin{equation}\label{tt-prolongation-ricci-tensor}
\scalebox{0.8}{$    \left\{\begin{array}{ll}
g_{rs,tt}&=\kappa R_{rs,t}=\kappa{[ (\partial g^{ij}.R_{rs})g_{ij,t}+(\partial g^{ij,p}.R_{rs})g_{ij,pt}+(\partial g^{ij,pq}.R_{rs})g_{ij,pqt}]}\\
\\
&=\kappa{[(\partial g^{ij}.R_{rs})\kappa R_{ij}+(\partial g^{ij,p}.R_{rs})\kappa R_{ij,p}+(\partial g^{ij,pq}.R_{rs})
\kappa R_{ij,pq}]}\\
   \end{array}
   \right\}_{1\le r,s,h,k\le 3}$}
    \end{equation}
   Then taking $t=0$, we get that $g_{rs,tt}|_{t=0}$ are functions that depend on $\gamma_{ij}$ and their derivatives up to fourth order, that we shortly denote by $\Phi_{ij}(\gamma)$.

   \begin{equation}\label{commutative-diagram-ricci-flow-equation}
\scalebox{0.8}{$\xymatrix{&&0&0&\\
  &&J{\it D}^2(E)\ar[u]\ar[r]^{\pi_{2,1}}&J{\it D}^1(E)\ar[u]\ar[r]&0\\
  0\ar[r]&S^0_2(\mathbb{R}\times  M)\otimes\widetilde{S^0_2(M)}\ar[ur]\ar[r]&vTJ{\it D}^2(E)\ar[u]\ar[r]&\pi_{2,1}^*vTJ{\it D}^1(E)\ar[ul]\ar[u]\ar[r]&0\\
  0\ar[r]&g_2\ar[dr]\ar@{^{(}->}[u]\ar[r]&vT(RF)\ar[d]\ar@{^{(}->}[u]\ar[r]&\pi^*_{2,1}vT(RF)_{-1}\ar[dl]\ar[d]\ar@{^{(}->}[u]\ar[r]&0\\
  &&(RF)\ar[d]\ar[r]_{\pi_{2,1}}&(RF)_{-1}\ar[d]\ar[r]&0\\
  &&0&0&}$}
  \end{equation}
 Similarly one can identify the Riemannian manifold $(S^3,\bar\gamma)$, with another space-like integral Cauchy $3$-dimensional manifold $N_1\subset(RF)_{t=t_1\not=0}$, diffeomorphic to its projection  $\pi_{2,0}(N_1)\equiv Y_1\subset W_{t=t_1}\cong \widetilde{S^0_2M}$. Let us emphasize that, fixed the space-like fiber $(RF)_t$, above integral Cauchy manifolds are uniquely identified by the Riemannian manifolds $(M,\gamma)$ and $(S^3,\bar\gamma)$, respectively. For such integral manifolds $N_0$ and $N_1$, necessarily pass solutions of $(RF)$, (hence they are admissible).\footnote{In general a steady state solution of $(RF)$ is not admitted since this should imply that $(M,\gamma)$ is Ricci flat. Furthermore, regular solutions $g_{ij}(t,x^k)=h(t)\bar g_{ij}(x^k)$, with separated time variable from space ones, imply that $(M,\gamma)$ is an Einstein manifold. In fact $R_{ij}(g)=R_{ij}(h\bar g)=R_{ij}(\bar g)$. The Ricci-flow equation becomes $h_t(t)\bar g_{ij}(x^k)=\kappa R_{ij}(\bar g)$. Therefore, must be $h_t(t)=\omega=\kappa R_{ij}(\bar g)/\bar g_{ij}(x^k)$, with $\omega \in\mathbb{R}$. By imposing the initial condition $g_{ij}(0,x^k)=\gamma_{ij}(x^k)$, we get that must be $h(t)=\omega t+1$ and $\bar g_{ij}(x^k)=\gamma_{ij}(x^k)$, hence $R_{ij}(\gamma)=\gamma_{ij}\omega/\kappa$. Vice versa, if the Ricci flow equation is considered only for Einstein manifolds, then solutions with separated variables like above, are admitted. This means that in general, i.e., starting with any $(M,\gamma)$, we cannot assume solutions $g_{ij}(t,x^k)$ with the above separated variables structure, even if these solutions ''arrive'' to $S^3$, that is just an Einstein manifold. The same results can be obtained by considering metrics $g(t,x^k)$, obtained deforming $\gamma(x^k)$, under a space-time flow $\phi_\lambda$ of $M\times\mathbb{R}$, i.e., $\{t\circ\phi_\lambda=t+\lambda, x^k\circ\phi_\lambda=\phi_\lambda^k(t,x^i)\}_{1\le i,k\le 3}$, and with initial condition $\{\phi_0^k=x^k\}_{k=1,2,3}$. In fact if we assume that $\phi_\lambda^r=h(\lambda)\bar\phi^r(x^k)$, then the Ricci flow equation becomes as given in (\ref{ricci-flow-equation-flow-separated}).
 \begin{equation}\label{ricci-flow-equation-flow-separated}
   \frac{\dot h^2(\lambda)}{h^2(\lambda)}=\kappa\frac{\bar\phi^r_i\bar\phi^r_i(R_{rs}(\gamma)\circ\phi_\lambda)}
   {\bar\phi^r_i\bar\phi^r_i(\gamma_{rs}(\gamma)\circ\phi_\lambda)}=\omega\in\mathbb{R}^{+}
   \Rightarrow\left\{
   \begin{array}{ll}
    {\rm(a)}& \hskip 2pt \dot h(\lambda)-\pm\sqrt{\omega}h(\lambda)=0\\
     {\rm(b)}& \hskip 2pt \bar\phi^r_i\bar\phi^r_i[R_{rs}(\gamma)-\frac{\omega}{\kappa}\gamma_{rs}]\circ\phi_\lambda=0\\
   \end{array}
   \right\}
 \end{equation}
 The integration of the equation (\ref{ricci-flow-equation-flow-separated}(a)) gives $h(\lambda)=Ce^{\pm\lambda\sqrt{\omega}}$. Furthermore, if $(M,\gamma)$ is an Einstein manifold, i.e., there exists $\mu\in\mathbb{R}$, such that $R_{rs}(\gamma)=\mu\gamma_{rs}$, then taking $\omega=\mu\kappa$, one has the following solution, uniquely identified by the initial condition, (up to rigid flows): $\phi_\lambda=e^{\pm\lambda\sqrt{\mu\kappa}}x^r$. If $(M,\gamma)$ is not Einstein, or equivalently, assuming $\omega\not=\kappa\mu$, we see that the solutions of equation (\ref{ricci-flow-equation-flow-separated}(b)) are $\bar\phi^r=a^r\in\mathbb{R}$, since the metric $\bar\gamma_{rs}\equiv R_{rs}(\gamma)-\frac{\omega}{\kappa}\gamma_{rs}$ is not degenerate, hence must necessarily be $\bar\phi^r_i=0$. But such a flow does not satisfy initial condition $\{\phi_0^r=x^r\}_{r=1,2,3}$. In conclusion a metric $g_{ij}(t,x^k)$, obtained deforming $\gamma$ with a space-time flow $\phi_\lambda$, where $\phi_\lambda^r=h(\lambda)\bar\phi^r(x^k)$, is a solution of the Ricci flow equation iff $(M,\gamma)$ is Einstein, or Ricci-flat. This last case corresponds to take $\omega=0$ in equation (\ref{ricci-flow-equation-flow-separated}) and has as solution the unique flow, up to rigid ones, $\phi_\lambda^r=x^r$.} We shall prove this by considering Lemma \ref{cauchy-criteria}(2). First let us note that, similarly to the Fourier's heat equation, the Ricci-flow equation is an affine subbundle of the affine bundle $\pi_{2,1}:J{\it D}^2(E)\to J{\it D}^1(E)$, since the fiber $\pi^{-1}_{2,1}(\bar q)\equiv (RF)_{\bar q}$, for $\bar q\in(RF)_{-1}$, is an affine space. In fact, $T_q(RF)_{\bar q}\cong (g_2)_q$, and the vector space $(g_2)_q$ is constant on the fiber $\bar q\in(RF)_{-1}$. The situation is shown by the exact commutative diagram in {\em(\ref{commutative-diagram-ricci-flow-equation})} and by the fact that a vector field $\zeta=\zeta_{ij}^{\alpha\beta}\partial g^{ij}_{\alpha\beta}$ belonging to the symbol must satisfy equation {\em(\ref{equation-condition-for-symbol-ricci-flow-equation})}. (One has used the usual space-time indexes numbering for the coordinates $(t,x^k)_{1\le k\le 3}=(x^\alpha)_{0\le \alpha\le 3}$ on $\mathbb{R}\times M$.) That equation is constant on the points of the fiber $(RF)_{\bar q}$, since it depends only on the derivatives of first order $g_{ij,h}$ and zero order $g_{ij}$.

\begin{equation}\label{equation-condition-for-symbol-ricci-flow-equation}
\zeta_{rs}^{\alpha\beta}(\partial g^{rs}_{\alpha\beta}.F_{ij})(q)=0,\hskip 2pt q\in(RF).
\end{equation}
This follows soon from the expression of the Ricci tensor as a differential polynomial. See equation {\em(\ref{ricci-tensor-differential-polynomial})}.

\begin{equation}\label{ricci-tensor-differential-polynomial}
\scalebox{0.7}{$R_{jn}=g^{rp}R_{rjnp}\hskip 2pt\left\{
\begin{array}{l}
 R_{rjnp}=\frac{1}{2}(g_{rp,jn}+g_{jn,rp}-g_{rn,jp}-g_{jp,rn})\\
 \hskip 5pt+g^{ts}({[jn,s][rp,t]-[jp,s][rn,t]})\\
\end{array}
\right\}\left\{
\begin{array}{l}
 {[ij,k]}=\frac{1}{2}(g_{ik,j}+g_{jk,i}-g_{ij,k})\\
 g^{rp}\equiv{[g_{rp}]}/|g|\\
 |g|\equiv\det(g_{rp})\\
  {[g_{rp}]}:\hskip 2pt\hbox{\rm algebraic complement of $g_{rp}$}.\\
\end{array}
\right.$}
\end{equation}
Since the derivatives of second order $g_{ij,hk}$ appear in {\em(\ref{ricci-tensor-differential-polynomial})} multiplied by zero-order terms only, it follows that equation {\em(\ref{equation-condition-for-symbol-ricci-flow-equation})} is constant on each fiber $(RF)_{\bar q}$. As a by-product it follows that we can apply the retraction method given in Lemma \ref{cauchy-criteria} directly on $(RF)$.

With this respect, let us see more explicitly, and less formally, how we can build, a smooth solution of the Cauchy problem $(RF)$ identified by the $3$-dimensional integral manifold given in {\em(\ref{3-dimensional-integral-cauchy-manifold-ricci-flow-equation-a})}. Since, fixed the Riemannian manifold $(M,\gamma)$, we uniquely identify a $3$-dimensional integral manifold $N_0\subset(RF)$, diffeomorphic to $M$, (resp. $S^3$), then in order to build a solution passing for $N_0$, (resp. $N_1$), it is enough to prove that we are able to identify a time-like vector field, tangent to $(RF)$, transversal to $N_0$ that besides the tangent space $TN_0$ generate an integral planes of $(RF)$. In fact, integral planes $L_{\widehat{q}}\subset (\mathbf{E}_2)_q\cong L_{\widehat{q}}\bigoplus (g_2)_q\subset T_q(RF)$, $\widehat{q}\in(RF)_{+4}$, $q=\pi_{6,2}(\widehat{q})$, are generated by the following {\em horizontal vectors} $\zeta_\alpha(q)=[\partial x_\alpha+\sum_{0\le|\beta|\le 2}g_{ij,\alpha\beta}\partial g^{ij,\beta}]_{\widehat{q}\in(RF)_{+4}}\in T_q(RF)$. {\em($|\beta|$ denotes the multi-index length and $0\le\alpha\le 3$.)} In {\em(\ref{time-like-characteristic-vector-field-of-ricci-flow-equation})} is given a more explicit expression of the horizontal vectors $\zeta_0$.

\begin{equation}\label{time-like-characteristic-vector-field-of-ricci-flow-equation}
\scalebox{0.8}{$\left\{\begin{array}{l}
\zeta_0=\partial t+g_{ij,t}\partial g^{ij}+g_{ij,th}\partial g^{ij,h}+g_{ij,tt}\partial g^{ij,t}
+g_{ij,thk}\partial g^{ij,hk}+g_{ij,tth}\partial g^{ij,th}+g_{ij,ttt}\partial g^{ij,tt}\\
\\
g_{ij,t}=\kappa R_{ij}(\gamma)\hskip 2pt,\hskip 2pt g_{ij,th}=\kappa R_{ij,h}(\gamma) \\
\\
  g_{ij,tt}=\kappa R_{ij,t}(\gamma)=\kappa\left[(\partial g^{rs}.R_{ij})\kappa R_{rs}(\gamma)
  +(\partial g^{rs,p}.R_{ij})\kappa R_{rs,p}(\gamma)+(\partial g^{rs,pq}.R_{ij})\kappa R_{rs,pq}(\gamma)\right]\\
  \\
  g_{ij,thk}=\kappa R_{ij,hk}(\gamma)\hskip 2pt,\hskip 2pt
  g_{ij,tth}=\kappa R_{ij,th}(\gamma) \hskip 2pt,\hskip 2pt
  g_{ij,ttt}=\kappa R_{ij,tt}(\gamma) \\
\end{array}
\right\}.$}
\end{equation}

For $\widehat{q}\in (N_0)^{(4)}=D^6\gamma(M)\subset[(RF)_{+4}]_{t=0}\subset(RF)_{+4}$, at the corresponding point $q\in N_0$, one has $\zeta=\zeta_{(S)}+\zeta_{(T)}\in (\mathbf{E}_2)_q$, with $\zeta_{(S)}=X^k\zeta_k(q)\in T_qN_0$, and $\zeta_{(T)}=X^0\zeta_0(q)$ transverse to $N_0$. Thus the $4$-dimensional integral manifold $V_q$, tangent to such an integral plane, is tangent to $N_0$ too, and has the vector $\zeta_0$ as characteristic time-like vector at $q\in N_0$. By varying $q$ in $N_0$ we generate a $4$-dimensional manifold $V$ that is the {\em envelopment manifold} of the family $\{V_q\}_{q\in N_0}$ of solutions of $(RF)$, formally defined as $V\equiv\bigcup_{q\in N_0}X_{0,q}$, where $X_{0,q}$ is the integral line transversal to $N_0$, starting from $q$, tangent to $\zeta_0(q)$ and future-directed. Thus $V$ is a line bundle over $N_0$, containing $N_0$, with $T_{q}N_0\subset T_{q}V$, $\forall q\in N_0$. This is enough to claim that $V$ is contained in $(RF)$ and it is a time-like integral manifold of $(RF)$, whose tangent space is an horizontal one. In other words, $V$ is not only a viscosity solution, but just a solution of the Cauchy problem identified by $N_0$ (resp. $N_1$).\footnote{Generalized solutions of PDE's, called {\em viscosity solutions} were introduced by Pierre-Louis Lions and Michael Crandall in the paper \cite{LIO-CRAN}. Such mathematical objects do not necessitate to be solutions, but are envelope manifolds of solutions.} (See Fig. \ref{solution-cauchy-problem-envelopment-solution} where the $3$-dimensional manifold $N_0$ is reduced to dimension $2$ (figure in the left-side) or dimension $1$ (figure in the right-side) for graphic necessities.)
\begin{figure}[t]
\centerline{\includegraphics[height=4cm]{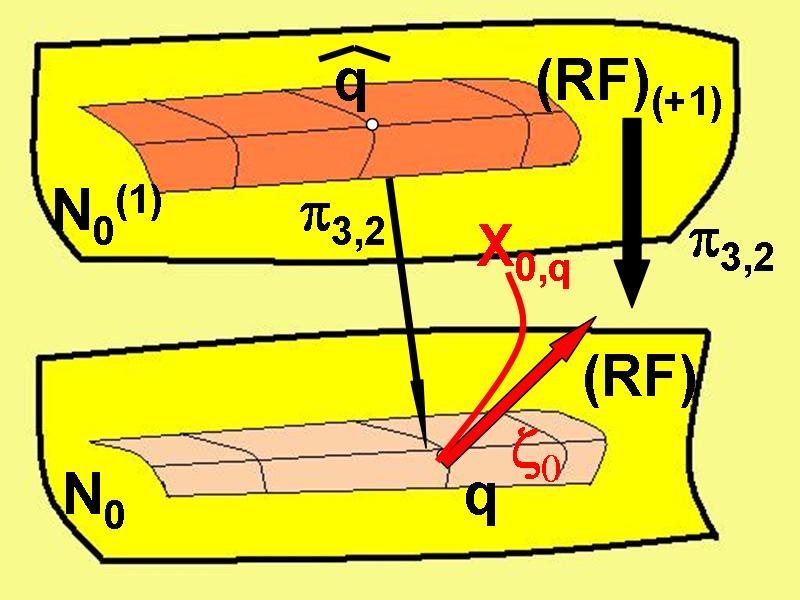}
\includegraphics[height=4cm]{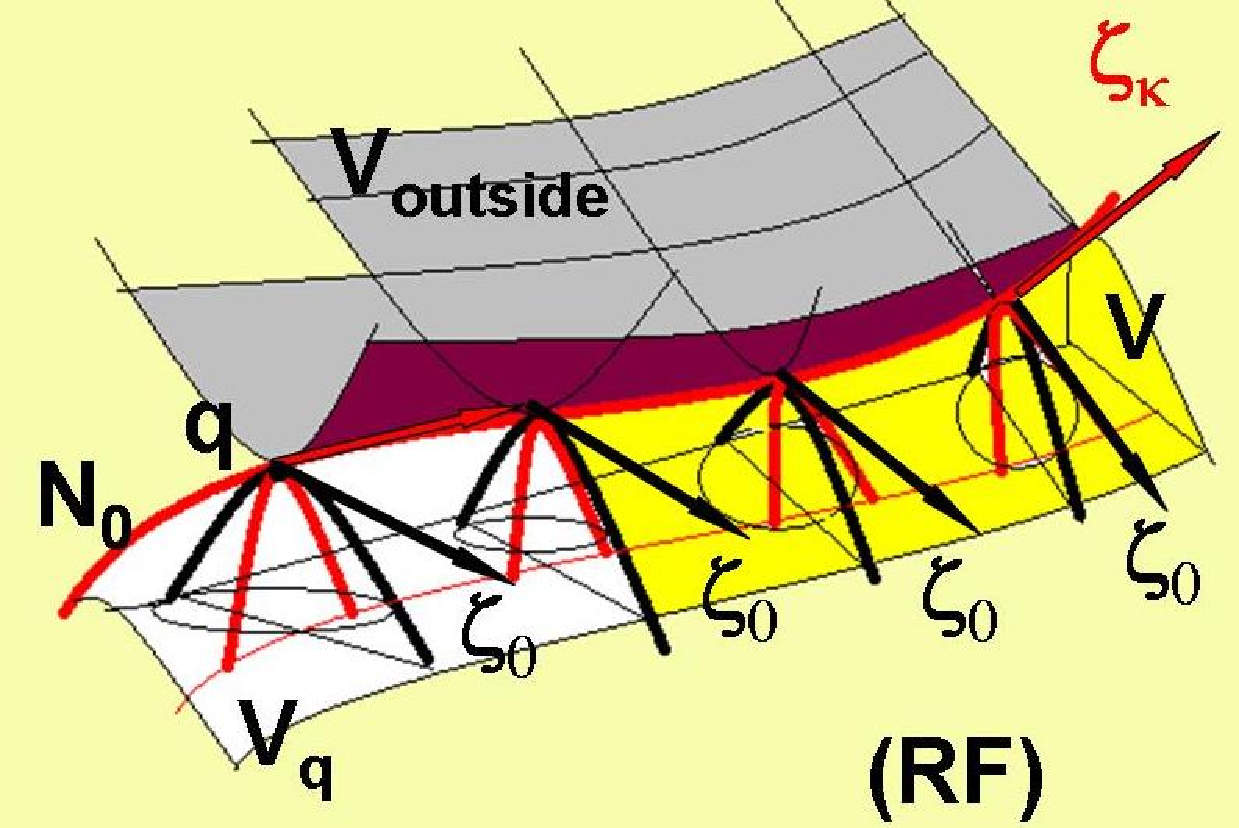}}
\caption{Solution Cauchy problem for $(RF)$, represented by the integral manifold $V$, {\em envelopment manifold}, generated by local solutions represented by manifolds $V_q$ tangent to $N_0$, $\forall q\in N_0$, and identified by means of $4$-dimensional integral planes $L_{\widehat{q}}$, for any $\widehat{q}\in N_0^{(1)}$. $V_{outside}\subset J{\it D}^2(E)$ represents a $4$-dimensional integral manifold, passing for $N_0$, but $V_{outside}\not\subset(RF)$. $V_{outside}$ can be deformed in $V$, taking fixed $N_0$. $(RF)$ and $(RF)_{+1}$ are represented with (yellow) frames, containing respectively $N_0$ and $N_0^{(1)}$.}
\label{solution-cauchy-problem-envelopment-solution}
\end{figure}

Let us emphasize that, on $N_0$, (resp. $N_1$), the components of the characteristic vector field $\zeta_0$ in {\em(\ref{time-like-characteristic-vector-field-of-ricci-flow-equation})} are uniquely characterized by the metric $\gamma_{ij}$ and its derivatives up to sixth order. This proves that in order to characterize a solution of the Cauchy problem, identified by $N_0\subset(RF)$, it is necessary to consider the prolongation $(RF)_{+4}$, of $(RF)$, up to fourth order. Then the $4$-dimensional manifold $V=\bigcup_{t\in[0,\epsilon[}\phi_t(N_0)\subset(RF)$, obtained by means of the local flow $\phi_t$ generated by $\zeta_0=\partial\phi$, is a local integral manifold, solution of the Cauchy problem identified by $N_0$. (So $\zeta_0$ is the characteristic vector field for such a solution.)\footnote{Let us emphasize also that an integral $4$-plane, where the components $g_{ij,\beta}$, $0\le|\beta|\le 3$, satisfy conditions in {\em(\ref{3-dimensional-integral-cauchy-manifold-ricci-flow-equation-a})}, but not all conditions in {\em(\ref{time-like-characteristic-vector-field-of-ricci-flow-equation})}, has the corresponding integral manifold, say $\widetilde{V}$, that passes for $N_0$, but it is not contained in $(RF)$. (This is represented by $V_{outside}$ in Fig. \ref{solution-cauchy-problem-envelopment-solution}.) Let us denote the corresponding metric with $\widetilde{g}$. Therefore the retraction method imposes also to $g_{ij,\beta}$, with $0\le|\beta|\le 3$, to satisfy conditions reported in {\em(\ref{time-like-characteristic-vector-field-of-ricci-flow-equation})}. The relation between $\widetilde{V}\equiv V_{outside}\subset J{\it D}^2(E)$, and $V\subset (RF)$, can be realized with a deformation, $g_{\lambda}$, connecting the corresponding metrics $\widetilde{g}$ and $g$. More precisely, $g_{\lambda}=\widetilde{g}+\lambda[g-\widetilde{g}]$, $\lambda\in[0,1]$. The integral manifolds $V_{\lambda}\subset J{\it D}^2(E)$, generated by $g_{\lambda}$, are not contained in $(RF)$ for any $\lambda\in[0,1]$, but all pass for $N_0$. In fact, since $(g_{\lambda})_{ij,\beta}=\widetilde{g}_{ij,\beta}+\lambda[g_{ij,\beta}-\widetilde{g}_{ij,\beta}]$, $|\beta|\ge 0$, we get that $(g_{\lambda})_{ij,\beta}|_{N_0}=\widetilde{g}_{ij,\beta}|_{N_0}=g_{ij,\beta}|_{N_0}$, for $0\le|\beta|\le 2$, but $(g_{\lambda})_{ij,\beta}|_{N_0}$ do not satisfy conditions in {\em(\ref{time-like-characteristic-vector-field-of-ricci-flow-equation})} for $|\beta|=3$. This is just the meaning of the retraction method considered in the proof of Lemma \ref{cauchy-criteria}.} Let us remark that such a solution does not necessitate to be smooth. In fact in this process we have used only the prolongation of $(RF)$ up to fourth-order, i.e., we have considered only the projections given in {\em(\ref{fourth-order-prolongation-tower-ricci-flow-equation})}.
\begin{equation}\label{fourth-order-prolongation-tower-ricci-flow-equation}
\xymatrix{(RF)_{+4}\ar[r]&(RF)_{+3}\ar[r]&(RF)_{+2}\ar[r]&(RF)_{+1}\ar[r]&(RF)}.
\end{equation}
However, we can obtain $V$ as a continuous manifold. In fact, let $g_{ij}(x^\alpha)$ denote a local solution identified by $L_{\widehat{q}}=T_qV_q\subset T_q(RF)$, and $T_qN_0\subset L_{\widehat{q}}$. Then the time-like integral line $X_{0,q}\subset V_q$, tangent to $\zeta_0$, is represented by a local curve $\xi_q:[0,\epsilon_q[\subset\mathbb{R}\to(RF)$, given by the parametric equation (\ref{parametric-equation-time-like-curve-local-solution-ricci-flow}) in $J{\it D}^2(E)$.
\begin{equation}\label{parametric-equation-time-like-curve-local-solution-ricci-flow}
\left\{
\begin{array}{l}
 t\circ\xi=t\\
 x^k\circ\xi=x^k(q)\\
 g_{ij}\circ\xi=g_{ij}(t,x^k(q))\\
 g_{ij,t}\circ\xi=\kappa R_{ij}(g(t,x^k))|_{x^k=x^k(q)}\\
 g_{ij,h}\circ\xi=g_{ij,h}(t,x^k)|_{x^k=x^k(q)}\\
 g_{ij,hk}\circ\xi=g_{ij,hk}(t,x^k)|_{x^k=x^k(q)}\\
 g_{ij,th}\circ\xi=\kappa R_{ij,h}(g(t,x^k))|_{x^k=x^k(q)}\\
 g_{ij,tt}\circ\xi=\kappa^2\sum_{0\le|p|\le 2}[(\partial g^{rs,p}.R_{ij})R_{rs,p}(g(t,x^k))]|_{x^k=x^k(q)}\\
\end{array}
\right.
\end{equation}
Since $\zeta_0$ continuously changes on $N_0$, it follows that curves $X_{0,q}$, $q\in N_0$, continuously change with $q\in N_0$, if the solutions $g_{ij}(t,x^k)$ continuously change with $q\in N_0$ too. This is surely realized by the fact that points $q'\in N_0$, near to $q$, come from points $\widehat{q}'\in N_0^{(1)}$, near to $\widehat{q}$, if $\pi_{3,2}(\widehat{q})=q$. In fact, $N_0\cong N_0^{(1)}$ and $\pi_{3,2}(N_0^{(1)})=N_0$, thus $\pi_{3,2}|_{N_0^{(1)}}:N_0^{(1)}\to N_0$ is necessarily a continuous mapping. Moreover, we can continuously transform any (local) solution $g_{ij}(t,x^k)$ of $(RF)$ into other ones by means of space-like (local) $1$-parameter group of diffeomorphisms $\phi_\lambda$, $\lambda\in[0,\epsilon[\subset\mathbb{R}$ of $M$. In fact, such transformations are locally represented by the following functions $\{t\circ\phi_\lambda=t,x^k\circ\phi_\lambda=\phi_\lambda^k(x^i)\}$. As a by-product it follows that $(RF)$ is invariant under such transformations. In fact $(\phi_\lambda^*g)_{ij,t}-\kappa R_{ij}(\phi_\lambda^*g)=\phi_\lambda^*[g_{ij,t}-\kappa R_{ij}(g)]=0$. Thus we can continuously transform an integral manifold $V_q$, into other solutions along the space-like coordinate lines of $N_0$, passing for $q$ and identify, in the points $q'\in N_0$, the time-like curves $X_{0,q'}$ that result continuously transformed of $X_{0,q}$. In this way
$V=\bigcup_{q\in N_0}X_{0,q}$ is a smooth manifold in a suitable tubular neighborhood of $N_0\times[0,\epsilon[$, that is the integral manifold of a metric $g_{ij}(t,x^k)$,  of class $C^3$, solution of the Ricci flow equation $(RF)$. This proves that the envelopment manifold $V$ is more regular than a viscosity solution.

In order to obtain a smooth solution of the Cauchy problem given by the integral manifold $N_0$ it is necessary to repeat the above process on the infinity prolongation $(RF)_{+\infty}\subset J{\it D}^{\infty}(E)$. In fact, also on $(RF)_{+\infty}$, the $3$-dimensional integral manifold $N_0\cong D^\infty\gamma(M)\equiv N_0^{(\infty)}\subset (RF)_{+\infty}$ is uniquely identified by $(M,\gamma)$, and identifies also a unique transversal time-like characteristic vector field $\zeta_0^{(\infty)}$, tangent to $(RF)_{+\infty}$. More precisely, $\zeta_0^{(\infty)}=\zeta_0+\sum_{|\alpha|>2}g_{ij,t\alpha}\partial g^{ij,\alpha}$, with $g_{ij,t\alpha}=\kappa R_{ij,\alpha}(\gamma)$, where $\alpha=0,1,2,3$.

The situation is resumed in the commutative diagram {\em(\ref{commutative-diagram-characteristic-vector-field-smooth-solution-cauchy-problem-ricci-flow-equation})}.
\begin{equation}\label{commutative-diagram-characteristic-vector-field-smooth-solution-cauchy-problem-ricci-flow-equation}
\scalebox{0.7}{$\xymatrix{&T(RF)_{+\infty}&&TN_0^{(\infty)}\ar@{^{(}->}[ll]&\\
T(RF)_{+\infty}\ar[d]&(RF)_{+\infty}\ar[l]^{\zeta_0^{(\infty)}}\ar[u]^{\zeta_0^{(\infty)}}\ar@{^{(}->}[l]\ar[d]&[(RF)_{+\infty}]_{t=0}\ar@{^{(}->}[l]\ar[d]&
N_0^{(\infty)}\ar@{=}[d]_{\wr}\ar[u]^{0}\ar@{^{(}->}[l]\ar@{^{(}->}[r]&J{\it D}^{\infty}(\widetilde{S^0_2M})\ar[d]\\
\vdots\ar[d]&\vdots\ar[d]&\vdots\ar[d]&\vdots\ar@{=}[d]_{\wr}&\vdots\ar[d]\\
T(RF)_{+r}\ar[d]&(RF)_{+r}\ar[l]^{\zeta_0^{(r)}}\ar@{^{(}->}[l]\ar[d]&[(RF)_{+r}]_{t=0}\ar@{^{(}->}[l]\ar[d]&N_0^{(r)}\ar@{=}[d]_{\wr}\ar@{^{(}->}[l]\ar@{^{(}->}[r]&J{\it D}^{2+2r}(\widetilde{S^0_2M})\ar[d]\\
\vdots\ar[d]&\vdots\ar[d]&\vdots\ar[d]&\vdots\ar@{=}[d]_{\wr}&\vdots\ar[d]\\
T(RF)&(RF)\ar[l]^{\zeta_0}\ar@{^{(}->}[l]\ar[dd]&[(RF)]_{t=0}\ar@{^{(}->}[l]\ar[dd]&N_0\ar@{=}[dd]_{\wr}\ar@{^{(}->}[l]\ar@{^{(}->}[r]&J{\it D}^{4}(\widetilde{S^0_2M})\ar[d]\\
&&&&\widetilde{S^0_2M}\ar[d]\\
&\mathbb{R}\times M&M\ar@{=}[r]\ar@{^{(}->}[l]&M\ar@{=}[r]&M\ar@/^1pc/[u]^{\gamma}\ar@/^3pc/[uu]^(0.7){D^4\gamma}\ar@/_3pc/[uuuu]_(0.8){D^{2+2r}\gamma}
\ar@/_5pc/[uuuuuu]_(0.8){D^{\infty}\gamma}}$}
\end{equation}
Similarly to what made in the Fourier's heat equation, (see Example \ref{fourier-heat-equation}), we can prove that to the above solution one can associate weak-singular ones by using perturbations of the initial Cauchy data. In fact, if $V\subset(RF)$ is a regular solution passing for the integral manifold $N_0$, we can consider perturbations like solutions $\nu:\mathbb{R}\times M\to E$, of the linearized equation $(RF)[V]\subset J{\it D}^2(E)$. Since also $(RF)[V]$ is formally integrable and completely integrable, in neighborhoods of $N_0$ there exist perturbations. These deform $V$ {\em(background solution)} giving some new solutions $\widetilde{V}\subset(RF)$. Then $\widehat{V}\equiv V\bigcup\widetilde{V}\subset(RF)$ is a weak-singular solutions of the type just considered in Example \ref{fourier-heat-equation} for the Fourier's heat equation.\footnote{This generalizes a previous result by Hamilton \cite{HAMIL2}, and after separately by De Turk \cite{DETU} and Chow \& Knopp \cite{CH-KN}, that proved existence and uniqueness of non-singular solution for Cauchy problem in some Ricci flow equation. Let us also emphasize that our approach to find solutions for Cauchy problems, works also when $N_0\subset (RF)$ is diffeomorphic to a $3$-dimensional space-like submanifold of $W$, that is not necessarily representable by a section of
$E_{t=0}\cong\widetilde{S^0_2M}\to M$.}

Now, for any of two of such integral manifolds, $N_0$ and $N_1$, we can find a
 smooth solution $V$ bording them, $V=N_0\DU N_1$ iff their integral characteristic
 numbers are equal, i.e. all the conservation laws of $(RF)$ valued
 on them give equal numbers.
By the way, under our assumptions we can consider $N_0$ and $N_1$
homotopy equivalent. Let $f:N_1\to N_0$ be such an homotopy
equivalence. Let $\omega$ be any conservation law for $(RF)$. Then
one has equation {\em(\ref{conservation-laws-relations})}.

\begin{equation}\label{conservation-laws-relations}
<[N_0],[\omega]>=\int_{N_0}\omega=\int_{N_1}f^*\omega=<[N_1],[f^*\omega]>=<[N_1],[\omega]>=\int_{N_1}\omega.
\end{equation}

So any possible integral characteristic number of $N_0$ must
coincide with ones with $N_1$ and vice versa. Thus we can say, that
with this meaning of admissibility ({\em full admissibility hypothesis}) on the Cauchy integral manifolds,
one has $cry(RF)=0$, i.e., $(RF)$ becomes a $0$-crystal. Therefore there
are not obstructions on the existence of smooth solutions $V$ of
$(RF)$ bording $N_0$ and $N_1$, $\partial V=N_0\DU N_1$, i.e.,
solutions without singular points.

This has as a by-product that $M$ and $S^3$ are homeomorphic
manifolds. Hence the Poincar\'e conjecture is proved. With this
respect we can say that this proof of the  Poincar\'e conjecture is
related to the fact that under suitable conditions of admissibility
for the Cauchy integral manifolds, the Ricci-flow equation becomes a
$0$-crystal PDE.\footnote{The proof of the Poincar\'e conjecture
given here refers to the Ricci-flow equation, according to some
ideas pioneered by Hamilton \cite{HAMIL1, HAMIL2, HAMIL3, HAMIL4,
HAMIL5}, and followed also by Perelman \cite{PER1, PER2}. However
the arguments used here are completely different from ones used by
Hamilton and Perelman. (For general informations on the relations
between Poincar\'e conjecture and Ricci-flow equation, see, e.g.,
Refs.\cite{ANDER, C-C-G-G-I-I-K-L-L-N1, C-C-G-G-I-I-K-L-L-N2} and
papers quoted there.) Here we used our general PDE's
algebraic-topological theory, previously developed in some works.
Compare also with our previous proof given in \cite{AG-PRA1,
AG-PRA2}, where, instead was not yet introduced the relation between
PDE's and crystallographic groups.}
\end{example}

\begin{example}[The d'Alembert equation $\frac{\partial^2\log f}{\partial x\partial y}=0$ on the
$2$-dimensional torus]\label{ex10} The  d'Alembert equation on a
$2$-dimensional manifold $M$ can be encoded by a second-order
differential equation
\begin{equation}
(d'A)_2\subset J{\it D}^2(W):\left\{uu_{xy}-u_xu_y=0\right\}
\end{equation}
with $W\equiv M\times\mathbb{R}$. In
{\em\cite{PRA15}} we have calculated the integral bordism groups of
such an equation. In particular for $M=T^2$, the $2$-dimensional
torus, one has
$\Omega_1^{(d'A)_2}\cong\mathbb{Z}_2\bigoplus\mathbb{Z}_2$. Taking
into account Example \ref{ex6} we see that $\Omega_1^{(d'A)_2}$ is
isomorphic to the crystal group $G(2)=\mathbb{Z}^2\rtimes D_4=p4m$.
Therefore, $(d'A)_2$ on $T^2$ is an extended crystal PDE, with
crystal dimension $2$.
\end{example}

\begin{example}[The Tricomi equation on $2$-dimensional manifolds]\label{ex11}
In {\em\cite{PRA15}} we have considered the integral bordism groups
of the {\em Tricomi equation} $(T)$:
\begin{equation}
(T)\subset J{\it D}^2(W):\left\{u_{yy}-yu_{xx}=0\right\}
\end{equation}

defined on a $2$-dimensional manifold $M$, i.e. with $W\equiv M\times\mathbb{R}$. For example, on the
$2$-dimensional torus $T^2$ one has $\Omega_1^{(T)}\cong{\mathbb
Z}_2\oplus{\mathbb Z}_2$. Furthermore on ${\mathbb R}P^2$ we obtains
$ \Omega_1^{(T)}\cong{\mathbb Z}_2$. Thus the Tricomi equation on
$T^2$, (resp. $S^2$), is an extended crystal PDE with crystal group
$p4m$, (resp. $p2$), and crystal dimension $2$, (resp. $2$).
\end{example}

\begin{example}[The Navier-Stokes PDE]\label{ex12} The non isothermal
Navier-Stokes equation can be encoded by a second order PDE
$\widehat{(NS)}\subset J{\it D}^2(W)$ on a $9$-dimensional affine
fiber bundle $\pi:W\to M$ on the Galilean space-time $M$. In Tab. \ref{completely-integrable-navier-stokes-equation}
is reported its polynomial differential structure.
\begin{table}[t]\centering
\caption{Completely integrable Navier-Stokes equation: $\widehat{(NS)}\subset J{\it D}^2(W)$ defined by differential polynomials.}
\label{completely-integrable-navier-stokes-equation}
\begin{tabular}{|l|}
\hline ${\rm(A)}: F^0\equiv\dot x^k
G^j_{jk}+\dot x^i_s\delta^s_i=0$ \hfill{\rsmall($\scriptstyle F^0\in A[\dot x^k,\dot x^i_s]$)}\\
\hfill\hbox{\bsmall(continuity equation)}\\
\hline ${\rm(B)}: F^0_\alpha\equiv\dot x^k (\partial
x_\alpha.G^j_{jk})+\dot x^k_\alpha G^j_{jk}+ \dot x^i_{s\alpha}\delta ^s_i=0$ \hfill{\rsmall($\scriptstyle F^0_\alpha\in A[\dot x^k,\dot x^i_\alpha,\dot x^i_{s\alpha}]$)}\\
\hfill\hbox{\bsmall(first prolonged continuity equation)}\\
\hline ${\rm(C)}:
                         F^j\equiv\dot x^s R^j_s+\dot x^s\dot x^i\rho G^j_{is}+\dot x^s\dot x^j_s
                         \rho+\rho\dot x^j_0+\dot x^k_s S^{js}_k+\dot x^j_{is}T^{is}$\\
\hskip 3cm $+p_i g^{ij}+\rho(\partial x_i.f)g^{ij}=0$ \hfill{\rsmall($\scriptstyle F^j\in A[\dot x^s,\dot x^i_s,\dot x^j_{is},p_i]$)}\\
\hfill\hbox{\bsmall(motion equation)}\\
\hline ${\rm(D)}: F^4\equiv\theta_0 \rho C_p+\rho C_p\dot x^k
\theta_k+\theta_{is} \overline{E}^{is}+\dot x^k\dot x^pW_{kp}+\dot
x^k\dot x^s_p {\overline W}^p_{ks}+ \dot
x^k_i\dot x^s_p Y^{ip}_{ks}=0$\\
 \hskip 4cm{\rsmall($\scriptstyle F^4\in A[\dot x^k,\dot x^s_p,\theta_0,\theta_k,\theta_{is}]$)}\hfill\hbox{\bsmall(energy equation)}\\
\hline \hline
\multicolumn{1}{|c|}{\bsmall Functions belonging to {\boldmath $\scriptstyle A\subset\mathbb{R}[[x^1,x^2,x^3]]$}}\\
\hline \multicolumn{1}{|l|}{$\scriptstyle R^j_s\equiv\chi[(\partial
x_p.
G^j_{is})+G^j_{pq} G^q_{is}-G^q_{pi}G^j_{qs}]g^{pi}$}\\
\multicolumn{1}{|l|}{$\scriptstyle S^{js}_k\equiv\chi[2G^j_{ik}g^{si}-G^s_{qi}\delta^j_kg^{qi}]$}\\
\multicolumn{1}{|l|}{$\scriptstyle T^{is}\equiv \chi g^{is}=T^{si}$}\\
\multicolumn{1}{|l|}{$\scriptstyle\overline{E}^{is}\equiv-\nu
g^{is}=\overline{E}^{si}$}\\
\multicolumn{1}{|l|}{$\scriptstyle W_{kp}\equiv-\chi
G^a_{jk}G^b_{sp}(g_{ba}g^{sj}+\delta^j_b\delta^s_a)=-2\chi G^s_{bk}G^b_{sp}=W_{pk}$}\\
\multicolumn{1}{|l|}{$\scriptstyle{\overline W}^p_{ks}\equiv -2\chi[
G^j_{ik}g_{js}g^{ip}+ G^p_{sk}]=-4\chi G^p_{sk}={\overline
W}^p_{sk}$}\\
\multicolumn{1}{|l|}{$\scriptstyle Y^{ip}_{ks}\equiv
-\chi[g_{ks}g^{ip}+\delta^p_k\delta^i_s]=Y^{pi}_{sk}$}\\
\hline \multicolumn{1}{l}{\rsmall $\scriptstyle 0\le\alpha\leq 3,
\hskip 2pt 1\leq i,j,k,p,s\leq 3. \hskip 2pt (\partial
x_0.G^j_{jk})=0$.}\\
\end{tabular}
\end{table}

There
$\{x^\alpha,\dot x^k,p,\theta\}$ are fibered coordinates on $W$,
adapted to the inertial frame, and $G^k_{ij}$ are the canonical
connection symbols on $M$. Furthermore
$A\subset\mathbb{R}[[x^1,x^2,x^3]]$ is the algebra of real valued
analytic functions of
$(x^k)$.\footnote{$\mathbb{R}[[x^1,\dots,x^n]]$ denotes the algebra
of formal series $\sum_{i_1,\dots,i_n}a_{i_1\cdots
i_n}(x^1)^{i_1}\cdots(x^n)^{i_n}$, with $a_{i_1\cdots
i_n}\in\mathbb{R}$. Real analytic functions in the indeterminates
$(x^1,\dots,x^n)$, are identified with above formal series having
non-zero converging radius. Thus real analytic functions belong to a
subalgebra of $\mathbb{R}[[x^1,\dots,x^n]]$. This last can be also
called the {\em algebra of real formal analytic functions} in the
indeterminates $(x^1,\dots,x^n)$.} We have proved in
Refs.{\em\cite{PRA10, PRA17}} that the singular integral bordism
groups of such an equation are trivial, i.e.,
$\Omega^{\widehat{(NS)}}_{3,s}\cong\Omega^{\widehat{(NS)}}_{3,w}\cong
0$. Furthermore, with respect to the notation used in diagram
{\em(\ref{comdiag2})}, one has that for the integral bordism group
for smooth solutions: $\Omega^{\widehat{(NS)}}_{3}\cong
K^{\widehat{(NS)}}_{3,w}$. Thus we can conclude that the
Navier-Stokes equation is an extended $0$-crystal PDE, but not a
$0$-crystal, i.e. $cry\widehat{(NS)}\not=0$. Note that if we
consider admissible only all the Cauchy integral manifolds
$X\subset\widehat{(NS)}$ such that all their integral characteristic
numbers are zero, {\em(full admissibility hypothesis)}, it follows
that $cry\widehat{(NS)}=0$. So, under this condition,
$\widehat{(NS)}$ becomes a $0$-crystal PDE.

\begin{figure}[h]
\centerline{\includegraphics[height=4.5cm]{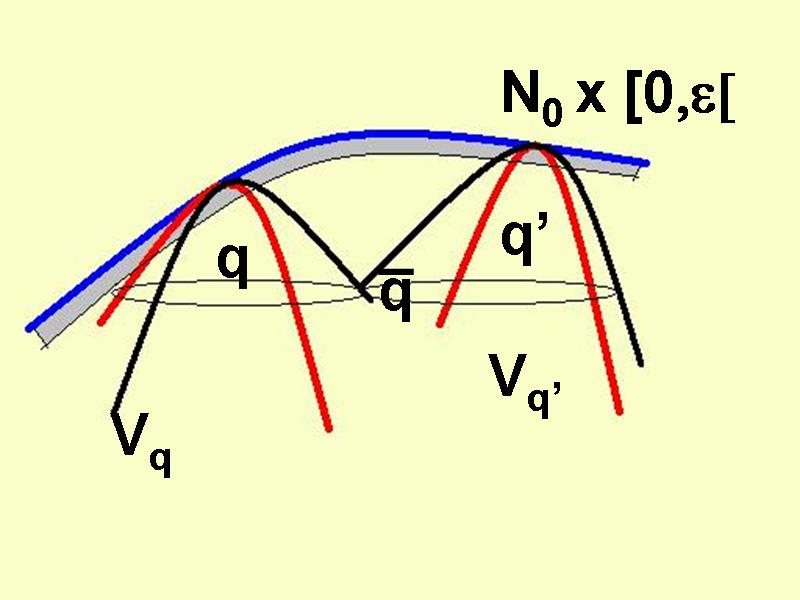}\includegraphics[height=4.5cm]{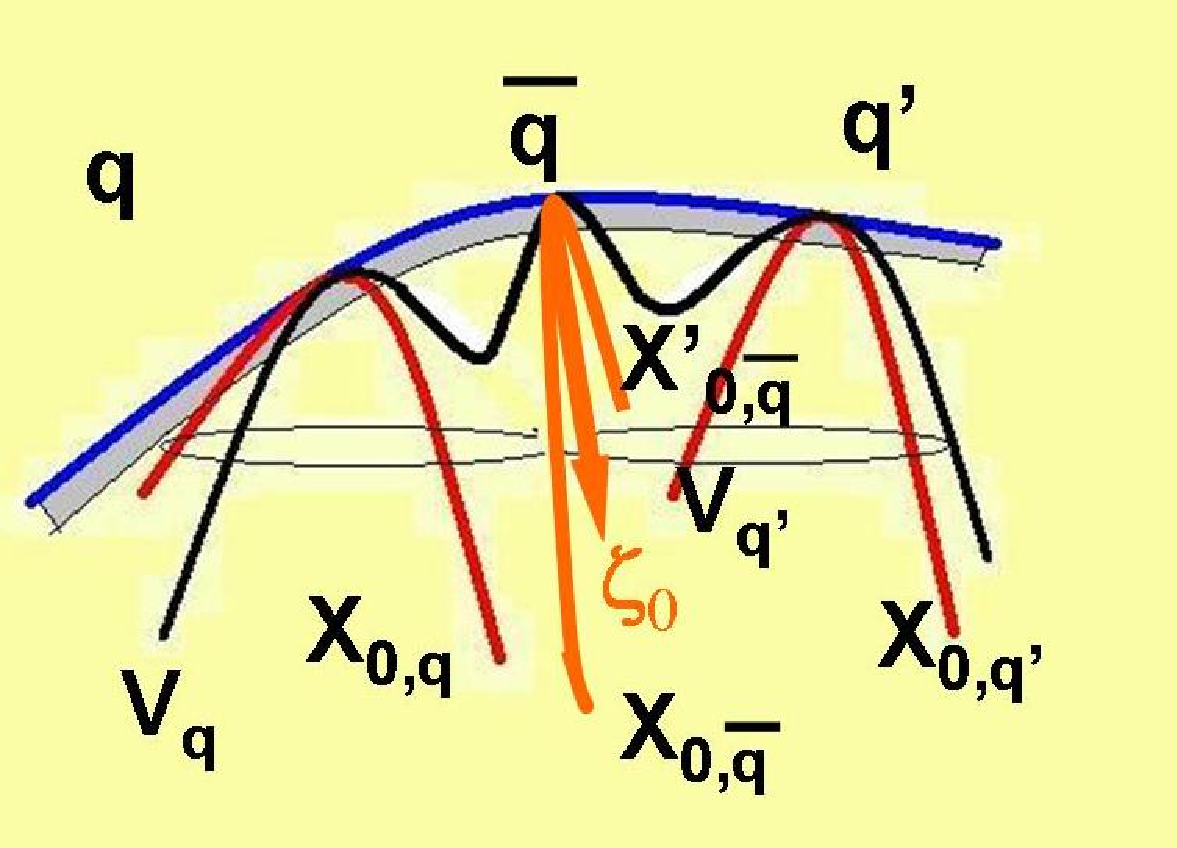}}
\caption{Construction of envelopment solution $V$ for Cauchy problem $N_0^{(\infty)}\subset\widehat{(NS)}_{+\infty}$, on the infinity prolongation of $\widehat{(NS)}_{+\infty}$ of $\widehat{(NS)}$. One has $V=\bigcup_{q\in N_0^{(\infty)}}X_{0,q}$, with $X_{0,q}$ the unique time-like curve passing for $q$, tangent to the vector $\zeta_0(q)$ and contained in the unique analytic solution passing for $q$.}
\label{envelopment-cauchy-solution-at-infinity-prolongation}
\end{figure}

Let us emphasize that, similarly to the Ricci flow equation, we can identify $3$-dimensional space-like smooth integral manifold $N_0\subset\widehat{(NS)}_{t}$, such that $N_0\cong M_t$, via the canonical projection $\pi_{2}:J{\it D}^2(W)\to M$, with some smooth space-like section $s_t:M_t\subset M\to W_t\subset W$, of the configuration bundle $\pi:W\to M$. In fact, all the coordinates in $\widehat{(NS)}$, containing time-derivatives, can be expressed by means of the other derivatives containing only space coordinates. (Warn for equation $\widehat{(NS)}$! Even if there are not restrictions on the time derivatives of pressure, between space derivatives of pressure arise some constraints, as well as there are constraints between space derivatives of velocity components, hence sections $s_t$ that identify above considered integral manifolds, are not arbitrary ones. For details see below.) Then we can solve the corresponding Cauchy problem applying Lemma \ref{cauchy-criteria}, similarly to what made in the Ricci flow equation. It is useful to remark that in order to build envelopment solutions $V=\bigcup_{q\in N_0}X_{0,q}$, it does not necessitate to handle with PDE's that admit any space-like symmetry. This of course does not happen for any smooth boundary value problem in $\widehat{(NS)}$. With this respect, it is useful also to underline that the popular request to maximize entropy cannot be an enough criterion to realize a smooth envelopment solution. (For complementary results on variational problems constrained by the Navier-Stokes equation see \cite{PRA17}.) The existence of such a smooth envelopment manifold, can be proved by working on $\widehat{(NS)}_{+\infty}$. In fact, let $q,q'\in N_0^{(\infty)}\subset\widehat{(NS)}_{+\infty}$, and let $V_q$ and $V_{q'}$ be two smooth solutions passing for the initial conditions $q$ and $q'$ respectively. We claim that their time-like integral curves $X_{0,q}$ and $X_{0,q'}$ cannot intersect for suitable short times, i.e., $t\in[0,\epsilon[$, if $q\not=q'$. Really if $\bar q\in V_{q}\bigcap V_{q'}\not=\emptyset$, then $T_{\bar q}V_q=T_{\bar q}V_{q'}=(\mathbf{E}_\infty)_{\bar q}$. This means that in such a point $\bar q$, $V_q$ and $V_{q'}$ must have a contact of infinity order with the Navier-Stokes equation and between them. We can assume that $\bar q$ is outside a suitable tubular neighborhood $N_0^{(\infty)}\times[0,\epsilon[$ of $N_0^{(\infty)}$, otherwise we should admit that $\bar q\in N_0^{(\infty)}$. In this last case $V_q$ and $V_{q'}$ should have in $\bar q$ time-like curves $X_{0,\bar q}$ and $X'_{0,\bar q}$, transversal to $N_0^{(\infty)}$, and with a common tangent vector $\zeta_0(\bar q)$. Furthermore , $V_q$ and $V_{q'}$ should be tangent to $N_{0}^{(\infty)}$ at $\bar q$: $T_{\bar q}V_q=T_{\bar q}V_{q'}=(\mathbf{E}_\infty)_{\bar q}\supset T_{\bar q}N_0^{(\infty)}$. Let us assume that in a suitable neighborhood $N_0^{(\infty)}\times[0,\epsilon[$, the manifold $V_q$ identifies two separated pieces, say and $V_{q|1}$ and  $V_{q|2}$, tangent to $N_0^{(\infty)}$ at $q$ and $\bar q$, respectively. In other words $V_q\bigcap N_0^{(\infty)}\times[0,\epsilon[=V_{q|1}\bigcup V_{q|2}$. (If this condition is not satisfied, then $V_q$ is necessarily a smooth solution of the Cauchy problem at least in a subset $Y\subset N_0^{(\infty)}$. Then we can write $V=V_1\bigcup V_2\equiv\{\bigcup_{q\in Y}X_{0,q}\}\bigcup\{\bigcup_{q\in Z}X_{0,q}\}$, $Z\equiv N_0^{(\infty)}\setminus Y$, and we can continue to similarly discuss about the $V_2$-part of $V$.) Then, for the time-like coordinate lines $X_{0,q|1}$ and $X_{0,q|2}$, one has $X_{0,q|1}\bigcap X_{0,q|2}=\emptyset$ in $N_0^{(\infty)}\times[0,\epsilon[$. The same circumstance can be verified between $V_{q'}$ and $V_{\bar q}$, eventually by reducing $\epsilon$. This proves that taking $\epsilon$ little enough, we can consider $V\equiv\bigcup_{q\in N_0^{(\infty)}}X_{0,q}\subset \widehat{(NS)}_{+\infty}$, a fiber bundle over $N_0^{(\infty)}$ that is a $4$-dimensional integral manifold representing a local solution of the Cauchy problem identified by the smooth $3$-dimensional integral manifold $N_0\subset \widehat{(NS)}$. (See Fig. \ref{envelopment-cauchy-solution-at-infinity-prolongation}.) On $V$ we can recognize a natural smooth fiber bundle structure. In fact, since $\widehat{(NS)}$ is an analytic equation and any point $q\in N_0^{(\infty)}$ identifies an analytic solution in a suitable neighborhood $U$ of $p\equiv\pi_{\infty}(q)\in M$, it follows that if two such solutions $s$ and $s'$ are defined at $q\in N_0^{(\infty)}$ they should coincide in a suitable neighborhood of $p$, denote again it
by $U$. Then we can assume that $s$ and $s'$ coincide also in $\Omega\equiv\pi_{\infty}^{-1}(U)\bigcap N_0^{(\infty)}$. Therefore, we can cover $N_0^{(\infty)}$ by means of a covering set $\{\Omega_\alpha\}_{\alpha\in J}$, where each $\Omega_\alpha$ is such that if $q,q'\in \Omega_\alpha$, the corresponding analytic solutions $s$ and $s'$ are well defined in $U_\alpha\subset M$. In this way each time-like curve $X_{0,q}$ is uniquely
identified for $q\in \Omega_\alpha$. Then we can define smooth functions ({\em transition functions}) $\psi_{\alpha\beta}:\Omega_\alpha\bigcap\Omega_\beta\to\mathbb{R}^\bullet\equiv\mathbb{R}\setminus\{0\}$, by $\tau_{\alpha}=\psi_{\alpha\beta}\tau_{\beta}$, where $\tau_{\alpha}$ are nowhere vanishing sections of $V\to N_0^{(\infty)}$ on $\Omega_\alpha\subset N_0^{(\infty)}$. These functions satisfy the following three conditions: (a) $\psi_{\alpha\alpha}=1$; (b) $\psi_{\alpha\beta}=\psi_{\beta\alpha}$; (c) $\psi_{\alpha\beta}\psi_{\beta\gamma}\psi_{\gamma\alpha}=1$ on $\Omega_\alpha\bigcap\Omega_\beta\bigcap\Omega_\gamma$ ({\em cocycle condition}). This is enough to claim that the line bundle $V\to N_0^{(\infty)}$ has a smooth structure. Let us emphasize that Lemma \ref{cauchy-criteria} can be applied here also to $3$-dimensional space-like smooth Cauchy integral manifolds $N\subset\widehat{(NS)}_t$, that are diffeomorphic to their projections $Y\equiv\pi_{2,0}(N)\subset W_t$, but are not necessarily holonomic images of smooth sections of $\pi:W\to M$. In fact, by using the embeddings $N\subset\widehat{(NS)}\subset J{\it D}^2(W)\subset J^2_4(W)$, we can repeat above construction to build smooth envelopment solutions of the Cauchy problem $N^{(\infty)}\subset\widehat{(NS)}_{+\infty}\subset J^\infty_4(W)$. Really, any smooth $3$-dimensional space-like submanifold $Y\subset W_t$, identified by a space-like smooth integral manifold $Z\subset\widehat{(NS)}_{+\infty}\subset J^\infty_4(W)$, such that $Z\cong Y$, via the canonical projection $\pi_{\infty,0}:J^\infty_4(W)\to W$, can be locally identified by some smooth implicit functions $\{f_I(x^k,y^j)=0\}_{1\le I\le 5}$, with Jacobian matrix of rank five, where $\{x^k,y^j\}_{1\le k\le 3; 1\le j\le 5}$ are coordinates in the $8$-dimensional affine space $W_t$, and we write $Z=Y^{(\infty)}$. Then, as a by-product we get also a cohomology criterion to classify envelopment solutions according to the first cohomology space $H^1(N_0,\mathbb{Z}_2)$. In fact, a line bundle $V\to N_0^{(\infty)}$ is classified by the first {\em Stiefel-Whitney class} of $V$, that belongs to $H^1( N_0^{(\infty)},\mathbb{Z}_2)\cong H^1( N_0,\mathbb{Z}_2)$. The classifying space is $\mathbb{R}P^\infty$ and the universal principal bundle is $S^\infty\to\mathbb{R}P^\infty\equiv S^\infty/\mathbb{Z}_2\cong \mathbb{R}^\infty/\mathbb{R}^\bullet\cong Gr_1(\mathbb{R}^\infty)$, where $Gr_1(\mathbb{R}^\infty)$ denotes the Grassmannian of $1$-dimensional vector subspaces in $\mathbb{R}^\infty$, and the nonzero element of $\mathbb{Z}_2$ acts by $v\mapsto -v$. Since $S^\infty$ is contractible one has $\pi_i(S^\infty)\cong\pi_i(\mathbb{R}P^{\infty})=0$, $i>1$ and $\pi_1(\mathbb{R}P^\infty)\cong \mathbb{Z}_2$. $\mathbb{R}P^\infty$ is the Eilenberg-Maclane space $K(\mathbb{Z}_2,1)$ \cite{PRA13}. Hence $[N_0,Gr_1]\cong H^1(N_0;\mathbb{Z}_2)$, $f\mapsto f^*\mu$, where $\mu$ is the generator of $H^1(\mathbb{R}P^\infty;\mathbb{Z}_2)\cong\mathbb{Z}_2$. Since one has the bijection $[N_0;Gr_1]\to \mathbf{V}_1(N_0)$, where $\mathbf{V}_1(N_0)$ denotes the set of $1$-dimensional vector bundles over $N_0$, we get the bijection $w_1:\mathbf{V}_1(N_0)\to H^1(N_0;\mathbb{Z}_2)$ that is just the Stiefel-Whitney class for line bundles over $N_0^{(\infty)}$. In conclusion an envelopment solution $V\supset N$, of an admissible Cauchy integral manifolds $N\subset\widehat{(NS)}$, can be identified with some cohomology class of $H^1(N;\mathbb{Z}_2)$. In particular $V$ is orientable iff its first Stiefel-Whitney cohomology class $w_1(V)=0$. Of course does not necessitate that all cohomology classes of $H^1(N;\mathbb{Z}_2)$ should be represented by some envelopment solution $V$ passing for $N$. However, when this happens we say that $N$ is a {\em wholly cohomologic Cauchy manifold} of $\widehat{(NS)}$. This is surely the case when $N\simeq D^3$, i.e., when the space-like, smooth, $3$-dimensional, integral manifold $N$ is homotopy equivalent to the $3$-dimensional disk $D^3$. In fact in such a case $H^1(N;\mathbb{Z}_2)=0$, hence there is an unique cohomology type of envelopment solution passing for $N^{(\infty)}$, (or $N$), the orientable one. Therefore, in such a case $N$ is a wholly cohomologic Cauchy manifold. For example, when $N$ is identified by a smooth space-like section $s_t:M_t\subset M\to W_t\subset W$, $N$ is a wholly cohomologic Cauchy manifold. It is important to remark also that, fixed some smooth Cauchy problem $N\subset\widehat{(NS)}$, it does not necessitate that a local smooth solution $V\epsilon$ should be unique. In fact, even if the Cartan distribution $\mathbf{E}_\infty\subset \widehat{(NS)}_{+\infty}$, is a $4$-dimensional involutive distribution, the manifold $\widehat{(NS)}$, is not finite dimensional, therefore for any point $q\in\widehat{(NS)}$ pass infinity $4$-dimensional integral manifolds tangent to $(\mathbf{E}_\infty)_{q}\subset T_q(\widehat{(NS)}_{+\infty})$. However, we can see that if there are two smooth solutions $V_\epsilon,V'_\epsilon\subset \widehat{(NS)}_{+\infty}$ for a fixed smooth Cauchy problem $N\subset\widehat{(NS)}$, their Stiefel-Whitney classes must coincide with a same cohomology class of $H^1(N,\mathbb{Z}_2)$. In fact, let us denote by $\mathbf{N}_\infty(N)_{\epsilon}$ the {\em integral $\epsilon$-normal bundle} of $N^{(\infty)}$, i.e., $\mathbf{N}_\infty(N)_{\epsilon}\equiv\bigcup_{q\in N^{(\infty)}}\mathbf{N}_\infty(N)_{\epsilon,q}$, with $\mathbf{N}_\infty(N)_{\epsilon,q}\equiv ((\mathbf{E}_\infty)_{q}/T_qN)_\epsilon\equiv\{[0,\epsilon[\zeta_0(q)\}$. Then one has the canonical isomorphism of line bundles over $N^{(\infty)}$: $\mathbf{N}_\infty(N)_{\epsilon}\cong V_\epsilon$, $\lambda\zeta_0(q)\mapsto X_{0,q}(\lambda)$, $\lambda\in[0,\epsilon[$. A similar isomorphism can be recognized between $\mathbf{N}_\infty(N)_{\epsilon}$ and $ V'_\epsilon$, taking into account that also $V'_\epsilon$ can be considered a line bundle: $V'_\epsilon\equiv\bigcup{q\in N^{(\infty)}}X'_{0,q}$. This proves that must necessarily be $\omega_1(V_\epsilon)=\omega_1(\mathbf{N}_\infty(N)_{\epsilon})=\omega_1(V'_\epsilon)$.

So we have shown that any space-like smooth $3$-dimensional integral manifold $N\subset\widehat{(NS)}$ that is diffeomorphic to its projection on $M_t\subset M$, (for some $t$), via the canonical projection $\pi_2:J{\it D}^2(W)_t\to M_t$, or is diffeomorphic to its projection on $W_t$, via the canonical projection $\pi_{2,0}:J{\it D}^2(W)_t\to W_t$, admits smooth solutions, $V\supset N$. Since such integral manifolds (Cauchy manifolds) are not arbitrary ones, but satisfy some constraints, here we shall explicitly prove that such Cauchy manifolds exist. From Tab. \ref{completely-integrable-navier-stokes-equation} we can write $\widehat{(NS)}$ in the form reported in (\ref{completely-integrable-navier-stokes-equation-compact-form})
\begin{equation}\label{completely-integrable-navier-stokes-equation-compact-form}
\scalebox{0.8}{$ \left\{
 \begin{array}{ll}
\hbox{\rm(A)}& \dot x^kG^j_{jk}+\dot x^s_s=0 \\
\hbox{\rm(B)}_0& \dot x^k_0G^j_{jk}+\dot x^s_{s0}=0 \\
\hbox{\rm(B)}_h& \dot x^k(\partial x_h.G^j_{jk})+\dot x^k_hG^j_{jk}+\dot x^s_{sh}=0 \\
\hbox{\rm(C)}& \dot x^j_0=\frac{1}{\rho}\left[-\dot x^s R^j_s-\dot x^s\dot x^i\rho G^j_{is}-\dot x^s\dot x^j_s
                         \rho-\dot x^k_s S^{js}_k-\dot x^j_{is}T^{is}-p_i g^{ij}-\rho(\partial x_i.f)g^{ij}\right]\\
&\equiv H^j\in A[[\dot x^k,\dot x^k_h,\dot x^k_{rs},p_i]] \\
\hbox{\rm(D)}& \theta_0=\frac{1}{\rho C_p}\left[-\rho C_p\dot x^k
\theta_k-\theta_{is} \overline{E}^{is}-\dot x^k\dot x^pW_{kp}-\dot
x^k\dot x^s_p {\overline W}^p_{ks}- \dot
x^k_i\dot x^s_p Y^{ip}_{ks}\right]\\
&\equiv K\in A[[\dot x^k,\dot x^k_h,\theta_k,\theta_{is}]] \\
\end{array}
 \right.$}
\end{equation}
From the prolongation of equation (\ref{completely-integrable-navier-stokes-equation-compact-form})(C) with respect to space coordinates $x^k$, we get
\begin{equation}\label{prolongation-completely-integrable-navier-stokes-equation-compact-form}
 \left\{
 \begin{array}{ll}
\hbox{\rm(C)}_q\hskip 2pt \dot x^j_{0q}&=\frac{1}{\rho}\left[-\dot x^s_q R^j_s-\dot x^s_q R^j_{s,q}-2\dot x^s_q\dot x^i\rho G^j_{is}-\dot x^s\dot x^i\rho G^j_{is,q}-\dot x^s_q\dot x^j_s
                         \rho-\dot x^s\dot x^j_{sq}
                         \rho\right.\\
                         &\left.-\dot x^k_{sq} S^{js}_k-\dot x^k_{s} S^{js}_{k,q}-\dot x^j_{isq}T^{is}-\dot x^j_{is}T^{is}_{,q}
-p_{iq} g^{ij}-p_{i} g^{ij}_{,q}\right.\\
&\left.-\rho(\partial x_{q}\partial x_i.f)g^{ij}-\rho(\partial x_i.f)g^{ij}_{,q}\right]\\
&\equiv H^j_{,q}\in A[[\dot x^k,\dot x^k_h,\dot x^k_{rs},\dot x^k_{rsq},p_i,p_{iq}]] \\
\end{array}
 \right.
\end{equation}
By using equations (\ref{completely-integrable-navier-stokes-equation-compact-form})(C) and (\ref{prolongation-completely-integrable-navier-stokes-equation-compact-form})$(C)_q$, we can rewrite equation (\ref{completely-integrable-navier-stokes-equation-compact-form})$(B)_0$ in the form reported in (\ref{equation-b-in-alternative-form}).
\begin{equation}\label{equation-b-in-alternative-form}
\hbox{\rm(E)} \hskip 10pt G^j_{jk}H^k(x^k,\dot x^s,\dot x^s_k,\dot x^s_{hk},p_i)+H^s_{,s}(x^k,\dot x^s,\dot x^s_k,\dot x^s_{hk},\dot x^s_{hkq},p_i,p_{iq})=0.
\end{equation}
Therefore, the parametric equations of a space-like integral analytic $3$-dimensional manifold $N\subset\widehat{(NS)}$, diffeomorphic to $M_t$, are given in (\ref{parametric-equation-space-like-analytic-cauchy-manifold-navier-stokes-pde}).
\begin{equation}\label{parametric-equation-space-like-analytic-cauchy-manifold-navier-stokes-pde}
 \scalebox{0.8}{$\left\{
 \begin{array}{ll}
 x^0&=t\\
 x^k&=x^k\\
 \dot x^j&=v^j(x^k)\\
 p&=p(x^k)\\
  \theta&=\theta(x^k)\\
   \dot x^j_h&=v^j_{,h}(x^k)\\
    \dot x^j_0&=H^j(x^k)\\
   p_0&=p_0(x^k)\\
    p_k&=p_{,k}(x^k)\\
    \theta_0&=K(x^k)\\
     \theta_k&=\theta_{,k}(x^k)\\
    \dot x^j_{0h}&=H^j_{,h}(x^k)\\
    \dot x^j_{hk}&=v^j_{,hk}(x^k)\\
    \dot x^j_{00}&=H^j_{,0}=(\partial \dot x_s.H^j)\dot x^s_0+(\partial \dot x_s^h.H^j)\dot x^s_{h0}+(\partial \dot x_s^{hk}.H^j)\dot x^s_{hk0}+(\partial p^i.H^j)p_{i0}\\
    &=(\partial \dot x_s.H^j)H^s+(\partial \dot x_s^h.H^j)(\partial x_h.H^s)+(\partial \dot x_s^{hk}.H^j)(\partial x_k\partial x_h.H^s)\\
    &+(\partial p^i.H^j)(\partial x_i.p_{0})\\
    p_{0k}&=(\partial x_k.p_0)\\
    p_{hk}&=(\partial x_h\partial x_k.p)\\
   p_{00}&=p_{00}(x^k)\\
    \theta_{hk}&=(\partial x_h\partial x_k.\theta)\\
    \theta_{0k}&=(\partial x_k.K)\\
     \theta_{00}&=(\partial \dot x_k.K)\dot x^k_0+(\partial \dot x_k^p.K)\dot x^k_{p0}+(\partial \theta^k.K)\theta_{k0}+(\partial\theta^{is}.K)\theta_{is0}\\
     &=(\partial \dot x_k.K)H^k+(\partial \dot x_k^p.K)\dot x_p.H^k)+(\partial \theta^k.K)(\partial x_k.K)+(\partial\theta^{is}.K)(\partial x_i\partial x_s.K)\\
 \end{array}
 \right.$}
\end{equation}
where $v^j=\dot x^j(x^k)$ are solutions of the continuity equation (\ref{completely-integrable-navier-stokes-equation-compact-form})(A). Therefore, equation (\ref{completely-integrable-navier-stokes-equation-compact-form})(A) can be considered a first order equation on the fiber bundle $E_t\equiv M_t\times\mathbf{I}\to M_t$, reported in (\ref{space-like-continuity-equation}).

\begin{equation}\label{space-like-continuity-equation}
E_1\subset J{\it D}(E_t):\hskip 5pt\{\dot x^k G^j_{jk}+\dot x^s_s=0\}.
\end{equation}
One can see that $E_1$ is an involutive, formally integrable and completely integrable PDE. In fact, one has $$\scalebox{0.8}{$\fbox{$\dim(E_1)_{+1}=29$}=\fbox{$\dim E_1=14$}+\fbox{$\dim (g_1)_{+1}=15$}$},$$ hence the mapping $(E_1)_{+1}\to E_1$ is surjective. Furthermore, we get $$\scalebox{0.8}{$\fbox{$\dim(g_1)_{+1}=15$}=\fbox{$\dim g_1=8$}+\fbox{$\dim (g_1)^{(1)}=5$}+\fbox{$\dim (g_1)^{(2)}=2$}+\fbox{$\dim (g_1)^{(3)}=0$}$}.$$ This is enough to state that $g_1$ is an involutive symbol. Therefore, $E_1$ is formally integrable and, since it is analytic it is  also completely integrable. As a by-product, we get that for any analytic solution $\dot x^k=\dot x^k(x^i)$ of $E_1$, we can write:
\begin{equation}\label{expressions-for-time-derivative-after-integration-continuity-equation}
 \left\{
 \begin{array}{ll}
 \dot x^j_0=A^j(x^s)-g^{ij}p_i\\
 \dot x^j_{q0}=A^j(x^s)_{,q}-\frac{1}{\rho}g^{ij}p_{iq}-\frac{1}{\rho}g^{ij}_{,q}p_{i}\\
\end{array}
 \right.
\end{equation}
where $A^j=A^j(x^s)$ are suitable analytic functions. So equation (\ref{equation-b-in-alternative-form})(E) can be rewritten as a second order equation $E_2$ for a function $p=p(x^k)$ as section of the trivial fiber bundle $F_t\equiv M_t\times \mathbb{R}\to M_t$:

\begin{equation}\label{space-like-equation-for-pressure}
E_2\subset J{\it D}^2(F_t):\hskip 5pt\{g^{is}p_{is}+p_iC^i(x^k)-B(x^k)=0\}
\end{equation}
where $C^i$ and $B$ are given analytic functions of $x^k$. This is an involutive, formally integrable PDE, hence it is completely integrable, since it is analytic. In fact, $$\scalebox{0.8}{$\fbox{$\dim(E_2)_{+1}=19$}=\fbox{$\dim(E_2)=12$}+\fbox{$\dim(g_2)_{+1}=7$}$}.$$ Therefore, the map $(E_2)_{+1}\to E_2$ is surjective. Furthermore, $$\scalebox{0.8}{$\fbox{$\dim(g_2)_{+1}=7$}=\fbox{$\dim(g_2)=5$}+\fbox{$\dim(g_2)^{(1)}=2$}+\fbox{$\dim(g_2)^{(2)}=0$}$},$$ hence $g_2$ is an involutive symbol. This concludes the proof. Thus, for any point $q\in E_1$, pass solutions of the continuity equation, and fixing a point $q$ on the infinity prolongation $(E_1)_{+\infty}$, one identifies an analytic solution defined in a suitable neighborhood of $p\equiv\pi_\infty(q)\in M_t$. (Of course also for equation $E_1$ we can identify smooth solutions by solving lower dimension Cauchy problems, by means of Lemma \ref{cauchy-criteria}, i.e., by using envelopment solutions.) Similar considerations can be applied to the equation $E_2$ to identify analytic and smooth solutions of the pressure functions $p=p(x^k)$. This assures that one can identify space-like smooth Cauchy manifolds in $\widehat{(NS)}$, that are diffeomorphic with their canonical projections on $W_t$, for any time $t$.

In {\em\cite{PRA17}} it is proved that the set of full admissible Cauchy integral manifolds is not empty. This result gives us a general criterion to characterize global smooth solutions of the Navier-Stokes equation and completely solves the well-known problem on the existence of global smooth solutions of the Navier-Stokes equation. (For complementary characterizations of the Navier-Stokes equation see also \cite{PRA17, PRA22, PRA23, PRA24, PRA25, PRA26}. There a geometric method to study stability of PDE's and their solutions, related to integral and quantum bordism groups of PDE's, has been introduced, and applied to the Navier-Stokes equation too.)\end{example}

\section{\bf EXTENDED CRYSTAL SINGULAR PDE's}
\vskip 0.5cm

Singular PDE's can be considered singular submanifolds of
jet-derivative spaces. The usual formal theory of PDE's works,
instead, on smooth or analytic submanifolds. However, in many
mathematical problems and physical applications, it is necessary to
work with singular PDE's. (See, e.g., the book by Gromov \cite{GRO}
where he talks of ''partial differential relations'', i.e., subsets
of jet-derivative spaces.) So it is useful to formulate a general
geometric theory for such more general mathematical structures. On
the other hand in order to build a formal theory of PDE's it is
necessary to assume some regularity conditions. So a geometric
theory of singular PDE's must in some sense weak the regularity
conditions usually adopted in formal theory and admit existence of
subsets where these regularity conditions are not satisfied. With
this respect, and by using our formulation of geometric theory of
PDE's and singular PDE's, we study criteria to obtain global
solutions of singular PDE's, crossing singular points. In
particular, some applications concerning singular MHD-PDE's,
encoding anisotropic incompressible nuclear plasmas dynamics, are
given following some our recent works on this subject. The origin of
singularities comes from the fact that there are two regions
corresponding to different components PDE's having different Cartan
distributions with different dimensions. However, by considering
their natural embedding into a same PDE, we can build physically
acceptable solutions, i.e., satisfying the second principle of the
thermodynamics, and that cross the nuclear critical zone of nuclear
energy production. A characterization of such solutions by means of
algebraic topological methods is given also.

The main result of this section is Theorem \ref{main-singular1} that
relates singular integral bordism groups of singular PDE's to global
solutions passing through singular points, and Example \ref{MHD-PDE}
that for some MHD-PDE's characterizes global solutions crossing the
nuclear critical zone and satisfying the entropy production energy
thermodynamics condition.

Let us, now, resume some fundamental definitions and results on the
geometry of PDE's in the category of commutative manifolds,
emphasizing some our recent results on the algebraic geometry of
PDE's, that allowed us to characterize singular PDE's.\footnote{For
general informations on the geometric theory of PDE's see,
e.g.,\cite{BOARD, B-C-G-G-G, C-R-R, GOL, GOL-SPE, GOL-GUIL, GRO,
KRA-LYC-VIN, KUR, L-P, THO1, THO2, THO3}. In particular, for
singular PDE's geometry, see the book \cite{PRA13} and the recent
papers \cite{AG-PRA2, PRA31} where many boundary value problems are
explicitly considered. For basic informations on differential
topology and algebraic topology see e.g., \cite{BOARD, GOL-GUIL,
HIR, M-M, M-S, PRA2, RUD, STO, S-W, SWI, THO1, THO2, THO3, WAL1,
WAL2}.}

\begin{definition}[Algebraic formulation of PDE's] Let $\pi:W\to M$ be a smooth fiber bundle,
$\dim W=m+n$, $\dim M=n$. We denote by $J^k_n(W)$  the {\em space of
all $k$-jets of submanifolds} of dimension $n$ of $W$ and by
$J^k(W)$ the {\em$k$-jet-derivative space of sections} of $\pi$.
Furthermore we denote by $J{\it D}^k(W)$ the {\em$k$-jet-derivative
space} for sections of $\pi$. One has $ J{\it D}^k(W)\cong
J^k(W)\subset J^k_n(W)$. $J^k(W)$ is an open subset of $J^k_n(W)$.
Let ${\frak A}_k$ be the sheaf of germs of differentiable functions
$J{\it D}^k(W)\rightarrow \mathbb{R}$. It is a sheaf of rings, but
also a sheaf of $\mathbb{R}$-modules. A subsheaf of ideals ${\frak
B}_k$ of ${\frak A}_k$ that is also a subsheaf of
$\mathbb{R}$-modules is a {\it PDE of order k} on the fiber bundle
$\pi :W\rightarrow M$. A {\it regular solution} of ${\frak B}_k$ is
a section $s:M\rightarrow W$ such that $f\circ D^ks=0, \forall f\in
{\frak B}_k$. The set of {\it integral points} of ${\frak B}_k$
(i.e., the zeros of ${\frak B}_k$ on $J{\it D}^k(W)$ is denoted by
$J({\frak B}_k)$. The first {\it prolongation} $({\frak B}_k)_{+1}$
of ${\frak B}_k$ is defined as the system of order $k+1$ on
$W\rightarrow M$, defined by the $f\circ \pi _{k,k-1}$ and
$f^{(1)}$, where $f^{(1)}$ on $D^{k+1}s(p)$ is defined by
$f^{(1)}(D^{k+1}s(p))=(
\partial x_\alpha .(f\circ D^ks(p)))$. In local coordinates
$(x^\alpha ,y^j,y^j_\alpha )$ the {\it formal derivative} $f^{(1)}$
is given by $f^{(1)}(x^\alpha ,y^j,y^j_\alpha )=(\partial x_\alpha
.f)+ \sum_{ [\beta ]\le k}y^j_{\beta \alpha }(\partial y^\beta
_j.f)$. The system ${\frak B}_k$ is said to be {\it involutive} at
an integral point $q\in J{\it D}^k(W)$ if the following two
conditions are satisfied: {\em(i)} ${\frak B}_k$ is a {\it regular
local equation} for the zeros of ${\frak B}_k$ at q (i.e., there are
local sections $F_1,...,F_t\in \Gamma (U,{\frak B}_k)$ of ${\frak
B}_k$ on an open neighborhood $U$ of $q$, such that the integral
points of ${\frak B}_k$ in $U$ are precisely the points $q'$ for
which $F_j(q' )=0$ and $dF_1\wedge \cdots\wedge dF_t(q)\not= 0$,
that is $F_1,\cdots,F_t$ are linearly independent at $q;$ and
{\em(ii)} there is a neighborhood $U$ of $q$ such that $\pi
^{-1}_{k+1,k}(U)\bigcap J(({\frak B}_k)_{+1})$ is a fibered manifold
over $U\bigcap J({\frak B}_k)$ (with projection $\pi _{k+1,k})$. For
a system ${\frak B}_k$ generated by linearly independent Pfaffian
forms $\theta ^1,\cdots,\theta ^k$ (i.e., a Pfaffian system) this is
equivalent to the involutiveness defined for distributions.
\end{definition}

\begin{theorem}[\cite{KUR}]
Let ${\frak B}_k$ be a system defined on $J{\it D}^k(W)$, and
suppose that ${\frak B}_k$ is involutive at $q\in J({\frak B}_k)$.
Then, there is a neighborhood $U$ of $q$ satisfying the following.
If $\widetilde q\in J(({\frak B}_k)_{+s})$ and $\pi
_{k+s,k}(\widetilde q)$ is in $U$, then there is a regular solution
$s$ of ${\frak B}_k$ defined on a neighborhood $p=\pi
_{k+s,-1}(\widetilde q)$ of $M$ such that $D^{k+s}s(p)=\widetilde
q$.
\end{theorem}

\begin{theorem}[Cartan-Kuraniski prolongation theorem \cite{KUR, PRA13}] Suppose that
there exists a sequence of integral points $q^{(s)}$ of $({\frak
B}_k)_{+s}, s=0,1,\cdots$, projecting onto each other, $\pi
_{k+s,k+s-1}(q^{(s)})=q^{(s-1)}$, such that: {\rm(a)} $({\frak
B}_k)_{+s}$ is a regular local equation for $J(({\frak B}_k)_{+s})$
at $q^{(s)};$ and {\rm(b)} there is a neighborhood $U^{(s)}$ of
$q^{(s)}$ in $J(({\frak B}_k)_{+s})$ such that its projection under
$\pi _{k+s,k+s-1}$ contains a neighborhood of $q^{(s-1)}$ in
$J(({\frak B}_k)_{+(s-1)})$ and such that $\pi
_{k+s,k+s-1}:U^{(s)}\rightarrow \pi _{k+s,k+s-1}(U^{(s)})$ is a
fibered manifold. Then, $({\frak B}_k)_{+s}$ is involutive at
$q^{(s)}$ for $s$ large enough.\end{theorem}

The algebraic characterization of singular PDE's can be given by
adopting the methods of the algebraic geometry, combined with the
differential algebra. (See e.g., \cite{PRA13}.) Let us go here in
some details about.

\begin{definition}
A {\em differential ring} is a ring $A$ with a finite number $n$ of
commutating derivations $d_1,\cdots,d_n$, $d_id_j-d_jd_i=0$,
$\forall i,j=1,\cdots,n$. A {\em differential ideal} is an ideal
$\mathfrak{a}\subset A$ which is stable by each $d_i$,
$i=1,\cdots,n$.
\end{definition}

A differential ring $(A,\{d_j\}_{1\le j\le n})$ identifies a subring
({\em subring of constants}): $C\equiv cst(A)\equiv\{a\in
A|d_ja=0,\hskip 2pt \forall j=1,\cdots, n\}\subset A$. We may extend
each $d_i$ to a derivation of the full ring of fractions, $Q(A)$,
still denoted by $d_i$ and such that $d_i(a/r)=(rd_ia-ad_ir)/r^2$,
for any $0\not=r$, $a\in A$.

\begin{example}
If $K$ is a differential field with derivations
$\partial_1,\cdots,\partial_\mu$ and $y^k$, $k=1,\cdots,m$, are
indeterminates over $K$, we set $y^k_0=y^k$. Then the polynomial
ring $K[y]_d=K[y^k_\mu,k=1,\cdots,m,\hskip
2pt\mu=\mu_1\cdots\mu_s,\hskip 2pt|\mu|\ge 0]$, can be endowed with
a structure of differential ring by defining the {\em formal
derivations} $d_i\equiv\partial_i+y^k_{\mu+1i}\partial y^\mu_k$. Of
course $K[y]_d$ is not a Noetherian ring. We write
$K[y_q]_d=K[y^k_\mu|k=1,\cdots,m; 0\le|\mu|\le q]$ and one has
$K(y_q)_d=Q(K[y_q]_d)$. We set also $K(y)_d=Q(K[y]_d)$.
\end{example}

\begin{definition}
A {\em differential subring} $A$ of a differential ring $B$ is a
subring which is stable under the derivations of $B$. Similarly we
can define a {\em differential extension} $L/K$ of differential
fields, and such an extension is said to be {\em finitely generated}
if one can find elements $\eta^1,\cdots,\eta^m\in L$ such that
$L=K(\eta^1,\cdots,\eta^m)$. Then the {\em evaluation epimorphism}
 is defined by $K[y]_d\to K[\eta]_d\subset L$, $y^k\mapsto \eta^k$. Its kernel is a prime differential ideal.
\end{definition}

\begin{proposition}\cite{PRA13}
Let $<S>_d$ denote the differential ideal generated by the subset
$S\subset A$, where $A$ is a differential ring. If $A$ is a
differential ring and $a,b\in A$, then one has the following:

{\em(i)} $a^{|\mu|+1}d_\mu b\in<d_\nu(ab)||\nu|\le|\mu|>$.

{\em(ii)} $(d_ia)^{2r-1}\in<a^r>_d$.

{\em(iii)} If $\mathfrak{a}$ is a radical differential ideal of the
differential ring $A$ and $S$ is any subset of $A$, then
$\mathfrak{a}:S\equiv\{a\in A|\mathfrak{a}S\subset\mathfrak{a}\}$ is
again a radical differential ideal of $A$.\footnote{If
$\mathfrak{a}$ is any ideal of $A$, the {\em radical} of
$\mathfrak{a}$ is the following ideal $r(\mathfrak{a})\equiv
\sqrt{\mathfrak{a}}\equiv\{x\in A | x^n\in\mathfrak{a}\hskip
3pt\hbox{for some $n>0$}\}\equiv rad(\mathfrak{a})$. If
$\sqrt{\mathfrak{a}}=\mathfrak{a}$, then $\mathfrak{a}$ is called
{\em radical ideal} or {\em perfect}. One has also that
$r(\mathfrak{a})$ is the intersection of all prime ideals
$\mathfrak{p}\subset A$, containing $\mathfrak{a}$. In particular,
the radical of the zero ideal $<0>$ is the {\em nilradical},
$nil(A)$, of $A$, i.e., the set of all nilpotent elements of $A$.
Therefore $nil(A)$ is the intersection of all prime ideals, (since
all ideals must contain $0$). One has also $nil(A)\subset rad(A)$,
where $rad(A)$ is the ideal of $A$ defined by intersection of all
maximal ideals $\mathfrak{m}\subset A$. If $\mathfrak{a}$ is a
radical ideal, then $A/\mathfrak{a}$ is {\em reduced}, i.e., the set
of its nilpotent elements is reduced to $\{0\}$. In particular
$A/nil(A)$ is reduced. If $\pi:A\to A/\mathfrak{a}$ is the canonical
projection, then $\pi^{-1}(nil(A/\mathfrak{a}))=r(\mathfrak{a})$.}

{\em(iv)} If $\mathfrak{a}$ is a differential ideal of a
differential ring $A$, then $rad(\mathfrak{a})$ is a differential
ideal too.

{\em(v)} One has the following inclusion: $a\hskip 2pt
rad<S>_d\subset rad<aS>_d$, $\forall a\in A$, and for all subset
$S\subset A$.

{\em(vi)}  If $S$ and $T$ are two subsets of a differential ring
$A$, then
$$rad<S>_d.rad<T>_d\subset rad<ST>_d=rad<S>_d\cap rad<T>_d.$$

{\em(vii)} If $S$ is any subset of a differential ring $A$, then we
have:
$$rad<S,a_1,\cdots,a_r>_d=rad<S,a_1>_d\cap\cdots\cap rad<S,a_r>_d.$$
\end{proposition}

\begin{definition}
A {\em differential vector space} is a vector space $V$ over a
differential field $(K,\partial_i)_{1\le i\le n}$ such that are
defined $n$ homomorphisms $d_i$, $i=1,\cdots,n$, of the additive
group $V$ such that: $d_i(av)=(\partial_ia)v+a(d_iv),\hskip
5pt\forall a\in K,\hskip 2pt \forall v\in V$. Then we say that $K$
is a {\em differential field of definition}. \end{definition}

\begin{proposition}\cite{PRA13}
Let $V$ be a differential vector space over a differential field
$K$, with derivations $d_i$, $i=1,\cdots,n$, and let $\{e_j\}_{j\in
I}$ be a basis of $V$. Then the field of definition $\kappa$ of a
differential subspace $W\subset V$ is a differential subfield of $K$
if it contains the field of definition of each
$d_1e_i,\cdots,d_ne_i$ with respect to $\{e_i\}$. \end{proposition}

\begin{definition}
A family $\eta=(\eta^1,\cdots,\eta^m)$ of elements in a differential
extension of the differential field $K$ is said to be {\em
differentially algebraically independent} (or a {\em family of
differential indeterminates}) over $K$, if the kernel of the
evaluation epimorphism $K[y]_d\to K[\eta]_d$ is zero. Otherwise the
family is said to be {\em differentially algebraically dependent}
(or {\em differentially algebraic}) over $K$. \end{definition}

\begin{proposition}
If $K/\kappa$ and $L/\kappa$ are two given differential extensions
with respective derivations $d_K$ and $d_L$, there always exists a
differential free composite field of $K$ and $L$ over $\kappa$.
\end{proposition}

\begin{proof}
The ring $K\bigotimes_\kappa L$ has a natural differential structure
given by $d(a\otimes b)=(d_Ka)\otimes b+a\otimes(d_Lb)$, as
$d_K|_\kappa=d_L|_\kappa=\partial$. On the other hand there is a
finite number of prime ideals ${\frak p}_i\subset K\bigotimes_\kappa
L$ such that $\bigcap_i\mathfrak{p}_i=0$ and
$\mathfrak{p}_i+\mathfrak{p}_j=<1>$, $\forall i\not=j$. Now we have
the following lemma.

\begin{lemma}
If $\mathfrak{a}_1,\cdots,\mathfrak{a}_r$ are ideals of a
differential ring $A$ such that $\mathfrak{a}_i+\mathfrak{a}_j=A$,
$\forall i\not=j$, and $\mathfrak{a}_1\cap\cdots\cap\mathfrak{a}_r$
is a differential ideal of $A$, then each $\mathfrak{a}_i$ is a
differential ideal too. \end{lemma}

Therefore we can conclude that each $\mathfrak{p}_i$ is a
differential ideal, hence the proposition is proved. \end{proof}

\begin{lemma}
A family $\eta$ is differentially algebraic over $K$ iff a
differential polynomial $P\in\mathfrak{p}$ exists such that
$(\partial y_P.P)\not\in\mathfrak{p}$, where $y_P$ is the highest
power of $y_p$ appearing in $P$. $S_P\equiv(\partial y_P.P)$ is
called the {\em separout} of $P$. (The {\em initial} of $P$ is the
coefficient of the highest power of $y_P$ appearing in $P$ and it is
denoted by $I_P$. More precisely one has $P=I_P(y_P)^r+\hbox{terms
of lower degree}$.) \end{lemma}

\begin{proposition}\cite{PRA13}
If $S$ is any subset of a differential ring $A$ and $r\ge 0$ is any
integer, we call {\em $r$-prolongation} of $S$, the ideal
$$(S)_{+r}=<d_\nu a|a\in S,\hskip 2pt 0\le|\nu|\le r>\subset A.$$
One has the following properties: {\em(i)}
$(S)_{+(r+s)}=((S)_{+s})_{+r}$. {\em(ii)} $(S)_{+\infty}=<S>_d$.

{\em(iii)} Let $\mathfrak{a}$ be a differential ideal of the
differential ring $K[y]_d$. We set $\mathfrak{a}_q=\mathfrak{a}\cap
K[y_q]_d$, $\mathfrak{a}_0={\frak a}\cap K[y]_d$,
$\mathfrak{a}_\infty=\mathfrak{a}$. We call the {\em
$r$-prolongation} of $\mathfrak{a}_q$, the following ideal:
$$(\mathfrak{a}_q)_{+r}=<d_\nu P|P\in\mathfrak{a}_q,\hskip 2pt 0\le|\nu|\le r>\subset K[y_{q+r}]_d.$$
One has:
$$(\mathfrak{a}_q)_{+r}\subseteq\mathfrak{a}_{q+r},\quad
(\mathfrak{a}_q)_{+\infty}\subseteq\mathfrak{a},\quad ({\frak
a}_q)_{+r}\cap K[y_q]_d=\mathfrak{a}_{q},\hskip 2pt\forall q,r\ge
0.$$ \end{proposition}

With algebraic sets it is better to consider radical ideals. Hence
if $\mathfrak{r}\subset K[y]_d$ is a radical differential ideal,
then $\mathfrak{r}_q$ is a radical ideal of $K[y_q]_d$, for all
$q\ge 0$. Then if $E_q=Z(\mathfrak{r}_q)$ is the algebraic set
defined over $K$ by $\mathfrak{r}_q=I(E_q)$, we call
{\em$r$-prolongation} of $E_q$ the following algebraic set:
$(E_q)_{+r}=Z(({\frak r}_q)_{+r})$. In general one has
$(\mathfrak{r}_q)_{+r}\subseteq\mathfrak{r}_{q+r}$, hence
$rad((\mathfrak{r}_q)_{+r})\subseteq \mathfrak{r}_{q+r}$. Therefore,
in general one has: $E_{q+r}\subseteq (E_q)_{+r}$.

\begin{proposition}\cite{PRA13}
Let $\mathfrak{p}\subset K[y]_d$ be a prime differential ideal. Then
we can identify each field $L_q=Q(K[y_q]_d/\mathfrak{p}_q)$ with a
non-differential subfield of $L=Q(K[y]_d/\mathfrak{p})$ and we have:
$K\subseteq L_0\subseteq\cdots\subseteq L_\infty=L$. Then there are
vector spaces $R_q$ over $L_q$ or $L$ defined by the following
linear system:
$$(\partial y^\mu_k.P_\tau)(\eta)v^k_\mu=0,\hskip 5pt \left\{1\le \tau\le
t,\quad 1\le k\le m,\quad |\mu|=q\right\},$$ where $\eta$ is a
generic solution of $\mathfrak{p}$ and $P_1,\cdots,P_t$ are
generating $\mathfrak{p}_q$. Such result does not depend on the
generating polynomials. We can also define the vector space $g_q$
({\em symbol}) over $L_q$ or $L$, by means of the linear system:
$$(\partial y^\mu_k.P_\tau)(\eta)v^k_\mu=0,\hskip 5pt \left\{1\le \tau\le
t,\quad 1\le k\le m,\quad 0\le |\mu|\le q\right\}.$$ For the
prolongations $(g_q)_{+r}$ one has, in general,
$g_{q+r}\subseteq(g_q)_{+r}$, $\forall q,r\ge 0$.
\end{proposition}

\begin{definition}
We say that $R_q$ or $g_q$ is {\em generic} over $E_q$, if one can
find a certain number of maximum rank determinants $D_\alpha$ that
cannot be all zero at a generic solution of $\mathfrak{p}$.
\end{definition}

\begin{proposition}
$R_q$ or $g_q$ is {\it generic} if we may find polynomials
$A_\alpha, B_\tau\in K[y_q]_d$ such that:
$$\sum_\alpha A_\alpha D_\alpha+\sum_\tau B_\tau P_\tau=1.$$
Furthermore, $R_q$ or $g_q$ are projective modules over the ring
$K[y_q]_d/\mathfrak{p}_q\subset K[y]_d/\mathfrak{p}$.
\end{proposition}

\begin{proof}
It follows directly from the Hilbert theorem of zeros. (See
\cite{PRA13}.)
\end{proof}

\begin{theorem}[Primality criterion \cite{PRA13}]
Let $\mathfrak{p}_q\subset K[y_q]_d$ and $\mathfrak{p}_{q+1}\subset
K[y_{q+1}]_d$ be prime ideals such that
$\mathfrak{p}_{q+1}=(\mathfrak{p}_q)_{+1}$ and
$\mathfrak{p}_{q+1}\cap K[y_q]_d=\mathfrak{p}_q$. If the symbol
$g_q$ of the variety $R_q$ defined by $\mathfrak{p}_q$ is
$2$-acyclic and its first prolongation $g_{q+1}$ is generic over
$E_q$, then $\mathfrak{p}=(\mathfrak{p}_q)_{+\infty}$ is a prime
differential ideal with $\mathfrak{p}\cap
K[y_{q+r}]_d=(\mathfrak{p}_q)_{+r}$, for all $r\ge 0$.

Let $\mathfrak{r}_q\subset K[y_q]_d$ and $\mathfrak{r}_{q+1}\subset
K[y_{q+1}]_d$ be radical ideals such that
$\mathfrak{r}_{q+1}=(\mathfrak{r})_{+1}$ and $\mathfrak{r}_{q+1}\cap
K[y_q]_d=\mathfrak{r}_q$. If the symbol $g_q$ of the algebraic set
$E_q$ defined by $\mathfrak{r}_q$ is $2$-acyclic and its first
prolongation $g_{q+1}$ is generic over $E_q$, then
$\mathfrak{r}=(\mathfrak{r}_q)_{+\infty}$ is a radical differential
ideal with $\mathfrak{r}\cap K[y_{q+r}]_d= (\mathfrak{r}_q)_{+r}$,
for all $r\ge 0$.
\end{theorem}

\begin{theorem}[Differential basis] If $\mathfrak{r}$ is a
differential ideal of $K[y]_d$, then $\mathfrak{r}=rad(({\frak
r}_q)_{+\infty})$ for $q$ sufficiently large.
\end{theorem}

\begin{proof}
In fact one has the following lemma.

\begin{lemma}
If $\mathfrak{p}$ is a prime ideal of $K[y]_d$, then for $q$
sufficiently large, there is a polynomial $P\in K[y_q]_d$ such that
$P\not\in \mathfrak{p}_q$ and $P\mathfrak{p}_{q+r}\subset
rad((\mathfrak{p}_q)_{+r})\subset\mathfrak{p}_{q+r}$, for all $r\ge
0$. \end{lemma}

After above lemma the proof follows directly. \end{proof}

Every radical differential ideal of $K[y]$ can be expressed in a
unique way as the non-redundant intersection of a finite number of
prime differential ideals. The smallest field of definition $\kappa$
of a prime differential ideal $\mathfrak{p}\subset K[y]$ is a
finitely generated differential extension of $\mathbb{Q}$.

\begin{example}
With $n=2$, $m=2$, $q=1$. Let us consider the differential
polynomial $P=y^1_1y^2_2-y^2_1y^1_2-1$. We obtain for the symbol
$g_1$: $y^1_1v^2_2+y^2_2v_1^1-y^1_2v_1^2-y^2_1v^1_2=0$. Setting
$v^k_i=y^k_lw^l_i$ we obtain
$(y^1_1y^2_2-y^2_1y^1_2)(w^2_2+w^1_1)=0$ and thus $w^2_2+w^1_1=0$ on
$E_1$. Hence $g_1$ is generic. One can also set $P_1=y^1_2$,
$P_2=y^2_1$ and we get the relation: $y^2_2P_1-y^1_2P_2-P\equiv 1$.
A similar result should hold for $E_1$. $g_1$ is involutive and the
differential ideal generated by $P$ in $\mathbb{Q}<y^1,y^2>$ is
therefore a prime ideal.
\end{example}

\begin{definition}
A {\em differentially algebraic extension} $L$ over of a
differential field $K$ is a differential extension over $K$ where
every element of $L$ is differentially algebraic over $K$.

The {\em differential transcendence degree} of a differential
extension $L/K$ is the number of elements of a maximal subset $S$ of
elements of $L$ that are differentially transcendental over $K$ and
such that $L$ becomes differentially algebraic over $K(S)$. We shall
denote such number by $trd_d(L/K)$.
\end{definition}

\begin{theorem}[\cite{PRA13}]
One has the following formula:
$$\dim(\mathfrak{p}_{q+r})=\dim(\mathfrak{p}_{q-1})+\sum_{1\le i\le n}{{(r+i)!}\over{r!i!}}\alpha^i_q,
\hskip 5pt \forall r\ge 0,$$ where $\alpha^i_q$ is the character of
the corresponding system of PDE's. The character $\alpha^i_q$ of a
$q$-order PDE $E_q\subset J{\it D}^q(W)$, $\pi:W\to M$, $\dim M=n$,
with symbol $ g_q$, is the integer
$\alpha^i_q\equiv\dim(g_q^{(i-1)})_p-\dim(g_q^{(i)})_p$, $p\in E_q$,
where
$(g^{(i)}_q)_p\equiv\{\zeta\in(g_q)_p|\zeta(v_1)=\cdots=\zeta(v_i)=0\}$,
where $(v_1,\cdots,v_n)$ is the natural basis in $ T_{\pi_k(p)}M$.

The character $\alpha^n_q$ and the smallest non-zero character only
depend on the differential extension $L/K$ and not on the
generators. In particular, one has: $trd_d(L/K)=\alpha^n_q$.

If $\zeta$ is differentially algebraic over $K(\eta)_d$ and $\eta$
is differentially algebraic over $K$, then $\zeta$ is differentially
algebraic over $K$.

If $L/K$ is a differential extension and $\xi,\eta\in L$ are both
differentially algebraic over $K$, then $\xi+\eta$, $\xi\eta$,
$\xi/\eta$, ($\eta\not=0$), and $d_i\xi$ are differentially
algebraic over $K$.
\end{theorem}

\begin{theorem}
Let $(A,\{\partial_j\}_{1\le j\le n})$ be a differential ring. The
set $D(A)$ of differential operators over $(A,\{\partial_j\}_{1\le
j\le n})$ is a non-commutative filtered ring and a filtered bimodule
over $A$. \end{theorem}

\begin{proof}
If $y$ is a differential indeterminate over $A$, we may introduce
the {\em formal derivatives} $d_1,\cdots,d_n$ which are such that
$d_id_j-d_jd_i=0$, $\forall i,j=1,\cdots,n$, and are defined by:
$d_i(ay)=(\partial_ia)y+a(d_iy)$. We shall write $d_iy=y_i$,
$d_iy_\mu=y_{\mu+1i}$, where $\mu$ is the multi-index
$\mu=(\mu_1,\cdots,\mu_n)$ with length $|\mu|=\mu_1+\cdots+\mu_n$.
If $y=(y^1,\cdots,y^m)$, we set
$d_\mu=(d_1)^{\mu_1}\cdots(d_n)^{\mu_n}$ and $d_\mu y^k=y^k_\mu$.
Any differential operator of order $q$ over $A$ can be written in
the form $P=\sum_{0\le\mu\le q}a^\mu d_\mu$, $a^\mu\in A$. Set
$ord(P)=q$. Then, we can write $D(A)\cong A[d_1,\cdots,d_n]\equiv
A[d]$ the {\em ring of partial differential operators} over $A$ with
derivatives $d_1,\cdots,d_n$. The addition rule is clear. The
multiplication rule comes from the Leibniz formula:

$$\left\{\begin{array}{l}
  \partial_\nu(ab)=\sum_{\lambda+\mu=\nu}{{\nu!}\over{\lambda!\mu!}}(\partial_\lambda a)(\partial_\mu b)\\
  \\
d_\nu(ay)=\sum_{\lambda+\mu=\nu}{{\nu!}\over{\lambda!\mu!}}(\partial_\lambda
a)d_\mu y\\
\end{array}\right\}\Rightarrow \hskip 3pt d_\nu
a=\sum_{\lambda+\mu=\nu}{{\nu!}\over{\lambda!\mu!}}(\partial_\lambda
a)d_\mu.$$

Here we have put $\mu!=\mu_1!\cdots\mu_n!$.  With these rules $D(A)$
becomes a non-commutative ring and a bimodule over $A$. In fact, the
previous formula defines the right action of $A$ on $D(A)$. The left
action of $A$ on $D(A)$ is simply the multiplication on the left by
$A$, that is $aP=a(\sum_{0\le\mu\le q}a^\mu d_\mu)=\sum_{0\le\mu\le
q}aa^\mu d_\mu$. Now, the filtration of $D(A)$ is naturally induced
by filtration of spaces of differential operators. More precisely
$D_q(A)=\{P\in D(A)| ord(P)\le q\}$, where $ord(P)=sup\{|\mu|
|a^\mu\not=0\}$. We set $D_{-1}(A)=0$ and $D_0(A)=A$. Then,
$D_q(A)\subset D_{q+1}(A)$, $D(A)=\bigcup_{q\ge 0}D_q(A)$ and
$D_q(A)D_p(A)\subseteq D_{p+q}(A)$. \end{proof}

\begin{theorem}[Algebraic criterion for formal
integrability \cite{PRA13}] Let $Z_q=Z(\mathfrak{p}_q)$ be the
variety defined by means of ideal $\mathfrak{p}_q\subset K[y_q]_d$
such that the following conditions are verified:

{\em(i)} $(\mathfrak{p}_q)_{+1}=\mathfrak{p}_{q+1}\subset
K[y_{q+1}]_d$ is also a prime ideal.

{\em(ii)} $\mathfrak{p}_{q+1}\cap K[y_{q}]_d=\mathfrak{p}_q$.

{\em(iii)} $g_{q+1}$ is generic over $E_q$.

{\em(iv)} $g_{q}$ is $2$-acyclic.

Then $(\mathfrak{p}_q)_{+\infty}=\mathfrak{p}\subset K[y]_d$ is a
prime differential ideal, where $\mathfrak{p}$ is the differential
ideal generated by a finite number of differential polynomials
$P_1,\cdots,P_t$, defining $E_q$, and $E_q$ is formally integrable.
If one of these conditions is not satisfied we get that
$\mathfrak{p}$ is not a prime ideal, hence we have a factorization
of $\mathfrak{p}$. In other words the PDE is not formally
integrable.
\end{theorem}

\begin{proof}
For a detailed proof see \cite{PRA13}.\end{proof}
\begin{table}[t]\centering
\caption{Examples of singular PDE's defined by differential polynomials.}
\label{examples-singular-pdes}
\begin{tabular}{|l|l|}
\hline
{\rm{\footnotesize Name}}&{\rm{\footnotesize Singular PDE}} \\
\hline
{\rm{\footnotesize  PDE with node and triple point}}&{\rm{\footnotesize  $ p_1\equiv(u^1_x)^4+(u^2_y)^4-(u^1_x)^2=0$}}\\
{\rm{\footnotesize  $ E_1\subset J{\it D}(E)$}}&{\rm{\footnotesize  $p_2\equiv(u^2_x)^6+(u^1_y)^6-u^2_xu^1_y=0$}}\\
\hline
{\rm{\footnotesize  PDE with cusp and tacnode}}&{\rm{\footnotesize  $q_1\equiv(u^1_x)^4+(u^2_y)^4-(u^1_x)^3+(u^2_y)^2=0$}}\\
{\rm{\footnotesize  $\bar E_1\subset J{\it D}(E)$}}&{\rm{\footnotesize  $q_2\equiv(u^2_x)^4+(u^1_y)^4-(u^2_x)^2(u^1_y)-(u^2_x)(u^1_y)^2=0$}}\\
\hline
{\rm{\footnotesize  PDE with conical double point,}}&{\rm{\footnotesize  $r_1\equiv(u^1)^2-(u^1_x)(u^2_y)^2=0$}}\\
{\rm{\footnotesize  double line and pinch point}} &{\rm{\footnotesize  $r_2\equiv(u^2)^2-(u^2_x)^2-(u^1_y)^2=0$}}\\
{\rm{\footnotesize  $\tilde E_1\subset J{\it D}(F)$}}&{\rm{\footnotesize  $r_3\equiv(u^3)^3+(u_y^3)^3+(u^2_x)(u^3_y)=0$}}\\
\hline
\multicolumn{2}{l}{\rm{\footnotesize  $\pi:E\equiv\mathbb{R}^4\to\mathbb{R}^2,\quad(x,y,u^1,u^2)\mapsto(x,y)$.\quad
$\bar\pi:F\equiv\mathbb{R}^5\to\mathbb{R}^2,\quad(x,y,u^1,u^2,u^3)\mapsto(x,y)$.}}\\
\multicolumn{2}{l}{\rm{\footnotesize  $\mathfrak{a}\equiv<p_1,p_2>\subset
A$, $\mathfrak{b}\equiv<q_1,q_2>\subset A$, $\mathfrak{c}\equiv<r_1,r_2,r_3>\subset B$.}}\\
\end{tabular}
\end{table}

\begin{example}[Some singular PDE's]\label{examples-singular-PDEs}
In Tab. \ref{examples-singular-pdes} we report some singular PDE's having  some algebraic
singularities. These are singular PDE's of first order defined on
$J{\it D}(E)\cong\mathbb{R}^{8}$ for the first two and on $J{\it
D}(F)\cong\mathbb{R}^{11}$ for the third. To the ideals
$\mathfrak{a}\equiv<p_1,p_2>\subset A$,
$\mathfrak{b}\equiv<q_1,q_2>\subset A$ and
$\mathfrak{c}\equiv<r_1,r_2,r_3>\subset B$, where
$A\equiv\mathbb{R}[u^1,u^2,u^1_x,u^1_y,u^2_x,u^2_y]$ and
$B\equiv\mathbb{R}[u^1,u^2,u^3,u^1_x,u^1_y,u^2_x,u^2_y,u^3_x,u^3_y]$,
one associates the corresponding algebraic sets
\begin{equation}\label{algebraic-sets-associated-ideals}
\left\{\begin{array}{l}
 E_1=\{q\in\mathbb{R}^8|f(q)=0, \forall f\in
\mathfrak{a}\}\subset\mathbb{R}^8\\
\\
\bar E_1=\{q\in\mathbb{R}^8|f(q)=0, \forall f\in
\mathfrak{b}\}\subset\mathbb{R}^8\\
\\
\tilde E_1=\{q\in\mathbb{R}^{11}|f(q)=0, \forall f\in
\mathfrak{c}\}\subset\mathbb{R}^{11}.\\
\end{array}\right.
 \end{equation}
 These algebraic sets are in
bijective correspondence with the corresponding radicals:
$r(\mathfrak{a})=\{g\in A|g(p)=0, \forall p\in
E_1\}\supset\mathfrak{a}$, $r(\mathfrak{b})=\{g\in A|g(p)=0, \forall
p\in \bar E_1\}\supset\mathfrak{b}$, $r(\mathfrak{c})=\{g\in
B|g(p)=0, \forall p\in \tilde E_1\}\supset\mathfrak{c}$. (This
follows from the Hilbert theorem of zeros \cite{PRA13}.)

Let us consider in some details the singular PDE $\tilde E_1\subset
J{\it D}(F)$, in order to see existence of global algebraic singular
solutions and characterize their stability. We have the following
representation: $\tilde E_1=\bigcup A_1\bigcup A_2\bigcup A_3$,
where
\begin{equation}\label{example-analytic-components-1-2}
A_j\equiv\left\{
\begin{array}{l}
1)\quad   u_x^1=\frac{(u^1)^2}{(u^2_y)^2}\\
  \\
2)\quad   u^1_y=s(j)\sqrt{(u^2)^2-(u^2_x)^2}\\
  \\
3)\quad   u^2_x=-\frac{(u^3)^3}{u^3_y}-(u^3_y)^2\\
\end{array}
\left|
\begin{array}{l}
(u^2_x)^2\le(u^2)^2\\
\\
u^2_y\not=0\\
\\
u^3_y\not=0\\
\end{array} \right.\right\}
\subset \tilde E_1\subset J{\it D}(F)
\end{equation}
\begin{equation}\label{example-analytic-components-3}
A_3\equiv \left\{q\in J{\it
D}(F)\cong\mathbb{R}^{11}|u^1=u^2=u^3=u^2_x=u^1_y=u^2_y=u^3_y=0\right\}.
\end{equation}

In {\em(\ref{example-analytic-components-1-2})} we put $s(1)=1$ and
$s(2)=-1$. $A_3\cong\mathbb{R}^4$ is the set of singular points of
$\tilde E_1$. Instead $A_1$ and $A_2$ are formally integrable and
completely integrable PDE's. In fact, one has the exact
commutative diagrams {\em(\ref{example-exact-commutative-diagrams-a-b})}, with $j=1,2$:
\begin{equation}\label{example-exact-commutative-diagrams-a-b}
\xymatrix{0\ar[d]&&\\
A_j\ar[d]\ar[r]&F\ar@{=}[d]\ar[r]&0\\
J{\it D}(F)\ar[r]_{\pi_{1,0}}&F\ar[r]&0}\hskip 5pt
\xymatrix{0\ar[d]&&\\
A_3\ar[d]\ar[r]&M\ar@{=}[d]\ar[r]&0\\
J{\it D}(F)\ar[r]_{\pi_1}&M\ar[r]&0}
\end{equation}
This can be seen rewriting equations in
{\em(\ref{example-analytic-components-1-2})} in the following
equivalent way:
\begin{equation}\label{example-analytic-components-1-2-second-form}
A_j\equiv\left\{
\begin{array}{l}
1)\quad   u_x^1=\frac{(u^1)^2}{(u^2_y)^2}\\
  \\
2)\quad   u^1_y=s(j)\sqrt{(u^2)^2-[\frac{(u^3)^3}{u^3_y}+(u^3_y)^2]^2}\\
  \\
3)\quad   u^2_x=-[\frac{(u^3)^3}{u^3_y}+(u^3_y)^2]\\
\end{array}
\left|
\begin{array}{l}
(u^2_x)^2\le(u^2)^2\\
\\
u^2_y\not=0\\
\\
u^3_y\not=0\\
\end{array} \right.\right\}
\subset \tilde E_1\subset J{\it D}(F).
\end{equation}

The first prolongation $(A_j)_{+1}$, $j=1,2$, is
given by the following equations:
\begin{equation}\label{first-prolongation-example}
\scalebox{0.9}{$\left\{\begin{array}{l}
  1)\quad r_1=0\\
  2)\quad  r_2=0\\
  3)\quad  r_3=0\\
  \\
  4)\quad  u^1_{xx}=[2u^1u^1_x(u^2_y)^2-(u^1)^22u^2_yu^2_{yx}]/(u^2_y)^4\\
  \\
  5)\quad  u^1_{xy}=[2u^1u^1_y(u^2_y)^2-(u^1)^22u^2_yu^2_{yy}]/(u^2_y)^4\\
  \\
  6)\quad  u^1_{xy}=s(j)\frac{a}{2}/\sqrt{(u^2)^2-[\frac{(u^3)^3}{u^3_y}+(u^3_y)^2]^2}\\
  \\
  7)\quad  u^1_{yy}=s(j)\frac{b}{2}/\sqrt{(u^2)^2-[\frac{(u^3)^3}{u^3_y}+(u^3_y)^2]^2}\\
  \\
  8)\quad  u^2_{xx}=-[\frac{3(u^3)^2u^3_xu^3_y-(u^3)^3u^3_{yx}}{(u^3_y)^2}+2u^3_yu^3_{yx}]\\
  \\
  9)\quad  u^2_{xy}=-[\frac{3(u^3)^2(u^3_y)^2-(u^3)^3u^3_{yy}}{(u^3_y)^2}+2u^3_yu^3_{yy}]\\
\end{array}\left|\begin{array}{l}
(u^2_x)^2\le(u^2)^2\\
\\
u^2_y\not=0\\
\\
u^3_y\not=0\\
\end{array} \right.\right\}\equiv (A_j)_{+1}\subset J{\it D}^2(F)$}
\end{equation}
where
\begin{equation}
\left\{\begin{array}{l}
a=2u^2u^2_x-2[\frac{(u^3)^3}{u^3_y}+(u^3_y)^2][\frac{3(u^3)^2u^3_xu^3_y-(u^3)^3u^3_{yx}}{(u^3_y)^2}+2u^3_yu^3_{yx}]\\
\\
b=2u^2u^2_y-2[\frac{(u^3)^3}{u^3_y}+(u^3_y)^2][\frac{3(u^3)^2(u^3_y)^2-(u^3)^3u^3_{yy}}{(u^3_y)^2}+2u^3_yu^3_{yy}]\\
\end{array}\right.
\end{equation}
Then by using {\em(\ref{first-prolongation-example})(9)} to
substitute $u^2_{yx}$ in {\em(\ref{first-prolongation-example})(1)},
and by using the two different expressions of $u^1_{xy}$ in
 {\em(\ref{first-prolongation-example})(5)} and
 {\em(\ref{first-prolongation-example})(6)} to obtain an explicit expression of $u^2_{yy}$, we get the following equations for the first prolongation of $(A_j)_{+1}$:

\begin{equation}\label{first-prolongation-example-2}
\scalebox{0.8}{$\left\{\begin{array}{l}
  1)\quad r_1=0\\
  2)\quad  r_2=0\\
  3)\quad  r_3=0\\
  \\
  4)\quad  u^1_{xx}=[2u^1u^1_x(u^2_y)^2-(u^1)^22u^2_yc]/(u^2_y)^4\\
  \\
  5)\quad  u^1_{xy}=s(j)\frac{a}{2}/\sqrt{(u^2)^2-[\frac{(u^3)^3}{u^3_y}+(u^3_y)^2]^2}\\
  \\
  6)\quad  u^1_{yy}=s(j)\frac{b}{2}/\sqrt{(u^2)^2-[\frac{(u^3)^3}{u^3_y}+(u^3_y)^2]^2}\\
  \\
  7)\quad  u^2_{xx}=-[\frac{3(u^3)^2u^3_xu^3_y-(u^3)^3u^3_{yx}}{(u^3_y)^2}+2u^3_yu^3_{yx}]\\
  \\
  8)\quad  u^2_{xy}=-[\frac{3(u^3)^2(u^3_y)^2-(u^3)^3u^3_{yy}}{(u^3_y)^2}+2u^3_yu^3_{yy}]\\
  \\
  9)\quad u^2_{yy}=[\frac{2u^1u^1_y}{(u^2_y)^2}-s(j)\frac{1}{2}\frac{a}{\sqrt{(u^2)^2-[\frac{(u^3)^3}{u^3_y}+(u^3_y)^2]^2}}]
  \frac{1}{2}\frac{(u^2_y)^3}{(u^1)^2}\\
\end{array}\left|\begin{array}{l}
(u^2_x)^2\le(u^2)^2\\
\\
u^2_y\not=0\\
\\
u^3_y\not=0\\
\end{array} \right.\right\}\equiv (A_j)_{+1}\subset J{\it D}^2(F)$}
\end{equation}
with
\begin{equation}
c=-[\frac{3(u^3)^2(u^3_y)^2-(u^3)^3u^3_{yy}}{(u^3_y)^2}+2u^3_yu^3_{yy}].
\end{equation}
Therefore also $(A_j)_{+1}$ are analytic submanifolds of $J{\it
D}^2(F)$, for $j=1,2$. Furthermore, since $(\dim
(A_j)_{+1}=11)=(\dim(A_j=8)+(\dim(g_1)_{+1}=3)$, we see that the
canonical projections $\pi_{2,1}:(A_j)_{+1}\to A_j$, $j=1,2$, are
affine subbundles of $J{\it D}^2(F)\to J{\it D}(F)$, with associated
vector bundles $(g_1)_{+1}\to A_j$, $j=1,2$. Finally the symbol
$\dim(g_1)_{q\in A_j}=3$, $j=1,2$, and $\dim(\partial x\rfloor
(g_1)_{q\in A_j})=0$. Therefore, $(g_1)_{q\in A_j}$ are involutive.
This is enough to conclude that $A_j$ are formally integrable, and
since they are also analytic they are completely integrable too.
Furthermore, taking into account that $\dim A_j=8>2\times 2+1=5$, we
can apply Theorem 2.15 in \cite{PRA15}, (here reported in Section
3), to calculate the weak and singular integral bordism groups of
$A_j$. One has $\Omega_{1,w}^{A_j}=\Omega_{1,s}^{A_j}=0$, $j=1,2$.
Therefore, $A_j$ are extended $0$-crystal PDE's. So for any two
admissible closed $1$-dimensional smooth integral manifolds $N_0,
N_1\subset A_j$, there exists a (singular) $2$-dimensional integral
manifold, solution $V\subset A_j$, such that $\partial V=N_0\DU
N_1$.\footnote{Let us emphasize that by using Lemma \ref{cauchy-criteria} we can identify admissible smooth $1$-dimensional integral manifolds in $A_j$, $j=1,2$. In fact, $A_j$ are formally integrable and completely integrable PDE's. So we can use Lemma \ref{cauchy-criteria}(1), but also Lemma \ref{cauchy-criteria}(3), since $\pi_{1,0}(A_j)=F$, and $A_j\to F$ are affine subbundles of $J{\it D}(F)\to F$, with associated vector bundle the symbol $g_1$.} Such a solution is smooth iff all the integral characteristic
numbers of $N_0$ are equal to ones of $N_1$.

The Cartan distribution on $\tilde E_1$ is given by the following
vector fields
\begin{equation}\label{cartan-distribution-example}
\zeta=X^x(\partial x+u^k_x\partial u_k)+X^y(\partial y+u^k_y\partial
u_k)+Y_k^x\partial u^k_x+Y_k^y\partial u^k_y
\end{equation}
such that the following equations are satisfied:
\begin{equation}\label{conditions-cartan-distribution-example}
\left\{
\begin{array}{l}
  Y^x_1(u^2_y)^2+Y^y_22u^1_xu^2_y-2u^1(u^1_xX^x+u^1_yX^y)=0\\
  \\
  Y^y_12u^1_y+Y^x_22u^2_x-2u^2(u^2_xX^x+u^2_yX^y)=0\\
  \\
Y^x_2u^3_y+Y^y_3(u^2_x+3(u^3_y)^2)+3(u^3)^2(u^3_xX^x+u^3_yX^y)=0.\\
\end{array}
\right.
\end{equation}
Therefore $\dim (\mathbf{E}_1)_{q\in A_j}=5$, $j=1,2$. For example
on $A_1$, one has the following expression of the Cartan vector
field:

\begin{equation}
\left\{\begin{array}{ll}
  \zeta &=X^x[\partial x+u^k_x\partial u_k+\frac{2u^1u^1_x}{(u^2_y)^2}\partial u^1_x+\frac{u^2u^2_x}{u^1_y}
  \partial u^1_y-\frac{3(u^3)^2u^3_x}{[u^2_x+3(u^3_y)^2]}\partial u^3_y]\\
  \\
  &+ X^y[\partial y+u^k_y\partial u_k+\frac{2u^1u^1_y}{(u^2_y)^2}\partial u^1_x+\frac{u^2u^2_y}{u^1_y}\partial u^1_y
  -\frac{3(u^3)^2u^3_y}{[u^2_x+3(u^3_y)^2]}\partial u^3_y] \\
  \\
 &+  Y^x_2[\partial u^2_x-\frac{u^2_x}{u^1_y}\partial u^1_y-\frac{u^3_y}{[u^2_x+3(u^3_y)^2]}\partial
  u^3_y]+Y^x_3\partial u^3_x+Y^y_2[-\frac{2u^1_x}{u^2_y}\partial u^1_x+\partial
  u^2_y].\\
\end{array}\right.
\end{equation}

Instead, since equations
{\em(\ref{conditions-cartan-distribution-example})} are identities
for $q\in A_3$, we get that $(\mathbf{E}_1)_{q\in
A_3}=\mathbf{E}_1(F)_{q\in A_3}$, i.e., in the singular points the
Cartan distribution of $\tilde E_1$ coincides with the Cartan
distribution of $J{\it D}(F)$ that is just given by vector fields
given in {\em(\ref{cartan-distribution-example})} for arbitrary
functions $X^x, X^y, Y^x_k, Y^y_k:J{\it D}(F)\to \mathbb{R}$,
$k=1,2,3$. Thus we can prolong any solution $V\subset A_j$,
$j=1,2$, $\partial V=N_0\DU N_1$ to a solution $Z$, such that
$\partial Z=N_1\DU \{q_0\}\cong N_1$, where $q_0\in A_3$ also. In
other words, $V'\equiv V\bigcup_{N_1}Z$ is an algebraic singular
solution of $\tilde E_1$. (See Fig. \ref{algebraic-singular-solution}.) In fact, we can always find a
solution $\tilde V\subset J^1_2(F)$ of the trivial equation
$J^1_2(F)\subseteq J^1_2(F)$, such that $\partial\tilde
V=N_1\DU\left\{q_0\right\}\cong N_1$, and such that there exists a
disk $D^2_{\epsilon}\subset\tilde V$, centered on $q_0$, with radius
$\epsilon$ and boundary $\partial D_\epsilon^2\equiv
N_\epsilon\subset A_j$. (Let us emphasize that $\dim F=5$, hence we
can embed in $F$ any $2$-dimensional smooth compact manifold. See,
e.g., \cite{HIR}.) Let $\hat V$ be the submanifold of $\tilde V$
such that $\partial\hat V=N_\epsilon\DU N_1$. Then, since $A_j$ is a
strong retract of $J^1_2(F)$, we can deform $\hat V$, obtaining a
solution $V''_\epsilon\subset A_j$, such that $\partial
V''_\epsilon=N_\epsilon\DU N_1$. By taking the limit $\epsilon\to
0$, we can see that the solution $V''_\epsilon$ identifies a
solution $V''\subset A_j$, such that
$V''\DU\left\{q_0\right\}=V'\subset\tilde E_1$ is just an algebraic
singular solution of the singular equation $\tilde E_1$.

\begin{figure}
\centerline{\includegraphics[width=5cm,height=6cm]{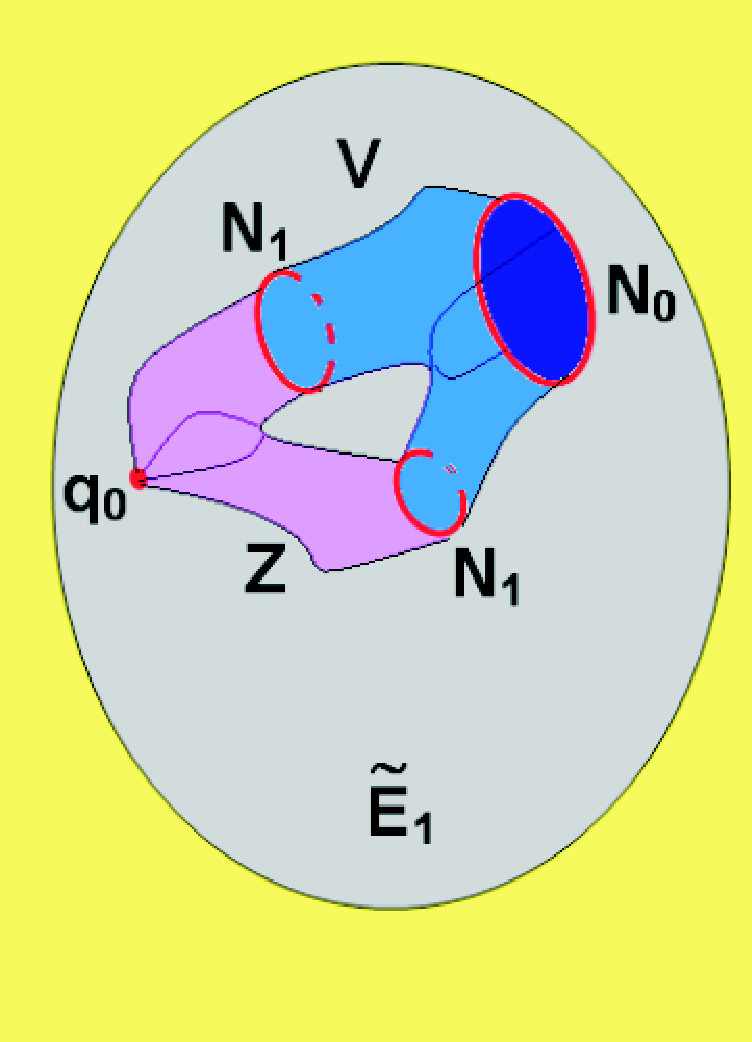}}
\caption{Algebraic singular solution
$V'=V\bigcup_{N_1}Z\subset\tilde E_1\subset J{\it
D}(F)$, passing through a singular point $q_0\in A_3$.}
\label{algebraic-singular-solution}
\end{figure}
Note that the symbol $(g_1)_{q\in A_3}=G_1\equiv TM\otimes F$, i.e.,
$(g_1)_{q\in A_3}$ coincides with the symbol of the trivial PDE
$J{\it D}(F)\subseteq J{\it D}(F)$. In fact the components of the
symbol on $\tilde E_1$, must satisfy the following equations:
\begin{equation}
\left\{
\begin{array}{l}
Y^x_1(u^2_y)^2+Y^y_22u^1_xu^2_y=0;\\
\\
Y^y_12u^1_y+Y^x_22u^2_x=0;\\
\\
Y^x_2u^3_y+Y^y_3[u^2_x+3(u^3_y)^2]=0.\\
\end{array}
\right\}
\end{equation}

These equations are identities on $q\in A_3$. As a by-product, we
get that above considered algebraic singular solutions are, in
general, unstable in finite times.
\end{example}

\begin{definition}\label{extended-crystal-singular-PDE}
We define {\em extended crystal singular PDE}, a singular PDE
$E_k\subset J^k_n(W)$ that splits in irreducible components $A_i$,
i.e., $E_k=\bigcup_i A_i$, where each $A_i$ is an extended crystal
PDE. Similarly we define {\em extended $0$-crystal singular PDE},
(resp. {\em $0$-crystal singular PDE}), an extended crystal singular
PDE where each component $A_i$ is an extended $0$-crystal PDE,
(resp. $0$-crystal PDE).
\end{definition}

\begin{definition}[Algebraic singular solutions of singular PDE's]\label{singular-algebraic-solution}
Let $E_k\subset J^k_n(W)$ be a singular PDE, that splits in
irreducible components $A_i$, i.e., $E_k=\bigcup_iA_i$. Then, we say
that $E_k$ admits an {\em algebraic singular solution} $V\subset
E_k$, if $V\bigcap A_r\equiv V_r$ is a solution  (in the usual
sense) in $A_r$ for at least two different components $A_r$, say
$A_i$, $A_j$, $i\not=j$, and such that one of following conditions
are satisfied: {\em(a)} ${}_{(ij)}E_k\equiv A_i\bigcap
A_j\not=\emptyset$; {\em(b)} ${}^{(ij)}E_k\equiv A_i\bigcup A_j$ is
a connected set, and ${}_{(ij)}E_k=\emptyset$. Then we say that the algebraic singular solution
$V$ is in the case {\em(a)}, {\em weak}, {\em singular} or {\em
smooth}, if it is so with respect to the equation ${}_{(ij)}E_k$. In
the case {\em(b)}, we can distinguish the following situations:
{\em(weak solution):} There is a discontinuity in $V$, passing from
$V_i$ to $V_j$; {\em(singular solution):} there is not discontinuity
in $V$, but the corresponding tangent spaces $TV_i$ and $TV_j$ do
not belong to a same $n$-dimensional Cartan sub-distribution of
$J^k_n(W)$, or alternatively $TV_i$ and $TV_j$ belong to a same
$n$-dimensional Cartan sub-distribution of $J^k_n(W)$, but the
kernel of the canonical projection  $(\pi_{k,0})_*:TJ^k_n(W)\to TW$,
restricted to $V$ is larger than zero; {\em(smooth solution):} there
is not discontinuity in $V$ and the tangent spaces $TV_i$ and $TV_j$
belong to a same $n$-dimensional Cartan sub-distribution of
$J^k_n(W)$ that projects diffeomorphically on $W$ via the canonical
projection $(\pi_{k,0})_*:TJ^k_n(W)\to TW$. Then we say that a
solution passing through a critical zone {\em
bifurcate}.\footnote{Note that the bifurcation does
 not necessarily imply that the tangent planes in the points of $V_{ij}\subset
 V$ to the components $V_i$ and $V_j$, should be different.}
\end{definition}

\begin{definition}[Integral bordism for singular PDE's]\label{integral-bordism-singular-PDE}
Let $N_1, N_2\subset E_k\subset J^k_n(W)$ be two
 $(n-1)$-dimensional admissible closed integral
manifolds. We say that $N_1$ {\em algebraic integral bords} with
$N_2$, if $N_1$ and $N_2$ belong to two different irreducible
components, say $N_1\subset A_i$, $N_2\subset A_j$, $i\not=j$, such
that there exists an algebraic singular solution $V\subset E_k$ with
$\partial V=N_1\DU N_2$.

In the integral bordism group $\Omega_{n-1}^{E_k}$ (resp.
$\Omega_{n-1,s}^{E_k}$, resp. $\Omega_{n-1,w}^{E_k}$) of a singular
PDE $E_k\subset J^k_n(W)$, we call {\em algebraic class} a class
$[N]\in\Omega_{n-1}^{E_k}$, (resp. $[N]\in\Omega_{n-1, s}^{E_k}$,
resp. $[N]\in\Omega_{n-1}^{E_k}$,), with $N\subset A_j$, such that
there exists a closed $(n-1)$-dimensional admissible integral
manifolds $X\subset A_i\subset E_k$, algebraic integral bording with
$N$, i.e., there exists a smooth (resp. singular, resp. weak)
algebraic singular solution $V\subset E_k$, with $\partial V=N\DU
X$.
\end{definition}

\begin{theorem}[Singular integral bordism group of singular
PDE]\label{main-singular1} Let $E_k\equiv\bigcup_i A_i\subset
J^k_n(W)$ be a singular PDE. Then under suitable conditions, {\em
algebraic singular solutions integrability conditions}, we can find
(smooth) algebraic singular solutions bording assigned admissible
closed smooth $(n-1)$-dimensional integral manifolds $N_0$ and $N_1$
contained in some component $A_i$ and $A_j$, $i\not= j$.
\end{theorem}

\begin{proof}
In fact, we have the following lemmas.
\begin{lemma}\label{lemma-main-singular1}
Let $E_k\equiv\bigcup_i A_i\subset J^k_n(W)$ be a singular PDE with
${}_{(ij)}E_k\equiv A_i\bigcap A_j\not=\emptyset$. Let us assume
that $A_i\subset J^k_n(W)$, $A_j\subset J^k_n(W)$ and
${}_{(ij)}E_k\subset J^k_n(W)$ be formally integrable and completely
integrable PDE's such that $\dim A_i>2n+1$, $\dim A_j>2n+1$, $\dim
{}_{(ij)}E_k>2n+1$. Then, one has the following isomorphisms:
\begin{equation}\label{singular-bordism-groups-singular-PDE-1}
    \Omega_{n-1,w}^{A_i}\cong\Omega_{n-1,w}^{A_j}\cong\Omega_{n-1,w}^{{}_{(ij)}E_k}.
\end{equation}
So we can find a weak algebraic singular solution $V\subset E_k$
such that $\partial V=N_0\DU N_1$, with $N_0\subset A_i$,
$N_1\subset A_j$, iff $N_1\in[N_0]$.

Furthermore, if $g_k(A_i)\not=0$, $g_{k+1}(A_i)\not=0$,
$g_k(A_j)\not=0$, $g_{k+1}(A_j)\not=0$, $g_k({}_{(ij)}E_k)\not=0$,
$g_{k+1}({}_{(ij)}E_k)\not=0$, then one has also the following
isomorphisms:
\begin{equation}\label{singular-bordism-groups-singular-PDE-2}
    \Omega_{n-1,s}^{A_i}\cong\Omega_{n-1,s}^{A_j}\cong\Omega_{n-1,s}^{{}_{(ij)}E_k}.
\end{equation}
So we can find a singular algebraic singular solution $V\subset E_k$
such that $\partial V=N_0\DU N_1$, $N_0\subset A_i$,  $N_1\subset
A_j$, iff $N_1\in[N_0]$.
\end{lemma}

\begin{proof}
In fact, under the previous hypotheses one has that we can apply
Theorem \ref{main1} to each component $A_i$, $A_j$ and ${}_{(ij)}E_k$ to
state that all their weak integral bordism groups of dimension
$(n-1)$ are isomorphic to
$\bigoplus_{r+s=n-1}H_r(W;\mathbb{Z}_2)\otimes_{\mathbb{Z}_2}\Omega_s$.

Furthermore, under the above hypotheses on nontriviality of symbols,
we can apply Theorem 2.1 in \cite{PRA10}. So we can state that weak
integral bordism groups are isomorphic to the corresponding singular
ones.
\end{proof}

\begin{lemma}\label{lemma-main-singular1}
Let $E_k=\bigcup_iA_i$ be a $0$-crystal singular PDE. Let
${}^{(ij)}E_k\equiv A_i\bigcup A_j$ be connected, and
${}_{(ij)}E_k\equiv A_i\bigcap A_j\not=\emptyset$. Then
$\Omega_{n-1,s}^{{}^{(ij)}E_k}=0$.\footnote{But, in general, it is
$\Omega_{n-1}^{{}^{(ij)}E_k}\not=0$.}
\end{lemma}

\begin{proof}
In fact, let $Y\subset{}_{(ij)}E_k$ be an admissible closed
$(n-1)$-dimensional closed integral manifold, then there exists a
smooth solution $V_i\subset A_i$ such that $\partial V_i=N_0\DU Y$
and a solution $V_j\subset A_j$ such that $\partial V_j=Y\DU N_1$.
Then, $V=V_i\bigcup_Y V_j$ is an algebraic singular solution of
$E_k$. This solution is singular in general.
\end{proof}

After above lemmas the proof of the theorem can be considered done
besides the algebraic singular solutions integrability conditions.
\end{proof}
\begin{table}[t]\centering
\caption{Nuclear energy producing magnetohydrodynamics equation $\{\widetilde{F}^{(s)}=0; 1\le s\le 7\}:
\widetilde{(MHD)}\subset J{\it D}^2(\widetilde{W})$.}
\label{nuclear-energy-producing-magnetohydrodynamics-equation}
\begin{tabular}{|l|l|}
\hline
{\rm{\footnotesize Maxwell}}&{\rm{\footnotesize $\widetilde{F}^{(1)}\equiv F^{(1)}\equiv {B^k}_{/k} = 0$\hskip 2pt  (no magnetic monopoles)}} \\
&{\rm{\footnotesize $ \widetilde{F}^{(2)}\equiv F^{(2)}\equiv {D^k}_{/k}-\bar\rho  4\pi=0$\hskip 2pt (Gauss's law of electrostatic)}}\\
&{\rm{\footnotesize $ \widetilde{F}^{(3)}\equiv F^{(3)}{}^i\equiv \epsilon ^{ijs} E_{j/s} + {{1}\over{c}} (\partial t \cdot  B^i )=0$\hskip 2pt(Faraday's law)}}\\
&{\rm{\footnotesize $\widetilde{F}^{(4)}\equiv F^{(4)}{}^i\equiv \epsilon
^{ijs} H_{j/s} - {{1}\over{c}} ( \partial t \cdot  D^i ) -
{{4\pi}\over{c}} I^i=0$\hskip 2pt(Ampere's law)}}\\
&{\rm{\footnotesize $ D^k = \epsilon^{ki} E_i$,\hskip 2pt($\epsilon=$dielectric permeability tensor)}} \\
&{\rm{\footnotesize $ \epsilon^{ik}=-\epsilon_0g^{ik}+\bar\epsilon(v_{r/s}+v_{s/r})g^{ri}g^{sk}$}}\\
&{\rm{\footnotesize $ B^k = \mu^{ki} H_i$,\hskip 2pt($\mu=$magnetic permeability tensor)}}\\
&{\rm{\footnotesize $ \mu^{ik}=-\mu_0g^{ik}+\bar\mu(v_{r/s}+v_{s/r})g^{ri}g^{sk}$}}\\
\hline
{\rm{\footnotesize Navier-Stokes}}&{\rm{\footnotesize $\widetilde{F}^{(5)}\equiv F^{(5)}\equiv {v^k}_{/k}=0$ (continuity equation)}}\\
&{\rm{\footnotesize $ \widetilde{F}^{(5)}\equiv
F^{(5)}_\alpha\equiv\dot x^k(\partial x_\alpha.G^j_{jk})+\dot
x^k_\alpha G^j_{jk}+\dot x^s_{s\alpha}=0$,}}\\
&{\rm{\footnotesize \hskip 1cm (first-prolonged
continuity equation)}}\\
&{\rm{\footnotesize $ \widetilde{F}^{(6)}\equiv
F^{(6)}{}^i\equiv\rho {{\delta
v^i}\over{\delta t}} -{P^{ik}}_{/k} - F_{(body)}^i=0$, (motion equation)}}\\
&{\rm{\footnotesize $ \widetilde{F}^{(7)}\equiv \rho
C_v\frac{\delta\theta}{\delta t}+\frac{\delta w}{\delta
t}-\nu(\theta_{/i})_{/k}g^{ik}+S^k_{/k}-[2\chi\dot e{}^{ik}+
\frac{1}{4\pi}(B^iB^k+E^iE^k)]v_{i/k}$}}\\
& {\rm{\footnotesize $+I^iE^jg_{ij}-\rho\bar h=0$, (energy equation)}} \\
&{\rm{\footnotesize $ q^k=-\nu\theta_{/i}g^{ik}$,  (heat flow)}}\\
&{\rm{\footnotesize $ P^{ik}={}^{(rh)}P^{ik}+M^{ik}=-pg^{ik}+\chi(v_{r/s}+v_{s/r})g^{ri}g^{sk}+M^{ik}$}}\\
&{\rm{\footnotesize $ \hskip 14pt=-g^{ik}[p+\frac{1}{8\pi}(B^sB_s+E^sE_s)]+\wp^{ik}$ (full stress tensor)}}\\
&{\rm{\footnotesize $ \wp^{ki}=\chi(v_{r/s}+v_{s/r})g^{ri}g^{sk}+\frac{1}{4\pi}(B^iB^j+E^iE^j)$ (deviatory stress)}}\\
&{\rm{\footnotesize (magnetic stress tensor) $ {}^{(B)}M^{ij}\equiv\frac{1}{4\pi}(-\frac{1}{2}B^kB_k g^{ij}+B^iB^j)$}}\\
&{\rm{\footnotesize (electric stress tensor) $ {}^{(E)}M^{ij}\equiv\frac{1}{4\pi}(-\frac{1}{2}E^kE_k g^{ij}+E^iE^j)$}}\\
&{\rm{\footnotesize (e.m. stress tensor) $ M^{ij}\equiv {}^{(B)}M^{ij}+{}^{(E)}M^{ij}$.}}\\
&{\rm{\footnotesize (body force): $ F_{(body)}^i=-\rho (\partial x_k.f)g^{ki}+\bar\rho E^i+\epsilon^{ijk}I_jB_k$}}\\
\hline \multicolumn{2}{l}{\rm{\footnotesize $ \bar h$ body energy source density.}}\\
\end{tabular}
\end{table}

\begin{example}[Extended crystal singular MHD-PDE's]\label{MHD-PDE}
In the recent paper \cite{PRA23} we introduced a new PDE of the
magnetohydrodynamics encoding the dynamic of anisotropic
incompressible plasmas able to describe nuclear energy production.
This equation is denoted $\widetilde{(MHD)}\subset J{\it
D}^2(\widetilde{W})$ and it is reported in Tab. \ref{nuclear-energy-producing-magnetohydrodynamics-equation}. The fiber bundle
there considered is $\pi:\widetilde{W}\to M$, over the Galilean
space-time $M$, with $\widetilde{W}\equiv M\times
\mathbf{I}\times\mathbf{S}^3\times\mathbb{R}^4$, where $\mathbf{I}$
is an affine $3$-dimensional space (time-like flow velocities space)
and $\mathbf{S}$ is a $3$-dimensional Euclidean vector space. A
section $s=(v,p,\theta,E,H,I,\bar\rho,\bar h)$ represents
flow-velocity, isobaric pressure, temperature, electric vector
field, magnetic vector field, electric current density, electric
charge density, nuclear energy production density. In that paper it
is also proved that equation $\widetilde{(MHD)}$ is formally
integrable, completely integrable and an extended $0$-crystal. Now,
from Theorem \ref{main1} we can see that for any two space-like
admissible Cauchy hypersurfaces $N_0\subset\widetilde{(MHD)}_{t_0}$,
$N_1\subset\widetilde{(MHD)}_{t_1}$, $t_1\not= t_2$, there exists a
(singular) solution $V\subset\widetilde{(MHD)}$, passing through
$N_0$ and $N_1$. The admissibility requires that $N_0$ and $N_1$ are
smooth $3$-dimensional regular manifolds with respect to the
embedding $\widetilde{(MHD)}\subset J^2_4(\widetilde{W})$, locally
satisfying the Cauchy problem and with orientable boundaries
$X_0\equiv\partial N_0$, $X_1\equiv\partial N_1$, bording by means
of a suitable $3$-dimensional time-like integral manifold $P$.
\begin{table}[h]\centering
\caption{Polynomial differential expression for $\theta^2\mathcal{R}$.}
\label{polynomial-differential-expression-for-entropy-production}
\begin{tabular}{|l|}
\hline
{\rm{\footnotesize $ \theta^2\mathcal{R}= \theta\left[\dot
x^a_b\dot x^c_d\mathcal{R}(1)^{bd}_{ac}+ \dot x^a_b\dot
x^c_dH^pH^q\mathcal{R}(2)^{bd}_{acpq}
+ \dot x^a_b\dot x^c\mathcal{R}(3)^{b}_{ac}+ \dot x^a_b\dot x^cH^pH^q\mathcal{R}(4)^{b}_{acpq}\right.$}} \\
{\rm{\footnotesize $ +\dot x^a\dot x^b\mathcal{R}(5)_{ab}+ \dot x^a\dot
x^b\mathcal{R}(6)_{abpq}+\dot x^m_sH^sH^r\mathcal{R}(7)_{mr}+ \dot
 x^m_sE^iE^j\mathcal{R}(8)^{s}_{mij}$}}\\
{\rm{\footnotesize $  + \dot x^tH^qH^r\mathcal{R}(9)_{tqr} +\dot
x^tE^iE^j\mathcal{R}(10)_{tij}+ \dot x^r_u\dot x^a_b\dot
x^c_dH^pH^q\mathcal{R}(11)^{ubd}_{tacpq}+
 \dot x^r_u\dot x^a_b\dot x^cH^pH^q\mathcal{R}(12)^{ub}_{racpq}$}}\\
 {\rm{\footnotesize $\left. + \dot x^r_u\dot x^a\dot x^bH^pH^q\mathcal{R}(13)^{u}_{rabpq}
 +\dot x^w\dot x^a\dot x^bH^pH^q\mathcal{R}(14)_{wabpq}
+\bar h\right]+\theta_j\theta_k\mathcal{R}(15)^{jk}$}}\\
\hline
{\rm{\footnotesize $\mathcal{R}(1)^{bd}_{ac}=\frac{\chi}{\rho}(g^{bd}g_{ac}+\delta^d_a\delta^b_c)$}}\\
{\rm{\footnotesize $\mathcal{R}(2)^{bd}_{acpq}=-\frac{\mu_0\bar\mu}{4\pi\rho}(\delta^d_a\delta^b_q
g_{pc}+\delta^b_q\delta^d_p g_{ca}+g^{db}g_{qa}g_{pc})$}}\\
{\rm{\footnotesize $\mathcal{R}(3)^{b}_{ac}=\frac{4\chi}{\rho}G^b_{ac}$}}\\
{\rm{\footnotesize $\mathcal{R}(4)^{b}_{acpq}=-\frac{\mu_0\bar\mu}{4\pi\rho}(\delta^b_q
g_{pk}G^k_{ac}+\delta^b_q g_{ka}G^k_{pc}+\delta^b_p g_{aj}G^j_{qc}+2g_{pa}G^b_{qc}+2g_{aq}G^b_{pc}+g_{qt}g_{aj}G^t_{kc}G^k_{pb})$}} \\
{\rm{\footnotesize $\mathcal{R}(5)_{ab}=\frac{2\chi}{\rho}G^i_{sa}G^s_{ib}$}}\\
{\rm{\footnotesize$\mathcal{R}(6)_{abpq}=-\frac{\mu_0\bar\mu}{4\pi\rho}(g_{kj}G^j_{qa}G^k_{pb}+g_{qt}G^t_{ka}G^k_{pb}+
2g_{pk}G^j_{qa}G^k_{jb})$}}\\
{\rm{\footnotesize $\mathcal{R}(7)_{mr}=\frac{\mu_0^2}{4\pi\rho}g_{rm}$}}\\
{\rm{\footnotesize $\mathcal{R}(8)^{s}_{mij}=\frac{1}{8\pi\rho}(\delta^s_j g_{im}+\delta^s_i g_{jm})$}}\\
{\rm{\footnotesize $\mathcal{R}(9)_{tqr}=\frac{\mu_0}{4\pi\rho}g_{mr}G^m_{qt}$}}\\
{\rm{\footnotesize $\mathcal{R}(10)_{tij}=\frac{1}{8\pi\rho}(g_{im}G^m_{jt}+g_{jm}G^m_{it})$}}\\
{\rm{\footnotesize $\mathcal{R}(11)^{ubd}_{tacpq}=\frac{\bar\mu^2}{8\pi\rho}(\delta^d_a\delta^b_q
g_{ic}+\delta^b_q\delta^d_p g_{ca}+g^{db}g_{qa}g_{ic})(\delta^i_r\delta^u_p+g^{iu}g_{rp})$}}\\
{\rm{\footnotesize $\mathcal{R}(12)^{ub}_{racpq}=\frac{\bar\mu^2}{4\pi\rho}[\frac{1}{2}(\delta^b_q
g_{ik}G^k_{ac}+ \delta^b_q g_{ka}G^k_{ic}+\delta^b_i
g_{aj}G^j_{qc}+2g_{ia}G^b_{qc}+2g_{aq}G^b_{ic}+g_{qt}g_{aj}G^t_{kc}G^k_{ib})$}}\\
{\rm{\footnotesize $(\delta^i_v\delta^u_p+g^{iu}g_{rp})$}}\\
{\rm{\footnotesize $\hskip 2cm+(\delta^u_a\delta^b_qg_{ir}+\delta^b_q\delta^u_pg_{ra}+g^{ab}g_{qa}g_{ir})G^i_{pc}]$}}\\
{\rm{\footnotesize $\mathcal{R}(13)^{u}_{rabpq}=\frac{\bar\mu^2}{4\pi\rho}[\frac{1}{2}(G_{qa}^jG_{ib}^kg_{kj}+
G_{ka}^tG_{ib}^kg_{qt}+2G_{qa}^jG_{jb}^kg_{ik})
(\delta^i_v\delta^u_p+g^{iu}g_{rp})$}}\\
{\rm{\footnotesize $\hskip 2cm+(\delta^u_q g_{ik}G^k_{rb}+\delta^u_q
g_{kr}G^k_{ib}
+\delta^u_i g_{rj}G^j_{qb}+2g_{ir}G^u_{qb}+2g_{rq}G^k_{ib}+g_{qt}g_{rj}G^t_{kb}G^k_{iu})G^i_{pa}]$}}\\
{\rm{\footnotesize $\mathcal{R}(14)_{wabpq}=\frac{\bar\mu^2}{4\pi\rho}(G^j_{qa}G^k_{ib}g_{kj}+
G^t_{ka}G^k_{ib}g_{qt}+2G^j_{qa}G^k_{jb}g_{ik}+)G^i_{pw}$}}\\
{\rm{\footnotesize $\mathcal{R}(15)^{jk}=\frac{\nu}{\rho}g^{kj}$}}\\
\hline
\end{tabular}
\end{table}

Then the solution $V$ has boundary $\partial
V=N_0\bigcup_{X_0}P\bigcup_{X_1}N_1$. (For details see
Refs.\cite{PRA24, PRA26} There can be also found the explicit
expressions of the differential polynomials defining
$\widetilde{(MHD)}$.) In \cite{PRA26} we proved that we can identify a sub-equation
$\widehat{\underline{(MHD)}}\subset\widetilde{(MHD)}$ such that in
some neighborhood of its points, there exist
entropy-regular-solutions passing from such initial conditions. More
precisely, we can represent $\theta^2\mathcal{R}$ as a polynomial
differential of first order on $J{\it D}(\widetilde{W})$, belonging
to $A[\dot x^a, \dot x^a_b, H^p,E^p,\theta,\theta_j,\bar h]$, where
$A\equiv[[x^1,x^2,x^3]]$. This can be seen taking into account that
\begin{equation}\label{entropy-production-function}
   \theta^2 \mathcal{R}(s)=\theta\left[\frac{1}{\rho 4\pi}(B^iB^j+E^iE^j)\dot e_{ij}+\frac{2\chi}{\rho}\dot e^{ij}\dot e_{ij}+
    \bar h\right]+\frac{\nu}{\rho}(\GRAD\theta)^2.
\end{equation}

and by calculating the corresponding explicit differential
polynomial expression for $\theta^2\mathcal{R}:J{\it
D}(\widetilde{W})\to\mathbb{R}$. (See Tab. \ref{polynomial-differential-expression-for-entropy-production}.)
$\widehat{\underline{(MHD)}}\subset\widetilde{(MHD)}$ is an {\em
extended $0$-crystal singular PDE}, in the sense of Definition
\ref{extended-crystal-singular-PDE}.

Let us resume the proof given in \cite{PRA26}, since this is
necessary to understand the further developments reported below. Let
us define
\begin{equation}\label{entropy-production-set-1}
    \underline{(MHD)}\equiv\left\{q\in\widetilde{(MHD)}|\theta(q)>0,\bar h(q)\ge 0,\mathcal{R}(q)\ge 0\right\}.
\end{equation}
Then, $\underline{(MHD)}$ is a connected, simply connected bounded
domain in $\widetilde{(MHD)}$. In fact, we can split $
\underline{(MHD)}$ in the following way
\begin{equation}\label{split-entropy-production-set}
\underline{(MHD)}={}_{(+,+)}\widetilde{Y}_1\bigcup{}_{(+,+)}\widetilde{Y}_2\bigcup\widetilde{Y}_3\bigcup
{}_{(+,+,+)}\underline{(MHD)},
\end{equation}
 where
\begin{equation}\label{entropy-production-set-02}
{}_{(+,+,+)}\underline{(MHD)}\equiv\left\{q\in\widetilde{(MHD)}|\bar
h(q)>0,\theta(q)>0,\mathcal{R}(q)>0\right\}
\end{equation}
is an open submanifold of $\widetilde{(MHD)}$, and

\begin{equation}\label{entropy-production-set-2}
\left\{\begin{array}{ll}
  {}_{(+,+)}\widetilde{Y}_1& \equiv\left\{
 q\in\widetilde{(MHD)}|\bar h(q)=0,\theta(q)>0,\mathcal{R}(q)>0\right\}\\
 \\
  {}_{(+,+)}\widetilde{Y}_2&\equiv\left\{
 q\in\widetilde{(MHD)}|\mathcal{R}(q)=0,\theta(q)>0,\bar
 h(q)>0\right\}
\end{array}\right.
\end{equation}
are codimension $1$ submanifolds of
 $\widetilde{(MHD)}$. Furthermore,
\begin{equation}\label{entropy-production-set-3}
 \widetilde{Y}_3\equiv\left\{
 q\in{}_{(+)}\widetilde{(MHD)}|\bar
 h(q)=0,\mathcal{R}(q)=0\right\}
\end{equation}
is a codimension $2$
 submanifold of $${}_{(+)}\widetilde{(MHD)}\equiv \left\{q\in \widetilde{(MHD)}|\theta(q)>0\right\}.$$
 Let us study the integrability properties of such submanifolds of
 $\widetilde{(MHD)}$. First note that since
 ${}_{(+,+,+)}\underline{(MHD)}$ is an open submanifold
 in $\widetilde{(MHD)}$ and in $J{\it
 D}^2({}_{(+,+)}\widetilde{W})$, where $${}_{(+,+)}\widetilde{W}\equiv M\times \mathbf{I}\times\mathbb{R}\times\mathbb{R}^+
 \times\mathbf{S}^3\times\mathbb{R}\times\mathbb{R}^+,$$ it follows that if
 $q\in{}_{(+,+,+)}\underline{(MHD)}$ there exist solutions belonging to ${}_{(+,+,+)}\underline{(MHD)}$,
 i.e., entropy-regular-solutions passing through $q$, just follows from the fact
 that ${}_{(+,+,+)}\underline{(MHD)}$ is an open
 submanifold of $\widetilde{(MHD)}$. So
 ${}_{(+,+,+)}\underline{(MHD)}$ is completely integrable. Furthermore,
 since
\begin{equation}
\left\{
\begin{array}{ll}
   {}_{(+,+,+)}\underline{(MHD)}_{+r}& =J{\it D}^r({}_{(+,+,+)}\underline{(MHD)})\bigcap
 J{\it D}^{2+r}({}_{(+,+)}\widetilde{W})\\
 \\
 &\subset J{\it D}^r(\widetilde{(MHD)})\bigcap J{\it D}^{2+r}(\widetilde{W})
  =\widetilde{(MHD)}_{+r}
\end{array}\right.
\end{equation}
 it follows that
 ${}_{(+,+,+)}\underline{(MHD)}$ is also formally integrable.
Let us consider the integrability properties of
${}_{(+,+)}\widetilde{Y}_1\subset\partial
 \underline{(MHD)}$. If $q\in{}_{(+,+)}\widetilde{Y}_1$
 it follows that $q\in{}_{(+)}\widehat{(MHD)}$. Since this last equation is formally
 integrable and completely integrable, and
 ${}_{(+,+)}\widetilde{Y}_1$ is an open submanifold of
 ${}_{(+)}\widehat{(MHD)}$, it follows that also ${}_{(+,+)}\widetilde{Y}_1$ is formally integrable and completely integrable.

 Let us, now, study the integrability
 properties of
 \begin{equation}
 {}_{(+,+)}\widetilde{Y}_2\subset J{\it
 D}^{2}({}_{(+,+)}\widetilde{W}):
 \left\{\mathcal{R}=0,\widetilde{F}^{(s)}=0\right\}.
\end{equation}
 One can see
 that this equation is not formally integrable, but becomes so if we
 add the first prolongation of $\mathcal{R}=0$. So we can prove that
 the following equation

 \begin{equation}\label{formally-integrable-2}
\widehat{{}_{(+,+)}\widetilde{Y}_2}\subset J{\it
D}^{2}({}_{(+,+)}\widetilde{W}):
 \left\{\mathcal{R}=0,\mathcal{R}_\alpha=0,\widetilde{F}^{(s)}=0\right\}
\end{equation}
is formally integrable and completely integrable. Similar considerations can be made on the part $\widetilde{Y}_3$,
and we can identify, the corresponding formally integrable PDE
$\widehat{\widetilde{Y}_3}$. We skip on the details. By conclusion
we get that
 \begin{equation}\label{formally-integrable-full}
\widehat{\underline{(MHD)}}\equiv{}_{(+,+)}\widetilde{Y}_1\bigcup\widehat{{}_{(+,+)}\widetilde{Y}_2}
\bigcup\widehat{\widetilde{Y}_3}
\bigcup{}_{(+,+,+)}\underline{(MHD)}\subset\widetilde{(MHD)}\end{equation}
is the formally integrable and completely integrable constraint in
$\widetilde{(MHD)}$, where for any initial condition, passes an
entropy-regular-solution for $\widetilde{(MHD)}$. Then we can apply
Theorem \ref{main1} to conclude that $\widehat{\underline{(MHD)}}$
is an extended crystal singular PDE. Furthermore, taking into
account that it results
$$\Omega^{\widehat{\underline{(MHD)}}}_{3,w}\cong\Omega^{\widehat{\underline{(MHD)}}}_{3,s}
\cong\Omega_3=0,$$ we get that $\widehat{\underline{(MHD)}}$ is also
an extended $0$-crystal singular PDE. This assures that for any two
space-like admissible Cauchy integral manifolds
$N_1\subset\widehat{\underline{(MHD)}}_{t_1}$,
$N_2\subset\widehat{\underline{(MHD)}}_{t_2}$, $t_1\not=t_2$, such
that if both $N_i$, $i=1,2$, belong to the same component, in the
split given in {\em(\ref{formally-integrable-full})}, there exist
(singular) entropy-regular-solutions $V$ such that $\partial
V=N_1\bigcup_{X_1}P\bigcup_{X_2} N_2$, and such that their
boundaries $X_1\equiv\partial N_1$, $X_2\equiv\partial N_2$, should
be orientable and propagating with an admissible $3$-dimensional
integral manifold $P$. (For details see \cite{PRA24}.)\footnote{Let us emphasize that admissible Cauchy integral manifolds can be found in each component of $\widehat{\underline{(MHD)}}$, thanks to Lemma \ref{cauchy-criteria}. More precisely, the proceeding followed for the Navier-Stokes equation in Example \ref{ex12}, to solve Cauchy problems there, can be applied also to $\widehat{\underline{(MHD)}}$. In fact, envelopment solutions can be built also for all the components of this last equation, since they are formally integrable and completely integrable PDE's.}
The stability properties of $\widehat{\underline{(MHD)}}$ and its
solutions, can be studied by utilizing our recent geometric theory
on the stability of PDE's \cite{PRA22, PRA23, PRA24, PRA25, PRA26}.
More precisely, $\widehat{\underline{(MHD)}}$ is a {\em functionally
stable singular PDE}, in the sense that it slpits in components that
are functionally stable PDE's. Furthermore, smooth
entropy-regular-solutions, i.e., smooth solutions of
$\widehat{\underline{(MHD)}}$, do not necessitate to be stable.
However, all they can be stabilized and the {\em stable extended
crystal singular PDE} of $\widehat{\underline{(MHD)}}$, i.e., a
singular PDE slpits in components that are stable extended crystal
PDE's. This last is just
${}^{(S)}\widehat{\underline{(MHD)}}=\widehat{\underline{(MHD)}}_{(+\infty)}$.
There, all smooth entropy-regular-solutions, belonging to only one
of the components, in the split representation
{\em(\ref{formally-integrable-full})} are stable at finite times.

Let us, now, find global solutions of $\widehat{\underline{(MHD)}}$,
crossing the singular sets corresponding to states without nuclear
energy production, to ones where $\bar h>0$, i.e., passing from
${}_{(+,+)}\tilde Y_1$ to ${}_{(+,+,+)}\underline{(MHD)}$.
Let $V\subset\widetilde{(MHD)}$ be a time-like solution such that
the following conditions are satisfied:

{\em (i)} $V$ is a regular entropy solution;

{\em (ii)} $V\bigcap {}_{(+,+)}\tilde Y_1\equiv
{}_{(+)}V\not=\emptyset$;

{\em (iii)} $V\bigcap {}_{(+,+,+)}\underline{(MHD)}\equiv
{}_{(+,+)}V\not=\emptyset$.

{\em (iv)} $\partial V=N_0\bigcup P\bigcup N_1$, $N_0\subset
{}_{(+)}V$, $N_1\subset {}_{(+,+)}V$, $\partial {}_{(+)}V=N_0\bigcup
{}_{(+)}P\bigcup Y$.

Then we can give the following split surgery representation of $V$:

\begin{equation}\label{crossing-nuclear-critical-zone}
    V={}_{(+)}V\bigcup_Y Z
\end{equation}
with $Z=Y\bigcup {}_{(+,+)}V$. We call a solution $V$ of
$\widetilde{(MHD)}$, such that holds the surgery property given in
{\em(\ref{crossing-nuclear-critical-zone})} just a {\em crossing
nuclear critical zone solution}. Then, for any two admissible
space-like Cauchy hypersurfaces
$N_1\subset\widehat{\underline{(MHD)}}_{t_1}$,
$N_2\subset\widehat{\underline{(MHD)}}_{t_2}$, $t_1\not=t_2$, such
that the following properties are satisfied:

{\em (a)} $N_1\subset{}_{(+,+)}\tilde Y_1$;

{\em (b)} $N_2\subset {}_{(+,+,+)}\underline{(MHD)}$;

there exists a crossing nuclear critical zone solution $V\subset \widetilde{\underline{(MHD)}}$,
such that $\partial V=N_1\bigcup P\bigcup N_2$, where $P$ is a
time-like admissible $3$-dimensional integral manifold of
$\widetilde{\underline{(MHD)}}$. In general such a solution is a
singular solution. Therefore, under above condition of
admissibility, one has algebraic classes in
$\Omega_{2,s}^{\widehat{\underline{(MHD)}}}$. Furthermore, $V$ is
represented by a smooth integral manifold at $Y$, with respect to
the split given in {\em(\ref{crossing-nuclear-critical-zone})}, if
the tangent space $T{}_{(+)}V|_Y=TZ|_Y\subset
\mathbf{E}_2\widetilde{(MHD)}|_Y$, where
$\mathbf{E}_2\widetilde{(MHD)}$ is the Cartan distribution of
$\widetilde{(MHD)}$. In fact the dimension of the Cartan distribution
$\mathbf{E}_2$ of $\widehat{\underline{(MHD)}}$ in the points $q\in
{}_{(+,+,+)}\underline{(MHD)}$ is higher than in the points
$q\in{}_{(+,+)}\tilde Y_1)$. This can be seen by direct computation.
Let us denote by
\begin{equation}\label{cartan0}
\zeta=X^\alpha[\partial x_\alpha+\sum_{0\le|\beta|\le
1}y^j_{\alpha\beta}\partial y_j^\beta]+Z^j_{\beta_1\beta_2}\partial
y_j^{\beta_1\beta_2}
\end{equation}
the generic vector field of the Cartan distribution
$\mathbf{E}_2(\widetilde{W})$, where $$\{y^j\}_{1\le j\le 16}=\{
v^i,p,\theta,E^i,H^i,I^i,\bar\rho,\bar h\}$$ are the vertical
coordinates of the fiber bundle $\pi:\widetilde{W}\to M$. The Cartan
distribution on the boundary ${}_{(+,+)}\tilde
Y_1\subset\widehat{\underline{(MHD)}}$ is given by the vector fields
{\em(\ref{cartan0})} such that the following equations are
satisfied: $\zeta.\bar F^I=0$, where $\bar F^I$ are the functions
defining equation ${}_{(+,+)}\tilde Y_1$. These are just the
functions defining ${}_{(+,+,+)}\underline{(MHD)}$, i.e.,
$\widetilde{F}^{(s)}$, $1\le s\le 7$, but with the condition $\bar
h=0$. So we get
 \begin{equation}\label{dim-cartan}
 \scalebox{0.8}{$\left\{\dim \mathbf{E}_2(\widetilde{MHD})= \dim \mathbf{E}_2({}_{(+,+,+)}\underline{(MHD)})
        =148>\dim \mathbf{E}_2({}_{(+,+)}\tilde Y_1)=138.\right\}.$}
 \end{equation}
Note that even if the points of ${}_{(+,+)}\tilde Y_1$ can be
considered singular one, with respect to the Cartan distribution of
$\widehat{\underline{(MHD)}}$, the embedding
$\widehat{\underline{(MHD)}}\subset\widetilde{(MHD)}$ allows us to
prolong a solution from ${}_{(+,+)}\tilde Y_1$ to
${}_{(+,+,+)}\underline{(MHD)}$, according to Theorem
\ref{main-singular1} and surgery representation
{\em(\ref{crossing-nuclear-critical-zone})}. Then such solution $V$
in general bifurcates along $Y$. By summarizing, we can say that the
crossing nuclear critical zone solution
$V\subset\widehat{\underline{(MHD)}}$ is represented by an integral
manifold $V\subset\widehat{\underline{(MHD)}}$ that is smooth in a
neighborhood of $X\subset Y$, iff for its tangent space $TV$, the
following condition is satisfied: $T{}_{(+)}V|_X=TZ|_X\subset
\mathbf{E}_2(\widetilde{MHD})|_X$.

Now, in order to surgery a smooth solution
$V'\subset{(+,+,+)}\underline{(MHD)}$, passing through a compact
smooth space-like $3$-dimensional manifold $N_1$, with
$X_1\equiv\partial N_1$ orientable, with suitable smooth solutions
of ${}_{(+,+)}\tilde Y_1$ it is enough that the following conditions
should be satisfied:
\begin{equation}\label{surgery-smoothness-condition-MHD}
    \left\{
    \begin{array}{l}
\lim_{\bar h\to 0}V'={}_{(+,+)}V\subset{{}_(+,+,+)}\underline{(MHD)}\\
\\
      T{}_{(+,+)}V\subset \mathbf{E}_2({(+,+,+)}\underline{(MHD)})_{\bar
      h=0}.\\
      \end{array}\right.
\end{equation}

Since ${}_{(+,+)}V$ is a smooth solution we can prolong it to
$\infty$, obtaining a solution ${}_{(+,+)}V^{(\infty)}\subset
{}_{(+,+,+)}\underline{(MHD)}_{(+\infty)}$. There it identifies an
horizontal $4$-plane $\mathbf{H}\subset
\mathbf{E}_2({}_{(+,+,+)}\underline{(MHD)}_{(+\infty)})$, contained
in the Cartan distribution of
${}_{(+,+,+)}\underline{(MHD)}_{(+\infty)}$. Then $\mathbf{H}$
identifies also an horizontal $4$-plane, that we continue to denote
with $\mathbf{H}$, in the Cartan distribution of
$\widetilde{(MHD)}_{+\infty}$, that on $({}_{(+,+)}\tilde
Y_1)_{+\infty}$, coincides with the Cartan distribution of this last
equation. Then we can smoothly prolong ${}_{(+,+)}V^{(\infty)}$
into $({}_{(+,+)}\tilde Y_1)_{+\infty}$, identifying there a smooth
solution ${}_{(+)}V^{(\infty)}$. The projection of this algebraic
singular solution on ${}_{(+,+)}\tilde
Y_1\bigcup{}_{(+,+,+)}\underline{(MHD)}$ identifies a smooth
solution $V$ that has the split surgery property
{\em(\ref{crossing-nuclear-critical-zone})}. Therefore, it is enough
to prove that solutions with the property
{\em(\ref{surgery-smoothness-condition-MHD})} exist in
${}_{(+,+,+)}\underline{(MHD)}$. Now, we can see that a Cartan
vector field $\zeta$ of ${}_{(+,+,+)}\underline{(MHD)}$ is given in
{\em(\ref{cartan0})} and subject to the condition $\zeta.\tilde
F^{(s)}=0$. Then conditions
{\em(\ref{surgery-smoothness-condition-MHD})} are satisfied iff
$X^\alpha {\bar h}_\alpha=0$. Since this condition can be satisfied
for suitable functions $X^\alpha$ on
${}_{(+,+,+)}\underline{(MHD)}$, we recognize that in the set of
smooth solutions of ${}_{(+,+,+)}\underline{(MHD)}$ exist solutions
that smoothly surgery with smooth solutions of ${}_{(+,+)}\tilde
Y_1$. Such solutions are not finite times stable in
$\widehat{\underline{(MHD}}$ since there the symbol is not trivial
in each components ${}_{(+,+)}\tilde Y_1$ and
${}_{(+,+,+)}\underline{(MHD)}$. However, by using their formal
integrability and complete integrability properties, we can state
that in $(\widehat{\underline{(MHD}})_{+\infty}$ such solutions are
finite times stable. So $\widehat{\underline{(MHD}})_{+\infty}$ is
the {\em stable extended crystal singular PDE} associated to
$\widehat{\underline{(MHD}}$.
\end{example}

\vfill\eject
\section{\bf APPENDIX A - THE AFFINE CRYSTALLOGRAPHIC GROUP
TYPES {\boldmath$[G(3)]$} AND {\boldmath$[G(2)]$}}

\begin{table}[h]
\caption{The 32 crystallographic point groups of the space-group types $[G(3)]$.}
\label{the-32-crystallographic-point-groups-of-the-space-group-types}
\scalebox{0.8}{$\begin{tabular}{|l|l|l|}
\hline
{\bf{\footnotesize Type}}&{\bf{\footnotesize Schoenflies Symbol}}&{\bf{\footnotesize  International Symbol}}\\
\hline
{\rm{\footnotesize  nonaxial (2)}}&{\rm{\footnotesize $ C_i=S_2$, $ C_s=C_{1h}$}}&{\rm{\footnotesize $ \bar 1$, $ m$}}\\
\hline
 {\rm{\footnotesize cyclic (5)}}&{\rm{\footnotesize $ C_n$, $ n=1,2,3,4,6$}}&{\rm{\footnotesize  $ n$, $n=1,2,3,4,6$}}\\
\hline
{\rm{\footnotesize  cyclic with horizontal planes (4)}}&{\rm{\footnotesize $
C_{nh}$, $ n=2,3,4,6$}}&{\rm{\footnotesize $ 2/m$,
$ \bar 6$, $ 4/m$, $ 6/m$}}\\
\hline
{\rm{\footnotesize  cyclic with vertical planes (4)}}&{\rm{\footnotesize $
C_{nv}$, $ n=2,3,4,6$}}&{\rm{\footnotesize $ mm2$,
$ 3m$, $ 4mm$, $ 6mm$}}\\
\hline
{\rm{\footnotesize  dihedral (4)}}&{\rm{\footnotesize $ D_{n}$, $
n=2,3,4,6$}}&{\rm{\footnotesize $ 222$, $ 32$, $
422$, $ 622$}}\\
\hline
{\rm{\footnotesize  dihedral with horizontal planes (4)}}&{\rm{\footnotesize $
D_{nh}$, $ n=2,3,4,6$}}&{\rm{\footnotesize $ mmm$,
$ \bar 6m2$, $ 4/mmm$, $ 6/mmm$}}\\
\hline
{\rm{\footnotesize  dihedral with planes between axes
(2)}}&{\rm{\footnotesize $ D_{nd}$, $ n=2,3$}}&{\rm{\footnotesize $
\bar 42m$, $ \bar 3m$}}\\
\hline
 {\rm{\footnotesize improper rotation (2)}}&{\rm{\footnotesize $ S_{n}$, $ n=4,6$}}&{\rm{\footnotesize $ \bar 4$, $ \bar 3$}}\\
\hline
{\rm{\footnotesize  cubic groups (5)}}&{\rm{\footnotesize  $T$, $
T_h$, $ T_d$, $ O$, $ O_h$}}&{\rm{\footnotesize
$ 23$, $ m\bar 3$, $ \bar 43m$, $ 432$, $ m\bar 3 m$}}\\
\hline \multicolumn{3}{l}{\rm{\footnotesize $ n$\hskip 2pt
($ n=1,2,3,4,6$):
rotations of $ 2\pi/n$ about a symmetry axis.}}\\
\multicolumn{3}{l}{\rm{\footnotesize $ \bar n$\hskip 2pt
($ n=1,2,3,4,6$):
rotation $ n$ composed with inversion about symmetry centre; $ m=\bar 2$.}}\\
\multicolumn{3}{l}{\rm{\footnotesize $
C_n=\mathbb{Z}_n\cong\{0,1,2,3,\cdots,n-1\}$ (cyclic
abelian groups); $ C_i=S_2$; $ C_s=m$; $ C_{3i}=S_6$}}\\
\multicolumn{3}{l}{\rm{\footnotesize $ D_n$: group with an
$ n$-fold axis plus a two-fold axis perpendicular to
that axis.}}\\
\multicolumn{3}{l}{\rm{\footnotesize $ D_n$ is a non abelian
group for $ n>2$. The group order is $ 2n$.}}\\
\multicolumn{3}{l}{\rm{\footnotesize $ D_{nh}$: $
D_n$ with a mirror plane symmetry perpendicular to the
$n$-fold axis.}}\\
\multicolumn{3}{l}{\rm{\footnotesize $ D_{nv}$: $
D_n$ with mirror plane symmetries parallel to the $
n$-fold axis.}}\\
\multicolumn{3}{l}{\rm{\footnotesize $ O$: symmetry group of the
octahedron. The group order is $ 24$. One has the isomorphism $ O\cong T_d$.}}\\
\multicolumn{3}{l}{\rm{\footnotesize $ O_h$: $ O$
with improper operations (those that change orientation). The
group order
is $ 48$.}}\\
\multicolumn{3}{l}{\rm{\footnotesize $ T$: symmetry group of the
tetrahedron, isomorphic to the alternating group $ A_4$. The group order is $ 12$.}}\\
\multicolumn{3}{l}{\rm{\footnotesize \hskip 2cm The group order is $ 12$.}}\\
\multicolumn{3}{l}{\rm{\footnotesize $ T_d$: $\scriptstyle T$
with improper operations.
Non abelian group of order $ 24$.}}\\
\multicolumn{3}{l}{\rm{\footnotesize $ T_h$: $ T$
with
the addition of an inversion.}}\\
\multicolumn{3}{l}{\rm{\footnotesize $ A_n$: group of even
permutations on a set of length $ n$. The group order is
$\frac{n!}{2}$.}}
\end{tabular}$}
\end{table}

\begin{table}[h]
\caption{Diehdral groups $D_{2m}$, $m\ge 1$.}
\label{diehdral-groups}
\scalebox{0.8}{$\begin{tabular}{|l|l|}
\hline
{\rm{\footnotesize $m$}}&{\rm{\footnotesize $D_{2m}$}}\\
\hline
{\rm{\footnotesize $1$}}&{\rm{\footnotesize $ D_{2}\cong\mathbb{Z}_2$}}\\
\hline
{\rm{\footnotesize $ m\ge 2$}}&{\rm{\footnotesize $ D_{2m}\cong\mathbb{Z}_m\rtimes\mathbb{Z}_2$}}\\
&{\rm{\footnotesize  The generator of $ \mathbb{Z}_2$ acts on
$ \mathbb{Z}_m$  as multiplication by $-1$}}\\
\hline
{\rm{\footnotesize $ 2$}}&{\rm{\footnotesize $
D_4\cong \mathbb{Z}_2\times\mathbb{Z}_2$ (Klein four-group)}}\\
\hline
\multicolumn{2}{l}{\rm{\footnotesize  Warn that in crystallography diehdral groups are usually denoted by $ D_m=D_{2m}$, $ m\ge 2$.}}\\
\multicolumn{2}{l}{\rm{\footnotesize This is just the notation used in Tab. \ref{the-32-crystallographic-point-groups-of-the-space-group-types}.}}\\
\end{tabular}$}
\end{table}

\begin{table}[h]
\caption{The 230 affine crystallographic space-group types $[G(3)]$.}
\label{230-affine-crystallographic-space-group-types-3-dim}
\scalebox{0.8}{$\begin{tabular}{|l|l|l|}
\hline
{\bf{\footnotesize Syngonies (7)}}&{\bf{\footnotesize Geometric classes (32=Point Groups)}}&{\bf{\footnotesize  Bravais types (14)}}\\
\hline
{\rm{\footnotesize Triclinic}}&{\rm{\footnotesize $ C_i=S_2 (1)$, $ C_1(1)$}}&{\rm{\footnotesize $ 2P$}}\\
\hline
{\rm{\footnotesize Monoclinic}}&{\rm{\footnotesize $ C_2(3)$, $
C_{1h}(4)$, $ C_{2h}(6)$}}&{\rm{\footnotesize $ 8P$, $ 5C$}}\\
\hline
{\rm{\footnotesize  Orthorhombic}}&{\rm{\footnotesize $ D_2(9)$, $
C_{2v}(22)$, $ D_{2h}(28)$}}&{\rm{\footnotesize $ 30P$, $ 15C$, $ 9I$, $ 5F$}}\\
\hline
{\rm{\footnotesize Tetragonal}}&{\rm{\footnotesize $ C_4(6)$, $
S_4(2)$, $ C_{4h}(6)$, $ D_4(10)$, $ C_{4v}(12)$, $ D_{2d}(12)$, $ D_{4h}(20)$}}&{\rm{\footnotesize $ 49P$,
 $ 19I$}}\\
\hline
{\rm{\footnotesize Trigonal}}&{\rm{\footnotesize $ C_6(4)$, $
S_6(2)$, $ D_3(7)$, $ C_{3v}(6)$, $ D_{3d}(6)$}}&{\rm{\footnotesize $ 18P$, $ 7R$}}\\
\hline
{\rm{\footnotesize Hexagonal}}&{\rm{\footnotesize $C_3(6)$, $
C_{3h}(1)$, $ C_{6h}(2)$, $ D_6(6)$, $ C_{6v}(4)$, $ D_{3h}(4)$, $
D_{6h}(4)$ }}&{\rm{\footnotesize $27P$}}\\
\hline
{\rm{\footnotesize Cubic}}&{\rm{\footnotesize $ T(5)$, $ T_h(7)$,
$ O(8)$, $ T_d(6)$, $O_h(10)$}}&{\rm{\footnotesize $ 15P$, $ 11F$, $ 9I$}}\\
\hline \multicolumn{3}{l}{\rm{\footnotesize (Crystal class systems
=Syngonies): 7. (Triclinic=Anorthic), (Trigonal=Rhombohedral).}}\\
\multicolumn{3}{l}{\rm{\footnotesize Bravais lattice centering types:
$ P=$primitive, $ I=$body, $F=$face, $ A/B/C=$side, $ R=$rhombohedral.}}\\
\multicolumn{3}{l}{\rm{\footnotesize (For the body-centerd case
($ B=I$), the infinite translation group is
$ \{\mathbb{Z}^3, \mathbb{Z}^3+(\frac{1}{2},\frac{1}{2},\frac{1}{2})\}$.)}}\\
\end{tabular}$}
\end{table}

\begin{table}[h]
\caption{The 17 affine crystallographic space-group types $[G(2)]$.}
\label{the-17-affine-crystallographic-space-group-types-2-dim}
\scalebox{0.8}{$\begin{tabular}{|l|l|l|}
\hline
{\bf{\footnotesize Syngony (5)}}&{\bf{\footnotesize International Symbols}}&{\bf{\footnotesize Geometric classes (10=Point Groups)}}\\
\hline
 {\rm{\footnotesize Oblique}}&{\rm{\footnotesize $ p1$, $p2$}}&{\rm{\footnotesize $ \mathbb{Z}_1$, $\mathbb{Z}_2$}}\\
\hline
{\rm{\footnotesize  Rectangular}}&{\rm{\footnotesize $ pm$, $
pg$, $ cm$, $ pmm$}}&{\rm{\footnotesize $ D_1$, $ D_1$, $ D_1$, $ D_2$}} \\
\hline
{\rm{\footnotesize  Rhombic}}&{\rm{\footnotesize $ pmg$, $
pgg$, $ cmm$}}&{\rm{\footnotesize $ D_2$, $ D_2$, $ D_2$}}\\
\hline
{\rm{\footnotesize Square}}&{\rm{\footnotesize $ p4$, $p4m$,
$ p4g$}}&{\rm{\footnotesize $ \mathbb{Z}_4$, $ D_4$,  $ D_4$}}\\
\hline
{\rm{\footnotesize Trigonal}}&{\rm{\footnotesize $ p3$, $ p3m1$,
$ p31m$}}&{\rm{\footnotesize $ \mathbb{Z}_3$, $ D_3$, $D_3$}}\\
{\rm{\footnotesize Hexagonal}}&{\rm{\footnotesize $ p6$, $ p6m$}}&{\rm{\footnotesize $ \mathbb{Z}_6$, $ D_6$}}\\
\hline
\multicolumn{3}{l}{\rm{\footnotesize  All planar crystallographic groups
$ G(2)$ are subgroups of $ p4m$ or $ p6m$, or both. (See Appendix D.)}}\\
\multicolumn{3}{l}{\rm{\footnotesize $ p=$primitive, $
c=$centered, $ m=$mirror plane, $ g=$glide reflection.}}\\
\multicolumn{3}{l}{\rm{\footnotesize A glide reflection is an isometry of the
Euclidean plane that combines reflection in a line}}\\
\multicolumn{3}{l}{\rm{\footnotesize with a translation along that line.}}\\
\end{tabular}$}
\end{table}

\begin{table}[h]
\caption{The 4 holohedries (lattice symmetries) in $[G(2)]$.}
\label{the-4-holohedries-in-2-dim}
\scalebox{0.8}{$\begin{tabular}{|l|l|}
\hline
{\bf{\footnotesize Holohedry}}&{\bf{\footnotesize Lattice system}}\\
\hline
{\rm{\footnotesize $ 2$}}&{\rm{\footnotesize  monoclinic=oblique}}\\
\hline
{\rm{\footnotesize $ 2mm$}}&{\rm{\footnotesize orthorohmbic=$\{$rectangular,rhombic$\}$}}\\
\hline
{\rm{\footnotesize $ 4mm$}}&{\rm{\footnotesize  tetragonal=square}}\\
\hline
{\rm{\footnotesize $ 6mm$}}&{\rm{\footnotesize  hexagonal=$\{$trigonal,hexagonal$\}$}}\\
\hline
\end{tabular}$}
\end{table}

\begin{table}[h]
\caption{The 7 holohedries (lattice symmetries) in $[G(3)]$.}
\label{the-7-holohedries-in-3-dim}
\scalebox{0.8}{$\begin{tabular}{|l|l|}
\hline
{\bf{\footnotesize Holohedry}}&{\bf{\footnotesize Lattice system}}\\
\hline
{\rm{\footnotesize $ \bar 1$}}&{\rm{\footnotesize triclinic}}\\
\hline
{\rm{\footnotesize $ 2/m$}}&{\rm{\footnotesize  monoclinic}}\\
\hline
{\rm{\footnotesize $ mmm$}}&{\rm{\footnotesize orthorohmbic}}\\
\hline
{\rm{\footnotesize $ 4/mmm$}}&{\rm{\footnotesize  tetragonal (square)}}\\
\hline
{\rm{\footnotesize $ \bar 3m$}}&{\rm{\footnotesize trigonal (rombohedric)}}\\
\hline
{\rm{\footnotesize $6/mmm$}}&{\rm{\footnotesize hexagonal}}\\
\hline
{\rm{\footnotesize $ m\bar 3m$}}&{\rm{\footnotesize cubic}}\\
\hline
\end{tabular}$}
\end{table}
\begin{figure}[b]\centering
\includegraphics[height=3cm]{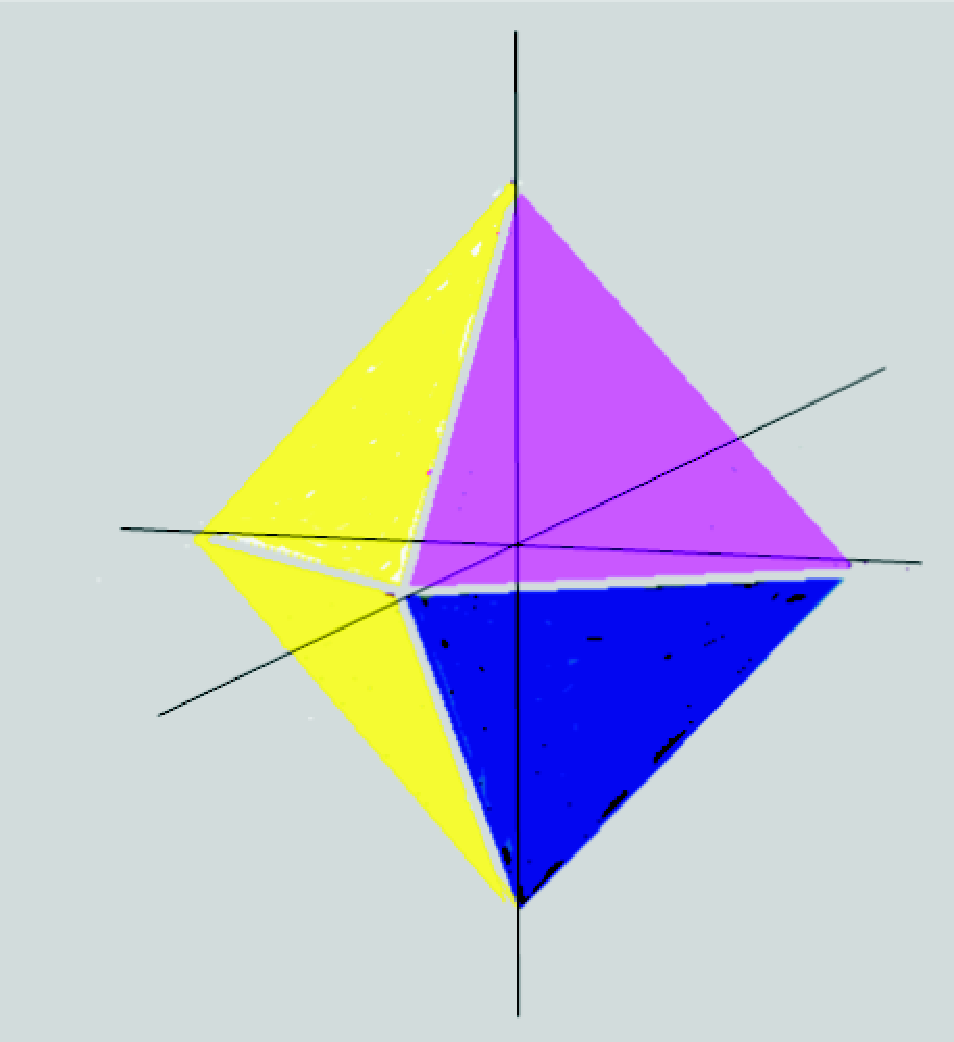}\includegraphics[height=3cm]{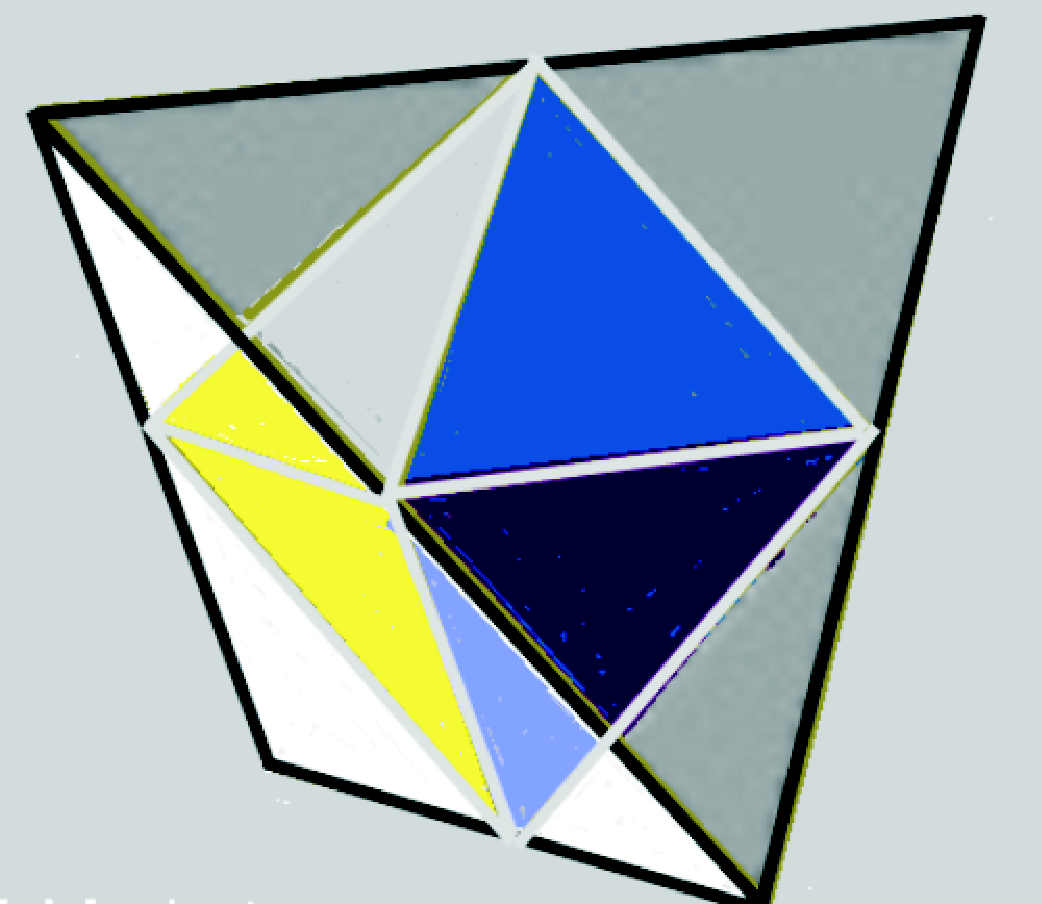}
\caption{Octahedron and Tetrahedron $3$-chains in
$\mathbb{R}^3$. (The tetrahedron is a $3$-chain in
$\mathbb{R}^3$, identified with a regular $4$-faced polyhedron,
where each face is an equilateral triangle. It can be derived from
the octahedron by extending alternate faces until they meet.)}
\label{octahedron-and-tetrahedron}
\end{figure}

\vfill\eject
\section{\bf APPENDIX B - THE SUBGROUPS OF THE AFFINE CRYSTALLOGRAPHIC POINT-GROUP
TYPES {\boldmath$[G(3)]$}}

$\begin{tabular}{|c|c|c|}
\hline \multicolumn{3}{|c|}{\bsmall $\scriptstyle C_i=S_2$}\\
\hline
\bsmall Subgroup&\bsmall Order&\bsmall Index\\
\hline
$\scriptstyle -1$& $\scriptstyle 2$&$\scriptstyle 1$\\
\hline
$\scriptstyle 1$& $\scriptstyle 1$&$\scriptstyle 2$\\
\hline
\end{tabular}
\begin{tabular}{|c|c|c|}
\hline \multicolumn{3}{|c|}{\bsmall $\scriptstyle C_s=C_{1h}$}\\
\hline
\bsmall Subgroup&\bsmall Order&\bsmall Index\\
\hline
$\scriptstyle m$& $\scriptstyle 2$&$\scriptstyle 1$\\
\hline
$\scriptstyle 1$& $\scriptstyle 1$&$\scriptstyle 2$\\
\hline
\end{tabular}
\begin{tabular}{|c|c|c|}
\hline \multicolumn{3}{|c|}{\bsmall $\scriptstyle C_1$}\\
\hline
\bsmall Subgroup&\bsmall Order&\bsmall Index\\
\hline
$\scriptstyle 1$& $\scriptstyle 1$&$\scriptstyle 1$\\
\hline
\end{tabular}$

$\begin{tabular}{|c|c|c|}
\hline \multicolumn{3}{|c|}{\bsmall $\scriptstyle C_2$}\\
\hline
\bsmall Subgroup&\bsmall Order&\bsmall Index\\
\hline
$\scriptstyle 2$& $\scriptstyle 2$&$\scriptstyle 1$\\
\hline
$\scriptstyle 1$& $\scriptstyle 1$&$\scriptstyle 2$\\
\hline
\end{tabular}
\begin{tabular}{|c|c|c|}
\hline \multicolumn{3}{|c|}{\bsmall $\scriptstyle C_3$}\\
\hline
\bsmall Subgroup&\bsmall Order&\bsmall Index\\
\hline
$\scriptstyle 3$& $\scriptstyle 3$&$\scriptstyle 1$\\
\hline
$\scriptstyle 2$& $\scriptstyle 2$&$\scriptstyle 3$\\
\hline
\end{tabular}
\begin{tabular}{|c|c|c|}
\hline \multicolumn{3}{|c|}{\bsmall $\scriptstyle C_4$}\\
\hline
\bsmall Subgroup&\bsmall Order&\bsmall Index\\
\hline
$\scriptstyle 4$& $\scriptstyle 4$&$\scriptstyle 1$\\
\hline
$\scriptstyle 2$& $\scriptstyle 2$&$\scriptstyle 2$\\
\hline
$\scriptstyle 1$& $\scriptstyle 1$&$\scriptstyle 4$\\
\hline
\end{tabular}$

$\begin{tabular}{|c|c|c|}
\hline \multicolumn{3}{|c|}{\bsmall $\scriptstyle C_6$}\\
\hline
\bsmall Subgroup&\bsmall Order&\bsmall Index\\
\hline
$\scriptstyle 6$& $\scriptstyle 6$&$\scriptstyle 1$\\
\hline
$\scriptstyle 3$& $\scriptstyle 3$&$\scriptstyle 2$\\
\hline
$\scriptstyle 2$& $\scriptstyle 2$&$\scriptstyle 3$\\
\hline
$\scriptstyle 1$& $\scriptstyle 1$&$\scriptstyle 6$\\
\hline
\end{tabular}
\begin{tabular}{|c|c|c|}
\hline \multicolumn{3}{|c|}{\bsmall $\scriptstyle C_{2h}$}\\
\hline
\bsmall Subgroup&\bsmall Order&\bsmall Index\\
\hline
$\scriptstyle 2/m$& $\scriptstyle 4$&$\scriptstyle 1$\\
\hline
$\scriptstyle 2$& $\scriptstyle 2$&$\scriptstyle 2$\\
\hline
$\scriptstyle m$& $\scriptstyle 2$&$\scriptstyle 2$\\
\hline
$\scriptstyle -1$& $\scriptstyle 2$&$\scriptstyle 2$\\
\hline
$\scriptstyle 1$& $\scriptstyle 1$&$\scriptstyle 4$\\
\hline
\end{tabular}
\begin{tabular}{|c|c|c|}
\hline \multicolumn{3}{|c|}{\bsmall $\scriptstyle C_{3h}$}\\
\hline
\bsmall Subgroup&\bsmall Order&\bsmall Index\\
\hline
$\scriptstyle -6$& $\scriptstyle 6$&$\scriptstyle 1$\\
\hline
$\scriptstyle 3$& $\scriptstyle 3$&$\scriptstyle 2$\\
\hline
$\scriptstyle m$& $\scriptstyle 2$&$\scriptstyle 3$\\
\hline
$\scriptstyle 1$& $\scriptstyle 1$&$\scriptstyle 6$\\
\hline
\end{tabular}$

$\begin{tabular}{|c|c|c|}
\hline \multicolumn{3}{|c|}{\bsmall $\scriptstyle C_{4h}$}\\
\hline
\bsmall Subgroup&\bsmall Order&\bsmall Index\\
\hline
$\scriptstyle 4/m$& $\scriptstyle 8$&$\scriptstyle 1$\\
\hline
$\scriptstyle 4$& $\scriptstyle 4$&$\scriptstyle 2$\\
\hline
$\scriptstyle -4$& $\scriptstyle 4$&$\scriptstyle 2$\\
\hline
$\scriptstyle 2/m$& $\scriptstyle 4$&$\scriptstyle 2$\\
\hline
$\scriptstyle 2$& $\scriptstyle 2$&$\scriptstyle 4$\\
\hline
$\scriptstyle m$& $\scriptstyle 2$&$\scriptstyle 4$\\
\hline
$\scriptstyle -1$& $\scriptstyle 2$&$\scriptstyle 4$\\
\hline
$\scriptstyle 1$& $\scriptstyle 1$&$\scriptstyle 8$\\
\hline
\end{tabular}
\begin{tabular}{|c|c|c|}
\hline \multicolumn{3}{|c|}{\bsmall $\scriptstyle C_{6h}$}\\
\hline
\bsmall Subgroup&\bsmall Order&\bsmall Index\\
\hline
$\scriptstyle 6/m$& $\scriptstyle 12$&$\scriptstyle 1$\\
\hline
$\scriptstyle -6$& $\scriptstyle 6$&$\scriptstyle 2$\\
\hline
$\scriptstyle 6$& $\scriptstyle 6$&$\scriptstyle 2$\\
\hline
$\scriptstyle -3$& $\scriptstyle 6$&$\scriptstyle 2$\\
\hline
$\scriptstyle 3$& $\scriptstyle 3$&$\scriptstyle 4$\\
\hline
$\scriptstyle 2/m$& $\scriptstyle 4$&$\scriptstyle 3$\\
\hline
$\scriptstyle m$& $\scriptstyle 2$&$\scriptstyle 6$\\
\hline
$\scriptstyle -1$& $\scriptstyle 2$&$\scriptstyle 6$\\
\hline
$\scriptstyle 1$& $\scriptstyle 1$&$\scriptstyle 12$\\
\hline
\end{tabular}
\begin{tabular}{|c|c|c|}
\hline \multicolumn{3}{|c|}{\bsmall $\scriptstyle C_{2v}$}\\
\hline
\bsmall Subgroup&\bsmall Order&\bsmall Index\\
\hline
$\scriptstyle mm2$& $\scriptstyle 4$&$\scriptstyle 1$\\
\hline
$\scriptstyle 2$& $\scriptstyle 2$&$\scriptstyle 2$\\
\hline
$\scriptstyle m$& $\scriptstyle 2$&$\scriptstyle 2$\\
\hline
$\scriptstyle 1$& $\scriptstyle 1$&$\scriptstyle 4$\\
\hline
\end{tabular}$

$\begin{tabular}{|c|c|c|}
\hline \multicolumn{3}{|c|}{\bsmall $\scriptstyle C_{3v}$}\\
\hline
\bsmall Subgroup&\bsmall Order&\bsmall Index\\
\hline
$\scriptstyle 3m$& $\scriptstyle 6$&$\scriptstyle 1$\\
\hline
$\scriptstyle m$& $\scriptstyle 2$&$\scriptstyle 3$\\
\hline
$\scriptstyle 1$& $\scriptstyle 1$&$\scriptstyle 6$\\
\hline
\end{tabular}
\begin{tabular}{|c|c|c|}
\hline \multicolumn{3}{|c|}{\bsmall $\scriptstyle C_{4v}$}\\
\hline
\bsmall Subgroup&\bsmall Order&\bsmall Index\\
\hline
$\scriptstyle 4mm$& $\scriptstyle 8$&$\scriptstyle 1$\\
\hline
$\scriptstyle 4$& $\scriptstyle 4$&$\scriptstyle 2$\\
\hline
$\scriptstyle mm2$& $\scriptstyle 4$&$\scriptstyle 2$\\
\hline
$\scriptstyle 2$& $\scriptstyle 2$&$\scriptstyle 4$\\
\hline
$\scriptstyle m$& $\scriptstyle 2$&$\scriptstyle 4$\\
\hline
$\scriptstyle 1$& $\scriptstyle 1$&$\scriptstyle 8$\\
\hline
\end{tabular}
\begin{tabular}{|c|c|c|}
\hline \multicolumn{3}{|c|}{\bsmall $\scriptstyle C_{6v}$}\\
\hline
\bsmall Subgroup&\bsmall Order&\bsmall Index\\
\hline
$\scriptstyle 6mm$& $\scriptstyle 12$&$\scriptstyle 1$\\
\hline
$\scriptstyle 6$& $\scriptstyle 6$&$\scriptstyle 2$\\
\hline
$\scriptstyle 3m$& $\scriptstyle 6$&$\scriptstyle 2$\\
\hline
$\scriptstyle mm2$& $\scriptstyle 4$&$\scriptstyle 3$\\
\hline
$\scriptstyle 2$& $\scriptstyle 2$&$\scriptstyle 6$\\
\hline
$\scriptstyle m$& $\scriptstyle 2$&$\scriptstyle 6$\\
\hline
$\scriptstyle 1$& $\scriptstyle 1$&$\scriptstyle 12$\\
\hline
\end{tabular}$

$\begin{tabular}{|c|c|c|}
\hline \multicolumn{3}{|c|}{\bsmall $\scriptstyle D_2$}\\
\hline
\bsmall Subgroup&\bsmall Order&\bsmall Index\\
\hline
$\scriptstyle 222$& $\scriptstyle 4$&$\scriptstyle 1$\\
\hline
$\scriptstyle 2$& $\scriptstyle 2$&$\scriptstyle 2$\\
\hline
$\scriptstyle 1$& $\scriptstyle 1$&$\scriptstyle 4$\\
\hline
\end{tabular}
\begin{tabular}{|c|c|c|}
\hline \multicolumn{3}{|c|}{\bsmall $\scriptstyle D_3$}\\
\hline
\bsmall Subgroup&\bsmall Order&\bsmall Index\\
\hline
$\scriptstyle 32$& $\scriptstyle 6$&$\scriptstyle 1$\\
\hline
$\scriptstyle 3$& $\scriptstyle 3$&$\scriptstyle 2$\\
\hline
$\scriptstyle 2$& $\scriptstyle 2$&$\scriptstyle 3$\\
\hline
$\scriptstyle 1$& $\scriptstyle 1$&$\scriptstyle 6$\\
\hline
\end{tabular}
\begin{tabular}{|c|c|c|}
\hline \multicolumn{3}{|c|}{\bsmall $\scriptstyle D_4$}\\
\hline
\bsmall Subgroup&\bsmall Order&\bsmall Index\\
\hline
$\scriptstyle 422$& $\scriptstyle 8$&$\scriptstyle 1$\\
\hline
$\scriptstyle 4$& $\scriptstyle 4$&$\scriptstyle 2$\\
\hline
$\scriptstyle 222$& $\scriptstyle 4$&$\scriptstyle 2$\\
\hline
$\scriptstyle 2$& $\scriptstyle 2$&$\scriptstyle 4$\\
\hline
$\scriptstyle 1$& $\scriptstyle 1$&$\scriptstyle 8$\\
\hline
\end{tabular}$

$\begin{tabular}{|c|c|c|}
\hline \multicolumn{3}{|c|}{\bsmall $\scriptstyle D_6$}\\
\hline
\bsmall Subgroup&\bsmall Order&\bsmall Index\\
\hline
$\scriptstyle 622$& $\scriptstyle 12$&$\scriptstyle 1$\\
\hline
$\scriptstyle 6$& $\scriptstyle 6$&$\scriptstyle 2$\\
\hline
$\scriptstyle 32$& $\scriptstyle 6$&$\scriptstyle 1$\\
\hline
$\scriptstyle 3$& $\scriptstyle 3$&$\scriptstyle 4$\\
\hline
$\scriptstyle 2$& $\scriptstyle 2$&$\scriptstyle 6$\\
\hline
$\scriptstyle m$& $\scriptstyle 2$&$\scriptstyle 6$\\
\hline
$\scriptstyle 1$& $\scriptstyle 1$&$\scriptstyle 12$\\
\hline
\end{tabular}
\begin{tabular}{|c|c|c|}
\hline \multicolumn{3}{|c|}{\bsmall $\scriptstyle D_{2h}$}\\
\hline
\bsmall Subgroup&\bsmall Order&\bsmall Index\\
\hline
$\scriptstyle mmm$& $\scriptstyle 8$&$\scriptstyle 1$\\
\hline
$\scriptstyle mm2$& $\scriptstyle 4$&$\scriptstyle 2$\\
\hline
$\scriptstyle 222$& $\scriptstyle 4$&$\scriptstyle 2$\\
\hline
$\scriptstyle 2/m$& $\scriptstyle 4$&$\scriptstyle 2$\\
\hline
$\scriptstyle 2$& $\scriptstyle 2$&$\scriptstyle 4$\\
\hline
$\scriptstyle m$& $\scriptstyle 2$&$\scriptstyle 4$\\
\hline
$\scriptstyle -1$& $\scriptstyle 2$&$\scriptstyle 4$\\
\hline
$\scriptstyle 1$& $\scriptstyle 1$&$\scriptstyle 8$\\
\hline
\end{tabular}
\begin{tabular}{|c|c|c|}
\hline \multicolumn{3}{|c|}{\bsmall $\scriptstyle D_{3h}$}\\
\hline
\bsmall Subgroup&\bsmall Order&\bsmall Index\\
\hline
$\scriptstyle -6$& $\scriptstyle 6$&$\scriptstyle 1$\\
\hline
$\scriptstyle 3$& $\scriptstyle 3$&$\scriptstyle 2$\\
\hline
$\scriptstyle m$& $\scriptstyle 2$&$\scriptstyle 3$\\
\hline
$\scriptstyle 1$& $\scriptstyle 1$&$\scriptstyle 6$\\
\hline
\end{tabular}$

$\begin{tabular}{|c|c|c|}
\hline \multicolumn{3}{|c|}{\bsmall $\scriptstyle D_{4h}$}\\
\hline
\bsmall Subgroup&\bsmall Order&\bsmall Index\\
\hline
$\scriptstyle 4/mmm$& $\scriptstyle 16$&$\scriptstyle 1$\\
\hline
$\scriptstyle -42m$& $\scriptstyle 8$&$\scriptstyle 2$\\
\hline
$\scriptstyle 422$& $\scriptstyle 8$&$\scriptstyle 2$\\
\hline
$\scriptstyle 4/m$& $\scriptstyle 8$&$\scriptstyle 2$\\
\hline
$\scriptstyle 4$& $\scriptstyle 4$&$\scriptstyle 4$\\
\hline
$\scriptstyle -4$& $\scriptstyle 4$&$\scriptstyle 4$\\
\hline
$\scriptstyle mmm$& $\scriptstyle 8$&$\scriptstyle 2$\\
\hline
$\scriptstyle mm2$& $\scriptstyle 4$&$\scriptstyle 4$\\
\hline
$\scriptstyle 222$& $\scriptstyle 4$&$\scriptstyle 4$\\
\hline
$\scriptstyle 2/m$& $\scriptstyle 4$&$\scriptstyle 4$\\
\hline
$\scriptstyle m$& $\scriptstyle 2$&$\scriptstyle 8$\\
\hline
$\scriptstyle 2$& $\scriptstyle 2$&$\scriptstyle 8$\\
\hline
$\scriptstyle -1$& $\scriptstyle 2$&$\scriptstyle 8$\\
\hline
$\scriptstyle 1$& $\scriptstyle 1$&$\scriptstyle 16$\\
\hline
\end{tabular}
\begin{tabular}{|c|c|c|}
\hline \multicolumn{3}{|c|}{\bsmall $\scriptstyle D_{6h}$}\\
\hline
\bsmall Subgroup&\bsmall Order&\bsmall Index\\
\hline
$\scriptstyle 6/mm$& $\scriptstyle 24$&$\scriptstyle 1$\\
\hline
$\scriptstyle -62m$& $\scriptstyle 12$&$\scriptstyle 2$\\
\hline
$\scriptstyle 6mm$& $\scriptstyle 12$&$\scriptstyle 2$\\
\hline
$\scriptstyle 622$& $\scriptstyle 12$&$\scriptstyle 2$\\
\hline
$\scriptstyle 6/m$& $\scriptstyle 12$&$\scriptstyle 2$\\
\hline
$\scriptstyle -3m$& $\scriptstyle 12$&$\scriptstyle 2$\\
\hline
$\scriptstyle -6$& $\scriptstyle 6$&$\scriptstyle 4$\\
\hline
$\scriptstyle 6$& $\scriptstyle 6$&$\scriptstyle 4$\\
\hline
$\scriptstyle 3m$& $\scriptstyle 6$&$\scriptstyle 4$\\
\hline
$\scriptstyle 32$& $\scriptstyle 6$&$\scriptstyle 4$\\
\hline
$\scriptstyle -3$& $\scriptstyle 6$&$\scriptstyle 4$\\
\hline
$\scriptstyle 6/m$& $\scriptstyle 12$&$\scriptstyle 2$\\
\hline
$\scriptstyle 3$& $\scriptstyle 3$&$\scriptstyle 8$\\
\hline
$\scriptstyle mmm$& $\scriptstyle 8$&$\scriptstyle 3$\\
\hline
$\scriptstyle mm2$& $\scriptstyle 4$&$\scriptstyle 6$\\
\hline
$\scriptstyle 222$& $\scriptstyle 4$&$\scriptstyle 6$\\
\hline
$\scriptstyle 2/m$& $\scriptstyle 4$&$\scriptstyle 6$\\
\hline
$\scriptstyle 2$& $\scriptstyle 2$&$\scriptstyle 12$\\
\hline
$\scriptstyle m$& $\scriptstyle 2$&$\scriptstyle 12$\\
\hline
$\scriptstyle -1$& $\scriptstyle 2$&$\scriptstyle 12$\\
\hline
$\scriptstyle 1$& $\scriptstyle 1$&$\scriptstyle 24$\\
\hline
\end{tabular}
\begin{tabular}{|c|c|c|}
\hline \multicolumn{3}{|c|}{\bsmall $\scriptstyle D_{2d}$}\\
\hline
\bsmall Subgroup&\bsmall Order&\bsmall Index\\
\hline
$\scriptstyle -42m$& $\scriptstyle 8$&$\scriptstyle 1$\\
\hline
$\scriptstyle -4$& $\scriptstyle 4$&$\scriptstyle 2$\\
\hline
$\scriptstyle mm2$& $\scriptstyle 4$&$\scriptstyle 2$\\
\hline
$\scriptstyle 222$& $\scriptstyle 4$&$\scriptstyle 2$\\
\hline
$\scriptstyle 2$& $\scriptstyle 2$&$\scriptstyle 4$\\
\hline
$\scriptstyle m$& $\scriptstyle 2$&$\scriptstyle 4$\\
\hline
$\scriptstyle 1$& $\scriptstyle 1$&$\scriptstyle 8$\\
\hline
\end{tabular}$

$\begin{tabular}{|c|c|c|}
\hline \multicolumn{3}{|c|}{\bsmall $\scriptstyle D_{3d}$}\\
\hline
\bsmall Subgroup&\bsmall Order&\bsmall Index\\
\hline
$\scriptstyle -3m$& $\scriptstyle 12$&$\scriptstyle 1$\\
\hline
$\scriptstyle 3m$& $\scriptstyle 6$&$\scriptstyle 2$\\
\hline
$\scriptstyle 32$& $\scriptstyle 6$&$\scriptstyle 2$\\
\hline
$\scriptstyle -3$& $\scriptstyle 6$&$\scriptstyle 2$\\
\hline
$\scriptstyle 3$& $\scriptstyle 3$&$\scriptstyle 4$\\
\hline
$\scriptstyle 2/m$& $\scriptstyle 4$&$\scriptstyle 3$\\
\hline
$\scriptstyle m$& $\scriptstyle 2$&$\scriptstyle 6$\\
\hline
$\scriptstyle -1$& $\scriptstyle 2$&$\scriptstyle 6$\\
\hline
$\scriptstyle 1$& $\scriptstyle 1$&$\scriptstyle 12$\\
\hline
\end{tabular}
\begin{tabular}{|c|c|c|}
\hline \multicolumn{3}{|c|}{\bsmall $\scriptstyle S_{4}$}\\
\hline
\bsmall Subgroup&\bsmall Order&\bsmall Index\\
\hline
$\scriptstyle -4$& $\scriptstyle 4$&$\scriptstyle 1$\\
\hline
$\scriptstyle 2$& $\scriptstyle 2$&$\scriptstyle 2$\\
\hline
$\scriptstyle 1$& $\scriptstyle 1$&$\scriptstyle 4$\\
\hline
\end{tabular}
\begin{tabular}{|c|c|c|}
\hline \multicolumn{3}{|c|}{\bsmall $\scriptstyle S_{6}=C_{3i}$}\\
\hline
\bsmall Subgroup&\bsmall Order&\bsmall Index\\
\hline
$\scriptstyle -3$& $\scriptstyle 6$&$\scriptstyle 1$\\
\hline
$\scriptstyle 3$& $\scriptstyle 3$&$\scriptstyle 2$\\
\hline
$\scriptstyle -1$& $\scriptstyle 2$&$\scriptstyle 3$\\
\hline
$\scriptstyle 1$& $\scriptstyle 1$&$\scriptstyle 6$\\
\hline
\end{tabular}$

$\begin{tabular}{|c|c|c|}
\hline \multicolumn{3}{|c|}{\bsmall $\scriptstyle T$}\\
\hline
\bsmall Subgroup&\bsmall Order&\bsmall Index\\
\hline
$\scriptstyle 23$& $\scriptstyle 12$&$\scriptstyle 1$\\
\hline
$\scriptstyle 3$& $\scriptstyle 3$&$\scriptstyle 4$\\
\hline
$\scriptstyle 222$& $\scriptstyle 4$&$\scriptstyle 3$\\
\hline
$\scriptstyle 2$& $\scriptstyle 2$&$\scriptstyle 6$\\
\hline
$\scriptstyle 1$& $\scriptstyle 1$&$\scriptstyle 12$\\
\hline
\end{tabular}
\begin{tabular}{|c|c|c|}
\hline \multicolumn{3}{|c|}{\bsmall $\scriptstyle T_h$}\\
\hline
\bsmall Subgroup&\bsmall Order&\bsmall Index\\
\hline
$\scriptstyle m-3$& $\scriptstyle 24$&$\scriptstyle 1$\\
\hline
$\scriptstyle 23$& $\scriptstyle 12$&$\scriptstyle 2$\\
\hline
$\scriptstyle -3$& $\scriptstyle 6$&$\scriptstyle 4$\\
\hline
$\scriptstyle 3$& $\scriptstyle 3$&$\scriptstyle 8$\\
\hline
$\scriptstyle mmm$& $\scriptstyle 8$&$\scriptstyle 3$\\
\hline
$\scriptstyle mm2$& $\scriptstyle 4$&$\scriptstyle 6$\\
\hline
$\scriptstyle 222$& $\scriptstyle 4$&$\scriptstyle 6$\\
\hline
$\scriptstyle 2/m$& $\scriptstyle 4$&$\scriptstyle 6$\\
\hline
$\scriptstyle 2$& $\scriptstyle 2$&$\scriptstyle 12$\\
\hline
$\scriptstyle m$& $\scriptstyle 2$&$\scriptstyle 12$\\
\hline
$\scriptstyle -1$& $\scriptstyle 2$&$\scriptstyle 12$\\
\hline
$\scriptstyle 1$& $\scriptstyle 1$&$\scriptstyle 24$\\
\hline
\end{tabular}
\begin{tabular}{|c|c|c|}
\hline \multicolumn{3}{|c|}{\bsmall $\scriptstyle T_d$}\\
\hline
\bsmall Subgroup&\bsmall Order&\bsmall Index\\
\hline
$\scriptstyle -43m$& $\scriptstyle 24$&$\scriptstyle 1$\\
\hline
$\scriptstyle 23$& $\scriptstyle 12$&$\scriptstyle 2$\\
\hline
$\scriptstyle 3m$& $\scriptstyle 6$&$\scriptstyle 4$\\
\hline
$\scriptstyle 3$& $\scriptstyle 3$&$\scriptstyle 8$\\
\hline
$\scriptstyle -42m$& $\scriptstyle 8$&$\scriptstyle 3$\\
\hline
$\scriptstyle -4$& $\scriptstyle 4$&$\scriptstyle 6$\\
\hline
$\scriptstyle mm2$& $\scriptstyle 4$&$\scriptstyle 6$\\
\hline
$\scriptstyle 222$& $\scriptstyle 4$&$\scriptstyle 6$\\
\hline
$\scriptstyle 2$& $\scriptstyle 2$&$\scriptstyle 12$\\
\hline
$\scriptstyle m$& $\scriptstyle 2$&$\scriptstyle 12$\\
\hline
$\scriptstyle 1$& $\scriptstyle 1$&$\scriptstyle 24$\\
\hline
\end{tabular}$

$\begin{tabular}{|c|c|c|}
\hline \multicolumn{3}{|c|}{\bsmall $\scriptstyle O$}\\
\hline
\bsmall Subgroup&\bsmall Order&\bsmall Index\\
\hline
$\scriptstyle 432$& $\scriptstyle 24$&$\scriptstyle 1$\\
\hline
$\scriptstyle 23$& $\scriptstyle 12$&$\scriptstyle 2$\\
\hline
$\scriptstyle 32$& $\scriptstyle 6$&$\scriptstyle 4$\\
\hline
$\scriptstyle 422$& $\scriptstyle 8$&$\scriptstyle 3$\\
\hline
$\scriptstyle 4$& $\scriptstyle 4$&$\scriptstyle 6$\\
\hline
$\scriptstyle 3$& $\scriptstyle 3$&$\scriptstyle 8$\\
\hline
$\scriptstyle 222$& $\scriptstyle 4$&$\scriptstyle 6$\\
\hline
$\scriptstyle 2$& $\scriptstyle 2$&$\scriptstyle 12$\\
\hline
$\scriptstyle 1$& $\scriptstyle 1$&$\scriptstyle 24$\\
\hline
\end{tabular}
\begin{tabular}{|c|c|c|}
\hline \multicolumn{3}{|c|}{\bsmall $\scriptstyle O_h$}\\
\hline
\bsmall Subgroup&\bsmall Order&\bsmall Index\\
\hline
$\scriptstyle m-3m$& $\scriptstyle 48$&$\scriptstyle 1$\\
\hline
$\scriptstyle -43m$& $\scriptstyle 24$&$\scriptstyle 2$\\
\hline
$\scriptstyle 432$& $\scriptstyle 24$&$\scriptstyle 2$\\
\hline
$\scriptstyle m-3$& $\scriptstyle 24$&$\scriptstyle 2$\\
\hline
$\scriptstyle 23$& $\scriptstyle 12$&$\scriptstyle 4$\\
\hline
$\scriptstyle -3m$& $\scriptstyle 12$&$\scriptstyle 4$\\
\hline
$\scriptstyle 3m$& $\scriptstyle 6$&$\scriptstyle 8$\\
\hline
$\scriptstyle 32$& $\scriptstyle 6$&$\scriptstyle 8$\\
\hline
$\scriptstyle -3$& $\scriptstyle 6$&$\scriptstyle 8$\\
\hline
$\scriptstyle 3$& $\scriptstyle 3$&$\scriptstyle 16$\\
\hline
$\scriptstyle 4/mmm$& $\scriptstyle 16$&$\scriptstyle 3$\\
\hline
$\scriptstyle -42m$& $\scriptstyle 8$&$\scriptstyle 6$\\
\hline
$\scriptstyle 4mm$& $\scriptstyle 8$&$\scriptstyle 6$\\
\hline
$\scriptstyle 422$& $\scriptstyle 8$&$\scriptstyle 6$\\
\hline
$\scriptstyle 4/m$& $\scriptstyle 8$&$\scriptstyle 6$\\
\hline
$\scriptstyle -4$& $\scriptstyle 4$&$\scriptstyle 12$\\
\hline
$\scriptstyle 4$& $\scriptstyle 4$&$\scriptstyle 12$\\
\hline
$\scriptstyle mmm$& $\scriptstyle 8$&$\scriptstyle 6$\\
\hline
$\scriptstyle 222$& $\scriptstyle 4$&$\scriptstyle 12$\\
\hline
$\scriptstyle 2/m$& $\scriptstyle 4$&$\scriptstyle 12$\\
\hline
$\scriptstyle 2$& $\scriptstyle 2$&$\scriptstyle 24$\\
\hline
$\scriptstyle m$& $\scriptstyle 2$&$\scriptstyle 24$\\
\hline
$\scriptstyle -1$& $\scriptstyle 2$&$\scriptstyle 24$\\
\hline
$\scriptstyle 1$& $\scriptstyle 1$&$\scriptstyle 48$\\
\hline
\end{tabular}$
\vfill\eject

\section{\bf APPENDIX C - AMALGAMATED FREE PRODUCTS IN THE AFFINE CRYSTALLOGRAPHIC SPACE-GROUP TYPES {\boldmath$[G(3)]$}}

$$\begin{tabular}{|c|}
\hline \multicolumn{1}{|c|}{\bsmall Amalgamated free products in $\scriptstyle [G(3)]$}\\
\hline
$\scriptstyle \mathbb{Z}_2\bigstar_{e}\mathbb{Z}_2=<\mathbf{x},\mathbf{y}>$\\
\hline
$\scriptstyle \mathbb{Z}_4\bigstar_{\mathbb{Z}_2}\mathbb{Z}_4=<4,\bar 4>$\\
\hline
$\scriptstyle \mathbb{Z}_4\bigstar_{\mathbb{Z}_2}D_2=<2,4,\mathbf{x}>$\\
\hline
$\scriptstyle \mathbb{Z}_6\bigstar_{\mathbb{Z}_3}D_3=<3,6,\mathbf{x}>$\\
\hline
$\scriptstyle D_4\bigstar_{D_2}D_4=<4,\bar 4,\mathbf{x}>$\\
\hline
$\scriptstyle D_2\times\mathbb{Z}_2\bigstar_{D_2}D_4=<2,4,\bar 1,\mathbf{x}>$\\
\hline
$\scriptstyle D_6\bigstar_{D_3}D_6=<6,\bar 6,\mathbf{x}>$\\
\hline
$\scriptstyle D_3\times\mathbb{Z}_2\bigstar_{D_3}D_6=<2,6,\bar 1,\mathbf{x}>$\\
\hline \multicolumn{1}{l}{\rsmall $\scriptstyle
\mathbf{x}=$reflection over $\scriptstyle x$-axis;
$\scriptstyle \mathbf{y}=$reflection over $\scriptstyle y$-axis.}\\
\end{tabular}$$

\section{\bf APPENDIX D - THE SUBGROUPS OF THE AFFINE CRYSTALLOGRAPHIC SPACE-GROUP TYPES {\boldmath$[G(2)]$}}

$$\begin{tabular}{|c|c|}
\hline \multicolumn{2}{|c|}{\bsmall $\scriptstyle p2$}\\
\hline
\bsmall Subgroup&\bsmall Index\\
\hline
$\scriptstyle p1$&$\scriptstyle 2$\\
\hline
\end{tabular}\begin{tabular}{|c|c|}
\hline \multicolumn{2}{|c|}{\bsmall $\scriptstyle pm$}\\
\hline
\bsmall Subgroup&\bsmall Index\\
\hline
$\scriptstyle pg$&$\scriptstyle 2$\\
\hline
$\scriptstyle cm$&$\scriptstyle 2$\\
\hline
$\scriptstyle p1$&$\scriptstyle $\\
\hline
\end{tabular}\begin{tabular}{|c|c|}
\hline \multicolumn{2}{|c|}{\bsmall $\scriptstyle pg$}\\
\hline
\bsmall Subgroup&\bsmall Index\\
\hline
$\scriptstyle p1$&$\scriptstyle 2$\\
\hline
\end{tabular}$$

$$\begin{tabular}{|c|c|}
\hline \multicolumn{2}{|c|}{\bsmall $\scriptstyle cm$}\\
\hline
\bsmall Subgroup&\bsmall Index\\
\hline
$\scriptstyle pg$&$\scriptstyle 2$\\
\hline
$\scriptstyle pm$&$\scriptstyle 2$\\
\hline
$\scriptstyle p2$&$\scriptstyle $\\
\hline
$\scriptstyle p1$&$\scriptstyle $\\
\hline
\end{tabular}\begin{tabular}{|c|c|}
\hline \multicolumn{2}{|c|}{\bsmall $\scriptstyle pmm$}\\
\hline
\bsmall Subgroup&\bsmall Index\\
\hline
$\scriptstyle cmm$&$\scriptstyle 2$\\
\hline
$\scriptstyle pmg$&$\scriptstyle 2$\\
\hline
$\scriptstyle pgg$&$\scriptstyle $\\
\hline
$\scriptstyle pm$&$\scriptstyle $\\
\hline
$\scriptstyle cm$&$\scriptstyle $\\
\hline
$\scriptstyle p2$&$\scriptstyle $\\
\hline
$\scriptstyle pg$&$\scriptstyle $\\
\hline
$\scriptstyle p1$&$\scriptstyle $\\
\hline
\end{tabular}\begin{tabular}{|c|c|}
\hline \multicolumn{2}{|c|}{\bsmall $\scriptstyle pgg$}\\
\hline
\bsmall Subgroup&\bsmall Index\\
\hline
$\scriptstyle p2$&$\scriptstyle 2$\\
\hline
$\scriptstyle pg$&$\scriptstyle 2$\\
\hline
$\scriptstyle p1$&$\scriptstyle $\\
\hline
\end{tabular}$$

$$\begin{tabular}{|c|c|}
\hline \multicolumn{2}{|c|}{\bsmall $\scriptstyle p4$}\\
\hline
\bsmall Subgroup&\bsmall Index\\
\hline
$\scriptstyle p2$&$\scriptstyle 2$\\
\hline
$\scriptstyle p1$&$\scriptstyle $\\
\hline
\end{tabular}\begin{tabular}{|c|c|}
\hline \multicolumn{2}{|c|}{\bsmall $\scriptstyle cmm$}\\
\hline
\bsmall Subgroup&\bsmall Index\\
\hline
$\scriptstyle pmm$&$\scriptstyle 2$\\
\hline
$\scriptstyle pmg$&$\scriptstyle 2$\\
\hline
$\scriptstyle cm$&$\scriptstyle 4$\\
\hline
$\scriptstyle pgg$&$\scriptstyle $\\
\hline
$\scriptstyle pm$&$\scriptstyle 2$\\
\hline
$\scriptstyle p2$&$\scriptstyle $\\
\hline
$\scriptstyle pg$&$\scriptstyle $\\
\hline
$\scriptstyle p1$&$\scriptstyle $\\
\hline
\end{tabular}\begin{tabular}{|c|c|}
\hline \multicolumn{2}{|c|}{\bsmall $\scriptstyle p4m$}\\
\hline
\bsmall Subgroup&\bsmall Index\\
\hline
$\scriptstyle p4g$&$\scriptstyle 2$\\
\hline
$\scriptstyle pmm$&$\scriptstyle 2$\\
\hline
$\scriptstyle cmm$&$\scriptstyle 2$\\
\hline
$\scriptstyle pmg$&$\scriptstyle 4$\\
\hline
$\scriptstyle p4$&$\scriptstyle 2$\\
\hline
$\scriptstyle pgg$&$\scriptstyle 4$\\
\hline
$\scriptstyle pm$&$\scriptstyle 4$\\
\hline
$\scriptstyle cm$&$\scriptstyle 4$\\
\hline
$\scriptstyle p2$&$\scriptstyle 4$\\
\hline
$\scriptstyle pg$&$\scriptstyle 8$\\
\hline
$\scriptstyle p1$&$\scriptstyle 8$\\
\hline
\end{tabular}$$

$$\begin{tabular}{|c|c|}
\hline \multicolumn{2}{|c|}{\bsmall $\scriptstyle p4g$}\\
\hline
\bsmall Subgroup&\bsmall Index\\
\hline
$\scriptstyle pmm$&$\scriptstyle 4$\\
\hline
$\scriptstyle cmm$&$\scriptstyle 2$\\
\hline
$\scriptstyle pmg$&$\scriptstyle $\\
\hline
$\scriptstyle p4$&$\scriptstyle 2$\\
\hline
$\scriptstyle cm$&$\scriptstyle $\\
\hline
$\scriptstyle pgg$&$\scriptstyle $\\
\hline
$\scriptstyle pm$&$\scriptstyle $\\
\hline
$\scriptstyle p2$&$\scriptstyle $\\
\hline
$\scriptstyle pg$&$\scriptstyle $\\
\hline
$\scriptstyle p1$&$\scriptstyle $\\
\hline
\end{tabular}\begin{tabular}{|c|c|}
\hline \multicolumn{2}{|c|}{\bsmall $\scriptstyle p6m$}\\
\hline
\bsmall Subgroup&\bsmall Index\\
\hline
$\scriptstyle pmm$&$\scriptstyle 6$\\
\hline
$\scriptstyle cmm$&$\scriptstyle 3$\\
\hline
$\scriptstyle pmg$&$\scriptstyle 6$\\
\hline
$\scriptstyle pgg$&$\scriptstyle 6$\\
\hline
$\scriptstyle p3m1$&$\scriptstyle 2$\\
\hline
$\scriptstyle p31m$&$\scriptstyle 2$\\
\hline
$\scriptstyle pm$&$\scriptstyle 12$\\
\hline
$\scriptstyle cm$&$\scriptstyle 6$\\
\hline
$\scriptstyle pg$&$\scriptstyle 12$\\
\hline
$\scriptstyle p6$&$\scriptstyle 2$\\
\hline
$\scriptstyle p3$&$\scriptstyle 4$\\
\hline
$\scriptstyle p2$&$\scriptstyle 6$\\
\hline
$\scriptstyle p1$&$\scriptstyle 12$\\
\hline
\end{tabular}\begin{tabular}{|c|c|}
\hline \multicolumn{2}{|c|}{\bsmall $\scriptstyle p3$}\\
\hline
\bsmall Subgroup&\bsmall Index\\
\hline
$\scriptstyle p1$&$\scriptstyle 3$\\
\hline
\end{tabular}$$

$$\begin{tabular}{|c|c|}
\hline \multicolumn{2}{|c|}{\bsmall $\scriptstyle p3m1$}\\
\hline
\bsmall Subgroup&\bsmall Index\\
\hline
$\scriptstyle pm$&$\scriptstyle 6$\\
\hline
$\scriptstyle cm$&$\scriptstyle 3$\\
\hline
$\scriptstyle p3$&$\scriptstyle 2$\\
\hline
$\scriptstyle p31m$&$\scriptstyle 3$\\
\hline
$\scriptstyle pg$&$\scriptstyle $\\
\hline
$\scriptstyle p1$&$\scriptstyle $\\
\hline
\end{tabular}\begin{tabular}{|c|c|}
\hline \multicolumn{2}{|c|}{\bsmall $\scriptstyle p31m$}\\
\hline
\bsmall Subgroup&\bsmall Index\\
\hline
$\scriptstyle p3m1$&$\scriptstyle 3$\\
\hline
$\scriptstyle p3$&$\scriptstyle 2$\\
\hline
$\scriptstyle p3$&$\scriptstyle 2$\\
\hline
$\scriptstyle pm$&$\scriptstyle 6$\\
\hline
$\scriptstyle cm$&$\scriptstyle 3$\\
\hline
$\scriptstyle pg$&$\scriptstyle $\\
\hline
$\scriptstyle p1$&$\scriptstyle $\\
\hline
\end{tabular}\begin{tabular}{|c|c|}
\hline \multicolumn{2}{|c|}{\bsmall $\scriptstyle p6$}\\
\hline
\bsmall Subgroup&\bsmall Index\\
\hline
$\scriptstyle p3$&$\scriptstyle 2$\\
\hline
$\scriptstyle p2$&$\scriptstyle 3$\\
\hline
$\scriptstyle p1$&$\scriptstyle $\\
\hline
\end{tabular}$$

\end{document}